\documentclass[11pt]{amsart}
\usepackage{amsxtra}
\usepackage{amssymb}
\theoremstyle{plain}
\newcommand{\cleqn}{\setcounter{equation}{0}}
\newcommand{\clth}{\setcounter{theorem}{0}}
\newcommand {\sectionnew}[1]{\section{#1}\cleqn\clth}
\newcommand{\nn}{\hfill\nonumber}
\newtheorem{theorem}{Theorem}[section]
\newtheorem{lemma}[theorem]{Lemma}
\newtheorem{definition-theorem}[theorem]{Definition-Theorem}
\newtheorem{proposition}[theorem]{Proposition}
\newtheorem{corollary}[theorem]{Corollary}
\newtheorem{definition}[theorem]{Definition}
\newtheorem{example}[theorem]{Example}
\newtheorem{remark}[theorem]{Remark}
\newtheorem{conjecture}[theorem]{Conjecture}
\newtheorem{notation}[theorem]{Notation}
\newcommand \bth[1] { \begin{theorem}\label{t#1} }
\newcommand \ble[1] { \begin{lemma}\label{l#1} }

\newcommand \bpr[1] { \begin{proposition}\label{p#1} }
\newcommand \bco[1] { \begin{corollary}\label{c#1} }
\newcommand \bde[1] { \begin{definition}\label{d#1}\rm }
\newcommand \bex[1] { \begin{example}\label{e#1}\rm }
\newcommand \bre[1] { \begin{remark}\label{r#1}\rm }
\newcommand \bcj[1] { \begin{conjecture}\label{j#1}\rm }

\newcommand \bnota[1] { \begin{notation}\label{n#1}\rm }
\renewcommand {\eth} { \end{theorem} }
\newcommand {\ele} { \end{lemma} }

\newcommand {\epr} { \end{proposition} }
\newcommand {\eco} { \end{corollary} }
\newcommand {\ede} { \end{definition} }
\newcommand {\eex} { \end{example} }
\newcommand {\ere} { \end{remark} }
\newcommand {\ecj} { \end{conjecture} }

\newcommand {\enota} { \end{notation} }
\newcommand \thref[1]{Theorem \ref{t#1}}
\newcommand \leref[1]{Lemma \ref{l#1}}
\newcommand \prref[1]{Proposition \ref{p#1}}
\newcommand \coref[1]{Corollary \ref{c#1}}

\newcommand \deref[1]{Definition \ref{d#1}}

\newcommand \reref[1]{Remark \ref{r#1}}
\newcommand \lb[1]{\label{#1}}

\def \Cset {{\mathbb C}}
\def \KK {{\mathbb K}}
\def \Tset {{\mathbb T}}
\def \Zset {{\mathbb Z}}
\def \Nset {{\mathbb N}}
\def \Qset {{\mathbb Q}}

\def \B  {{\mathcal{B}}}
\def \II {{\mathcal{I}}}

\def \UU {{\mathcal{U}}}
\def \RR {{\mathcal{R}}}
\def \SS {{\mathcal{S}}}
\def \SLl {{\mathcal{SL}}}
\def \Kl {{\mathcal{K}}}
\def \LL {{\mathcal{L}}}
 
\def \HH {{\mathcal{H}}}
\def \De {\Delta}   
\def \de {\delta}
\def \al {\alpha}
\def \be {\beta}
\def \la {\lambda}
\def \La {\Lambda}
\def \Om {\Omega} 
\def \om {\omega}
\def \ga {\gamma}
\def \de {\delta}
\def \Ga {\Gamma}
\def \Sig {\Sigma}
\def \sig {\sigma}

\def \ep {\epsilon}
\def \sig{\sigma}

\def \mt  {\mapsto}

\def \hra {\hookrightarrow}
\def \lha {\leftharpoonup}
\def \rha {\rightharpoonup}

\def \ci  {\circ}

\def \rcor {\rangle}
\def \lcor {\langle}
\def \o  {\otimes}
\def \ol {\overline}
\def \wt {\widetilde}
\def \wh {\widehat}

\def \id { {\mathrm{id}} }
\def \Stab { {\mathrm{Stab}} }
\def \sign { {\mathrm{sign}} }

\def \Lie { {\mathrm{Lie \,}} }

\def \g  {\mathfrak{g}}   

\def \h  {\mathfrak{h}}

\def \n  {\mathfrak{n}}

\def \b  {\mathfrak{b}}

\def \bfw {\bf{\mathrm{w}}}

\DeclareMathOperator \Span { {\mathrm{Span}} }
\DeclareMathOperator \Prim { {\mathrm{Prim}} }
\DeclareMathOperator \Fract { {\mathrm{Fract}} }
\DeclareMathOperator \Sympl { {\mathrm{Sympl}} }
\DeclareMathOperator \Aut { {\mathrm{Aut}} }
\DeclareMathOperator \red { {\mathrm{red}} }
\DeclareMathOperator \GKdim { {\mathrm{GK \, dim}}}

\DeclareMathOperator \htt { {\mathrm{ht}} }
\DeclareMathOperator \Diag { {\mathrm{Diag}} }
\DeclareMathOperator \Supp { {\mathrm{Supp}} }

\DeclareMathOperator \Ker { {\mathrm{Ker}} }

\DeclareMathOperator \gr { {\mathrm{gr}} }

\DeclareMathOperator \Irr { {\mathrm{Irr}} }

\DeclareMathOperator \congto  {{ \stackrel{\cong}{\to} }}

\newcommand \Spec { {\mathrm{Spec}} }
\newcommand \Max  { {\mathrm{Max}} }
\begin{document}
\title[On the spectra of quantum groups]
{On the spectra of quantum groups}
\author[Milen Yakimov]{Milen Yakimov}
\address{
Department of Mathematics \\
Louisiana State University \\
Baton Rouge, LA 70803 \\
U.S.A 
}
\email{yakimov@math.lsu.edu}
\date{}
\begin{abstract}
Joseph and Hodges--Levasseur (in the A case) described the 
spectra of all quantum function algebras $R_q[G]$ on simple 
algebraic groups in terms of the centers of certain localizations 
of quotients of $R_q[G]$ by torus invariant prime ideals,
or equivalently in terms of orbits of finite groups.
These centers were only known up to finite extensions.
We determine the centers explicitly under the general 
conditions that the deformation parameter is not a root 
of unity and without any restriction on the characteristic 
of the ground field. From it we deduce a more explicit description 
of all prime ideals of $R_q[G]$ than the previously known 
ones and an explicit parametrization of $\Spec R_q[G]$.
We combine the latter with a result of Kogan 
and Zelevinsky to obtain in the complex case a torus 
equivariant Dixmier type map from the symplectic foliation of the 
group $G$ to the primitive spectrum of $R_q[G]$. Furthermore, 
under the general assumptions on the ground field and deformation
parameter, we prove a theorem for separation of variables for 
the De Concini--Kac--Procesi algebras $\UU^w_\pm$, and 
classify the sets of their homogeneous normal elements and
primitive elements. We apply those results to obtain explicit 
formulas for the prime and especially the primitive ideals of 
$\UU^w_\pm$ lying in the Goodearl--Letzter stratum over the 
$\{0\}$-ideal. This is in turn used to prove that all Joseph's 
localizations of quotients of $R_q[G]$ by torus invariant prime
ideals are free modules over their subalgebras generated by 
Joseph's normal elements.
From it we derive a classification of the maximal 
spectrum of  $R_q[G]$ and use it to resolve a question of Goodearl  
and Zhang, showing that all maximal ideals of $R_q[G]$ have 
finite codimension. We prove that $R_q[G]$ possesses a stronger 
property than that of the classical catenarity: all maximal chains in 
$\Spec R_q[G]$ have the same length equal to $\GKdim R_q[G]= \dim G$,
i.e. $R_q[G]$ satisfies the first chain condition for 
prime ideals in the terminology of Nagata.
\end{abstract}
\maketitle
\tableofcontents
\sectionnew{Introduction}
\lb{intro}
In the last 20 years the area of quantum groups attracted a lot of 
attention in ring theory, since it supplied large families of concrete 
algebras on which general techniques can be developed and tested. One of the most 
studied family is the one of the quantum function algebras $R_q[G]$ 
on simple groups. In works from the early 90's Joseph \cite{J1, J} and 
Hodges--Levasseur \cite{HL0,HL} (and jointly with Toro \cite{HLT}) 
made fundamental contributions to the problem of determining 
their spectra by describing the spectra set theoretically and 
the topology of each stratum in a finite stratification of $\Spec R_q[G]$.
Despite the fact that a lot of research has been done since then
(see \cite{BG,GK,GL0,GLen0,GL,GZ,G,J2,LLR,LeSt,LWZ}),  
many ring theoretic questions for the algebras $R_q[G]$ remain open. 

In this paper we describe solutions of ring theoretic problems for $R_q[G]$ and 
the related De Concini--Kac--Procesi algebras $\UU^w_\pm$,
which range from the older question of determining explicitly the 
Joseph strata of $\Spec R_q[G]$ and setting up a torus equivariant Dixmier type map 
between the symplectic foliation of $G$ and $\Prim R_q[G]$,
to newer ones such as the question raised 
by Goodearl and Zhang \cite{GZ} on whether all maximal ideals of $R_q[G]$ 
have finite codimension and the question of classifying $\Max R_q[G]$.

In order to describe in concrete terms the problems addressed in this paper,
we introduce some notation on quantum groups.
We start with the quantized universal enveloping algebra $\UU_q(\g)$ 
of a simple (finite dimensional) Lie algebra $\g$ of rank $r$. Throughout the paper our assumption is 
that the base field $\KK$ is arbitrary and the deformation 
parameter $q \in \KK^*$ is not a root of unity (except the
small part on the Dixmier map and Poisson geometry where 
$\KK= \Cset$ and $q \in \Cset^*$ is not a root of unity). We 
do not use specialization at any point and thus do not need the 
ground field to have characteristic 0 and $q$ to be transcendental 
over $\Qset$. The quantized algebra of functions $R_q[G]$ on the 
``connected, simply connected'' group $G$ is the Hopf subalgebra 
of the restricted dual $(\UU_q(\g))^\ci$ consisting of the matrix coefficients 
of all finite 
dimensional type 1 representations of $\UU_q(\g)$. (Here $G$ 
is just a symbol, since all arguments are carried out over an arbitrary 
base field $\KK$. The only restriction is that $\KK$ is not finite, 
since $q \in \KK^*$ should not be a root of unity.)

All finite dimensional type 1 $\UU_q(\g)$-modules are completely reducible.
The irreducible ones are parametrized by the dominant integral weights 
$P^+$ of $\g$ and have $q$-weight space decompositions 
$V(\la) = \oplus_{\mu \in P} V(\la)_\mu$, where $P$ is the weight lattice of 
$\g$. The matrix coefficient of $\xi \in V(\la)^*$ and $v\in V(\la)$ will be 
denoted by $c^\la_{\xi, v} \in R_q[G]$. The algebra $R_q[G]$ is 
$P \times P$-graded by 
\begin{equation}
\label{grR}
R_q[G]_{\nu, \mu} = \Span \{ c^\la_{\xi, v} \mid \la \in P^+, 
\xi \in (V(\la)^*)_{\nu}, V(\la)_\mu \}  
\end{equation}
and has two distinguished subalgebras $R^\pm$ 
which are spanned by matrix coefficients $c^\la_{\xi, v}$ for 
$v \in V(\la)_\la$ and $v \in V(\la)_{ w_0 \la}$, respectively,
where $w_0$ denotes the longest element of the Weyl group $W$ 
of $\g$. Joseph proved \cite{J1} that $R_q[G] = R^+ R^- = R^- R^+$. 
For $w \in W$ one defines the Demazure modules 
\[
V^+_w(\la) = \UU_+ V(\la)_{w \la} \subseteq V(\la) \; \; 
\mbox{and} \; \; 
V^-_w(\la) = \UU_- V( - w_0 \la)_{- w \la} \subseteq V( - w_0 \la),
\]
where $\UU_\pm$ 
are the subalgebras of $\UU_q(\g)$ generated by the positive 
and negative Cartan generators.
To $w_\pm \in W$, Joseph 
associated certain ideals $I^\pm_{w_\pm}$ of $R^\pm$ by considering the span 
of those matrix coefficients for which the vector $\xi$ is 
orthogonal to $V^\pm_{w_\pm}(\la)$, see \S 2.3 for details. Those 
are combined into the following ideals of $R_q[G]$: 
\[
I_{\bfw} = I^+_{w_+} R^- + R^+ I^-_{w_-}, \quad 
{\bfw}= (w_+, w_-),  
\]
which are a key building block in Joseph's analysis 
of $R_q[G]$. The other part is a set of normal elements 
of the quotients $R_q[G]/I_{\bfw}$. Up to an 
appropriate normalization, for $\la \in P^+$ one defines
\begin{equation}
\label{celem}
c^+_{w_+, \la} = c^\la_{- w_+ \la, \la} \in R^+ \; \; 
\mbox{and} \; \; 
c^-_{w_-, \la} = c^{- w_0 \la}_{w_- \la, -\la} \in R^-,
\end{equation}
where vectors in $\UU_q(\g)$-modules are substituted with 
weights for the (one dimensional) weight spaces 
to which they belong, see \S \ref{2.4}.
The multiplicative subsets of $R^\pm$ generated by the first and second 
kind of elements will be denoted by $E^\pm_{w_\pm}$,
and their product by $E_{\bfw} \subset R_q[G]$. (As it is customary 
we will denote by the same symbols the images of elements
of $R_q[G]$ in its quotients.) The sets 
$E_{\bfw} \subset R_q[G]/I_{\bfw}$ consist of regular normal 
elements, thus one can localize 
$R_{\bfw}:= (R_q[G]/I_{\bfw})[E_{\bfw}^{-1}]$. Joseph \cite{J1, J} and 
Hodges--Levasseur (in the A case \cite{HL} and jointly with Toro
in the multiparameter case \cite{HLT}) proved that one can break down
\begin{equation}
\label{strata}
\Spec R_q[G] = \bigsqcup_{ {\bfw} \in W \times W } \Spec_{\bfw} 
R_q[G], 
\end{equation}
in such a way that $\Spec_{\bfw} R_q[G]$ is homeomorphic
to the spectrum of the center $Z_{\bfw} : = Z(R_{\bfw})$ via:
\begin{equation}
\label{iot}
J^0 \in \Spec Z_{\bfw} \mt \iota_{\bfw}(J^0) \in \Spec R_q[G] \; \; 
\mbox{so that} \; \; 
\iota_{\bfw}(J^0)/I_{\bfw} = (R_{\bfw} J^0) \cap (R_q[G]/I_{\bfw}) ,
\end{equation}
see \thref{J-thm} for details.
Joseph's original formulation of the parametrization of $\Spec_{\bfw} R_q[G]$ 
is in slightly different 
terms using an action of $\Zset_2^{\times r}$, see
\cite[Theorem 8.11]{J1}, \cite[Theorem 10.3.4]{J}. 
In this form the parametrization of $\Spec_{\bfw} R_q[G]$
is stated in Hodges--Levasseur--Toro 
\cite[Theorem 4.15]{HLT}.

Joseph and Hodges--Levasseur--Toro
proved that for all ${\bfw}=(w_+, w_-) \in W \times W$, 
$Z_{\bfw}$ is a Laurent polynomial ring over $\KK$ 
of dimension $\dim \ker (w_+ - w_-)$, 
that $Z_{\bfw}$ contains a particular 
Laurent polynomial ring, and that it is a free module over it 
of rank at most $2^r$. They also proved that the closure of 
each stratum is a union of strata given in terms of the inverse 
Bruhat order on $W \times W$, but the nature of the gluing of the 
strata $\Spec_{\bfw} R_q [G]$ inside $\Spec R_q[G]$ with the Zariski topology 
is unknown. 

One needs to know the explicit structure 
of $Z_{\bfw}$ to begin studying the Zariski 
topology of the space $\Spec R_q[G]$ in the sense of the  
interaction between the different strata in \eqref{strata}.
This is also needed for the construction of an equivariant Dixmier 
type map from the symplectic foliation of the underlying Poisson 
Lie group $G$ to $\Prim R_q[G]$ (when $\KK= \Cset$), since 
one cannot determine the stabilizers of the primitive ideals of $R_q[G]$
with respect to the natural torus actions 
(see \eqref{Tract}--\eqref{Tract0} below) from Joseph's theorem.
The centers $Z_{\bfw}$ are explicitly known only for 
$\g = {\mathfrak{sl}}_2$ and $\g= {\mathfrak{sl}}_3$
due to Hodges--Levasseur \cite{HL0} and 
Goodearl--Lenagan \cite{GLen0}. The first problem that we solve in 
this paper is the one of the explicit description of $Z_{\bfw}$ 
in full generality. At the time of the writing of \cite{J1,HL} a similar 
problem existed on the symplectic side. Hodges and Levasseur
proved \cite{HL0} that in the complex case the double Bruhat cells 
$G^{\bfw} \subset G$ are torus orbits 
of symplectic leaves and that certain intersections 
are finite unions of at most $2^r$ symplectic leaves, but no results 
were available on their connected components. Nine years later this
problem was solved by Kogan and Zelevinsky \cite{KZ}
by combinatorial methods using the theory of generalized 
minors. The problem of determining $Z_{\bfw}$ is the ring theoretic 
counterpart of the problem which they solved.

For an element ${\bfw} = (w_+, w_-)$ denote by $\SS({\bfw})$
the subset of $\{1, \ldots, r\}$ consisting of all 
simple reflections $s_i$ which appear either in a reduced 
expression of $w_+$ or $w_-$ (one thinks of it as of the 
support of ${\bfw})$. Denote by $\II({\bfw})$ its complement 
in $\{1, \ldots, r\}$, which is the set 
of all fundamental weights $\om_i$ fixed by $w_+$ and $w_-$. 
We denote 
\[
P_{\SS({\bfw})} = \oplus_{i \in \SS({\bfw})} \Zset \om_i \; \; 
\mbox{and} \; \; 
\wt{\LL}_{\red}( {\bfw} ) = \ker (w_+ - w_-) \cap P_{\SS({\bfw})}. 
\] 
It is easy to see that $\wt{\LL}_{\red}({\bfw})$ is a lattice 
of rank $\dim \ker (w_+ - w_-) - |\II({\bfw})|$. 
One extends the definition \eqref{celem} 
to $\la \in P$ to obtain elements in the localizations $R_{\bfw}$.
We have:
\bth{1} Assume that $\KK$ is an arbitrary base field
and $q \in \KK^*$ is not a root of unity. For any of the 
quantum function algebras $R_q[G]$ and  
${\bfw} =(w_+, w_-) \in W \times W$ the center $Z_{\bfw}$ 
of Joseph's localization $R_{\bfw}$ coincides with the Laurent polynomial 
algebra over $\KK$ of dimension $\dim \ker (w_+ - w_-)$ with generators
\[
\{ c^+_{w_+, \om_i} \mid i \in \II({\bfw}) \} \sqcup
\{ c^+_{w_+, \la^{(j)}} (c^-_{w_-, \la^{(j)}})^{-1} \}_{j=1}^k,
\]
where $k = \dim \ker (w_+ - w_-) - | \II({\bfw})|$ and 
$\la^{(1)}, \ldots, \la^{(k)}$ is a basis of $\wt{\LL}_{\red}({\bfw})$.
\eth

Joseph proved that for all $\la \in \wt{\LL}({\bfw}) = \Ker(w_+- w_-) \cap P$, 
\begin{equation}
\label{Jel}
c^+_{w_+, \la} (c^+_{w_+, \la} )^{-1} \in Z_{\bfw}.
\end{equation}
The center $Z_{\bfw}$ is a free module over it of rank equal 
to $2^{|\II({\bfw})|}$. The difficulty in the proof of the above 
theorem is not the guess of the exact form of the center $Z_{\bfw}$
(which can be interpreted as taking square roots of some of the elements
\eqref{Jel}), but the proof of the fact that $Z_{\bfw}$ does not
contain additional elements.

The proof of \thref{1} appears in Section \ref{cent}. We make use 
of a model of the algebra $R_{\bfw}$ due to Joseph. We expect that 
this model will play an important role in the future study of 
$\Spec R_q[G]$. Joseph \cite{J1,J} defines the algebras 
$S^\pm_{w_\pm}$ as the $0$ components of the localizations 
$(R^\pm/I^\pm_{w_\pm})[(E^\pm_{w_\pm})^{-1}]$ with respect to the 
second grading in \eqref{grR} (induced to the localization).
He then defines an algebra $S_{\bfw}$ which is a kind 
of bicrossed product of $S^+_{w_+}$ and $S^-_{w_+}$, and proves 
that $R_{\bfw}$ is isomorphic to a smash product 
of a Laurent polynomial ring and a localization of $S_{\bfw}$ 
by a set of normal elements. We refer the reader to 
\cite[\S 9.1-9.3 and \S 10.3]{J} and \S \ref{3.4} for details.

De Concini, Kac and Procesi defined a family of subalgebras 
$\UU^w_\pm \subseteq \UU_\pm \subset \UU_q(\g)$, which are 
parametrized by the elements $w$ of the Weyl group of $\g$. They
can be viewed as deformations of the universal enveloping 
algebras $\UU(\n_\pm \cap w (\n_\mp))$, where $\n_\pm$ are the 
nilradicals of a pair of opposite Borel subalgebras. The 
algebras $\UU^w_\pm$ 
are defined in terms of the Lusztig's root vectors of $\UU_q(\g)$.
Our approach to determining the center $Z_{\bfw}$ of $R_{\bfw}$ 
is to make use of a family of (anti)isomorphisms 
$\varphi^\pm_{w_\pm} \colon S^\pm_{w_\pm} \to \UU^{w_\pm}_\mp$, 
which was a main ingredient in our work \cite{Y} on the torus invariant 
spectra of $\UU^{w_\pm}_\mp$ (see \thref{isom} below).
With the help of these (anti)isomorphisms
we study $R_{\bfw}$ using on one side the De Concini--Kac--Procesi 
PBW bases of $\UU^{w_\pm}_\mp$ and the Levendorskii--Soibelman
straightening rule. On the other side we use techniques
from quantum function algebras which produce good supplies 
of normal elements and exploit the $R$-matrix type commutation 
relations inside the algebras $S^\pm_{w_\pm}$ and between 
them in the ``bicrossed product'' $S_{\bfw}$. The weight lattice 
$P$ of $\g$ acts in a natural way on the algebras $S^\pm_{w_\pm}$ 
and $S_{\bfw}$ by algebra automorphisms. Using the above mentioned 
techniques we investigate the set of homogeneous $P$-normal elements of 
the algebras $S^\pm_{w_\pm}$ and $S_{\bfw}$, and obtain from that
\thref{1} for the center of $R_{\bfw}$.

There is a natural action of 
the torus $\Tset^r \times \Tset^r = (\KK^*)^{\times 2 r}$ on 
$R_q[G]$ by algebra automorphisms, which quantizes the 
left and right regular actions of 
the maximal torus of $G$ on the coordinate ring of $G$, 
see \eqref{Tract2}. Joseph proved that, 
if the base field $\KK$ is algebraically closed, then the stratum
of primitive ideals 
$\Prim_{\bfw} R_q[G] = \Prim R_q[G] \cap \Spec_{\bfw} R_q[G]$ is a 
single $\Tset^r$-orbits with respect to the action of each 
component. Each stratum $\Prim_{\bfw} R_q[G]$ is preserved by 
$\Tset^r \times \Tset^r$. In Section \ref{Dixmier} we apply 
the results of \thref{1} to determine the exact structure of 
$\Prim_{\bfw} R_q[G]$ as a $\Tset^r \times \Tset^r$-homogeneous 
space. 

Now let us restrict ourselves to the case when $\KK = \Cset$ and $q \in \Cset^*$ 
is not a root of unity. The connected, simply connected complex 
algebraic group $G$ 
corresponding to $\g$ is equipped with the so called standard 
Poisson structure $\pi_G$. It follows from the Kogan--Zelevinsky
results \cite{KZ} that all symplectic leaves of $\pi_G$ are locally 
closed in the Zariski topology. We denote by $\Sympl (G, \pi_G)$ 
the symplectic foliation of $\pi_G$ with the topology 
induced from the Zariski topology of $G$. Joseph established 
that there is a bijection from $\Sympl (G, \pi_G)$ to $\Prim R_q[G]$,
that maps symplectic leaves in a double Bruhat cell $G^{\bfw}$ 
to the primitive ideals in $\Prim_{\bfw} R_q[G]$, thus settling 
a conjecture of Hodges and Levasseur \cite{HL0}. 
In Section 
\ref{Dixmier} we make this picture 
$\Tset^r \times \Tset^r$-equivariant. The torus
$\Tset^r \times \Tset^r$ acts on $(G, \pi_G)$ by Poisson 
maps. (One identifies $\Tset^r$ with a maximal torus of $G$ 
and uses the left and right regular actions.) This induces a 
$\Tset^r \times \Tset^r$-action on $\Sympl(G, \pi_G)$.   
Combining \thref{1} with the results of 
Kogan and Zelevinsky \cite{KZ}, we prove:
\bth{2} For each connected, simply connected, complex simple algebraic 
group 
$G$ and $q \in \Cset^*$ which is not a root of unity, there exists a
$\Tset^r \times \Tset^r$-equivariant map 
\begin{equation}
\label{DG}
D_G \colon \Sympl (G, \pi_G) \to \Prim R_q[G].
\end{equation}
\eth
This map is explicitly constructed in \S \ref{4.5}. The Hodges--Levasseur 
idea \cite{HL0} for an orbit method for $R_q[G]$ now can be formulated 
more concretely by conjecturing that \eqref{DG} is a homeomorphism.

For the rest of this introduction we return to the general assumptions 
on $\KK$ and $q$.
In order to be able to compare prime ideals in different strata
$\Spec_{\bfw} R_q[G]$ for ${\bfw} = (w_+, w_-) \in W \times W$, 
one needs to investigate the maps $\iota_{\bfw}$ 
from \eqref{iot}. For this one needs to know the structure of the 
algebras $R_{\bfw}$ as modules over their subalgebras generated by 
Joseph's normal sets $E_{\bfw}^{ \pm 1}$. Because Joseph's model for 
$R_{\bfw}$ is based on the algebras $S^\pm_{w_\pm}$, one first needs 
to investigate the module structure of $S^\pm_{w_\pm}$ over their 
subalgebras $N^\pm_{w_\pm}$, generated by the following 
normal elements of $S^\pm_{w_\pm}$:
\begin{equation}
\label{delem}
d^\pm_{w_\pm, \la} = (c^\pm_{w_\pm, \la})^{-1} c^\pm_{1, \la}, \quad  
\la \in P^+_{\SS(w_\pm)},
\end{equation}
cf. \eqref{IS}. The classical theorems for separation of variables 
of Kostant \cite{K} and Joseph--Letzter \cite{JL} prove that 
$\UU(\g)$ and $\UU_q(\g)$ are free modules over their centers and establish 
a number of properties of the related decomposition.
In Section \ref{free1} we prove a result for separation of 
variables for the algebras $S^\pm_{w_\pm}$. The difference here is 
that the algebras $S^\pm_{w_\pm}$ behave like universal enveloping 
algebras of nilpotent Lie algebras since they are (anti)isomorphic 
to the De Concini--Kac--Procesi algebras $\UU^{w_\pm}_{\mp}$
which are deformations of $\UU(\n_\mp \cap w_\pm(\n_\pm))$.
Because of this, generally they have small centers compared to their 
localizations by the multiplicative set of scalar multiples of the 
elements \eqref{delem}, see \reref{cent} for a precise comparison.
Because of this and for the ultimate purposes of classifying $\Max R_q[G]$, 
we consider the structure of $S^\pm_{w_\pm}$ as modules over their
subalgebras $N^\pm_{w_\pm}$. 
\bth{3} 
For an arbitrary base field $\KK$, $q \in \KK^*$ not a root of unity 
and $w_\pm \in W$, the algebras $S^\pm_{w_\pm}$ are free left and 
right modules over their subalgebras $N^\pm_{w_\pm}$ (generated 
by the normal elements \eqref{delem}), in which $N^\pm_{w_\pm}$ 
are direct summands.
\eth
Moreover we construct explicit bases of $S^\pm_{w_\pm}$ as 
$N^\pm_{w_\pm}$-modules using the PBW bases of 
$\UU^{w_\pm}_{\mp}$, see \thref{freeS2} for details. 
\thref{3} and some detailed analysis of the normal elements of 
$S^\pm_{w_\pm}$ lead us to the following classification 
result, proved in Section \ref{normal}.
\bth{4} For an arbitrary base field $\KK$, 
$q \in \KK^*$ not a root of unity and 
$w \in W$, the nonzero homogeneous normal 
elements of $S^\pm_w$ are precisely the nonzero
scalar multiples of the elements \eqref{delem}.
All such elements are distinct and even more 
the elements \eqref{delem} are linearly 
independent. All normal elements of the algebras $S^\pm_w$ 
are linear combinations of homogeneous normal elements.
\eth
The De Concini--Kac--Procesi algebras $\UU^w_\pm$ attracted 
a lot of attention from a ring theoretic perspective in recent 
years. The reason is that they contain as special cases various 
important families such as the algebras of quantum matrices.
Moreover, they are the largest known family of 
Cauchon--Goodearl--Letzter extensions which are a kind of iterated skew 
polynomial extensions for which both the Goodearl--Letzter  
stratification theory \cite{GL} and the Cauchon theory of deleted 
derivations \cite{Cau} work. A number of results were obtained for their torus 
invariant prime ideals. M\'eriaux and Cauchon classified them as 
a set \cite{MC}. The author described all inclusions between them 
and obtained an explicit formula for each torus invariant 
prime ideal \cite{Y}, using results of Gorelik \cite{G}. However 
there are no results explicitly describing prime ideals which are 
not torus invariant, except in some very special cases. 
Theorems \ref{t3} and \ref{t4} have many applications 
in this direction. 
Firstly, \thref{4} classifies the sets of all
normal elements of $\UU^\pm_{w_\pm}$, see 
Section \ref{normal} for details. Only the case 
when $w$ equals the longest element of the Weyl group 
of $\g$ was previously known due
to Caldero \cite{Ca}. Secondly, using Theorems \ref{t3} and \ref{t4} in 
Section \ref{normal} we obtain 
an explicit formula 
for the prime and especially the primitive ideals in the 
Goodearl--Letzter stratum of $\Spec \UU^{w_\pm}_\pm$
over the $\{0\}$-ideal. A result of Launois, Lenagan and Rigal 
\cite[Theorem 3.7]{LLR} implies that the algebras 
$\UU^{w_\pm}_\pm$ are noetherian unique factorization domains
(see \S \ref{5.2} for background). Therefore one is interested in 
knowing the sets of their prime elements. We classify all prime 
elements of the algebras $\UU^{w_\pm}_\pm$ in Theorems 
\ref{tnormal2} (ii) and \ref{tinhompr}. (For the latter result 
one needs to translate the results from the setting of $S^\mp_{w_\pm}$ 
to the one of $\UU^{w_\pm}_{\pm}$ via the (anti)isomorphisms 
of \thref{isom}, which is straightforward and is not 
stated separately.) As a corollary we obtain 
explicit formulas for all height one prime ideals of 
$\UU^{w_\pm}_\pm$.  

We return to the problem of describing the structure of the 
localizations $R_{\bfw}$ as modules
over their subalgebras generated by the Joseph sets of normal 
elements $(E_{\bfw})^{ \pm 1}$ for ${\bfw} =  (w_+, w_-) \in W \times W$.
We denote the latter subalgebras 
of $R_{\bfw}$ by $L_{\bfw}$. One cannot deduce the module structure 
of $R_{\bfw}$ over $L_{\bfw}$ immediately from the one of $S^\pm_{w_\pm}$
over $N^\pm_{w_\pm}$ (recall \thref{3}), because the former 
is not a tensor product of the latter in Joseph's model 
for $R_{\bfw}$. To overcome this difficulty, in Section \ref{Module},
we define nontrivial
$Q \times Q$ filtrations of the algebras $S_{\bfw}$, where
$Q$ is the root lattice of $\g$. (Note that
the $P \times P$-grading \eqref{grR} of $R_q[G]$ only induces 
a $Q$-grading on $S_{\bfw}$.) The associated graded of the new
$Q \times Q$-filtration of $S_{\bfw}$ breaks down in a certain 
way the ``bicrossed product'' of $S^+_{w_+}$ and $S^-_{w_-}$,
and then one can apply \thref{3}. In Section \ref{Module} we prove:

\bth{5} For all base fields $\KK$, $q \in \KK^*$ not a root of unity 
and ${\bfw} \in W \times W$, the algebras $L_{\bfw}$ are quantum tori
of dimension $r + |\SS({\bfw})|$. Moreover, Joseph's localizations 
$R_{\bfw}$ are free (left and right) modules over $L_{\bfw}$ in which 
$L_{\bfw}$ are direct summands.
\eth

In addition, in \thref{f1} we construct an explicit $L_{\bfw}$-basis 
of each $R_{\bfw}$. In Theorems \ref{t3} and \ref{t5} we do not 
obtain further representation theoretic properties of the bases 
for those free modules as Kostant and Joseph--Letzter did in the 
cases of $\UU(\g)$ and $\UU_q(\g)$. We think that this is 
an important problem which deserves future study.

For the purposes of the study of catenarity and homological properties
of $R_q[G]$ and its Hopf algebra quotients, Goodearl and Zhang \cite{GZ}
raised the question whether all maximal ideals of $R_q[G]$
have finite codimension. So far this was known in only two cases,
$\g = {\mathfrak{sl}}_2$ and ${\mathfrak{sl}}_3$ due 
to Hodges--Levasseur \cite{HL0} and Goodearl--Lenagan \cite{GLen0}.
In Section \ref{Max} we classify $\Max R_q[G]$ and settle affirmatively
the question of Goodearl and Zhang in full generality.

\bth{6} For all base fields $\KK$, $q \in \KK^*$ not a root 
of unity and a simple Lie algebra $\g$ the maximal spectrum 
of $R_q[G]$ is
\[
\Max R_q[G] = \Prim_{(1,1)} R_q[G].
\] 
The maximal spectrum of $R_q[G]$ is homeomorphic 
to the maximal spectrum of an $r$ dimensional 
Laurent polynomial ring. All maximal ideals
of $R_q[G]$ have finite codimension.
\eth

In addition, \thref{mid} provides an explicit formula 
for all maximal ideals of $R_q[G]$, see also \coref{mid1} 
for the case when $\KK$ is algebraically closed. 
The difficult part of \thref{6} is to show that 
$\Max R_q[G] \subset \Spec R_q[G]$.
Our approach is to consider the projection 
$\pi_{\bfw} \colon R_{\bfw} \to L_{\bfw}$ along the 
direct complement from \thref{5}. We use the formula for 
primitive ideals $J \in \Prim_{\bfw} R_q[G]$ from Section 
\ref{Dixmier} to study the projection $\pi_{\bfw}(J)$. 
We compare it to $\pi_{\bfw}(I_{(1,1)})$, to deduce 
that for ${\bfw} \neq (1,1)$, 
$J+ I_{(1,1)} \neq R_q[G]$.

A ring $R$ satisfies {\em{the first chain condition 
for prime ideals}} if all maximal chains in $\Spec R$
have the same length equal to $\GKdim R$. 
This is a stronger property than catenarity.
It was introduced by Nagata \cite{N} in the commutative 
case. Combining \thref{6} and results of Goodearl and Zhang 
\cite{GZ}, in Section \ref{Homology} we prove:
\bth{7} For all base fields $\KK$, $q \in \KK^*$ not a root 
of unity, and Hopf ideals $I$ of $R_q[G]$, the quotient
$R_q[G]/I$ satisfies the first chain condition 
for prime ideals and Tauvel's height formula holds.     
\eth

Gei\ss, Leclerc and Schr\"oer \cite{GLS} have recently proved 
that $\UU^w_+$ are quantum cluster algebras 
for symmetric Kac--Moody algebras $\g$ 
and base field $\KK = \Qset(q)$.
It will be very interesting if cluster algebra 
and ring theoretic techniques can be combined 
in the study of $\Spec \UU^w_-$ and $\Spec R_q[G]$.

{\bf Acknowledgements.} We are grateful to Ken Goodearl for 
many helpful discussions, email correspondence and advice on 
the available literature, and to Tom Lenagan for communicating 
to us the proof of the second part of \prref{factor}.
We also thank the referee whose comments helped us
to improve the exposition.

The research of the author was supported by NSF grants DMS-0701107 
and DMS-1001632.
\sectionnew{Previous results on spectra of quantum function algebras}
\lb{qalg}
\subsection{Quantized universal enveloping algebras}
\label{2.1}
In this section we collect background material on
quantum groups and some previous work on their spectra,
which will be used in the paper.

We fix a base field $\KK$ and 
$q \in \KK^* = \KK \backslash \{ 0 \}$
which is not a root of unity. Let $\g$ be a simple 
(finite dimensional) Lie algebra of rank $r$ with Cartan matrix
$(c_{ij})$.  Its quantized universal 
enveloping algebra $\UU_q(\g)$ over $\KK$ 
with deformation parameter $q$  
is a Hopf algebra over $\KK$ with generators
\[
X^\pm_i, K_i^{\pm 1}, \; i=1, \ldots, r
\]
and relations
\begin{gather*}
K_i^{-1} K_i = K_i K^{-1}_i = 1, \, K_i K_j = K_j K_i, 
\\
K_i X^\pm_j K^{-1}_i = q_i^{\pm c_{ij}} X^\pm_j,
\\
X^+_i X^-_j - X^-_j X^+_i = \de_{i,j} \frac{K_i - K^{-1}_i}
{q_i - q^{-1}_i},
\\
\sum_{k=0}^{1-c_{ij}} (-1)^k
\begin{bmatrix} 
1-c_{ij} \\ k
\end{bmatrix}_{q_{i}}
      (X^\pm_i)^k X^\pm_j (X^\pm_i)^{1-c_{ij}-k} = 0, \, i \neq j,
\end{gather*}
where $q_i = q^{d_i}$ and $\{d_i\}_{i=1}^r$ are
the positive relatively 
prime integers such that $(d_i c_{ij})$ is symmetric.  
The coproduct of $\UU_q(\g)$ is given by:
\begin{gather*}
\De(K_i)   = K_i \o K_i,
\\
\De(X^+_i) = X^+_i \o 1 + K_i \o X^+_i,
\\
\De(X^-_i) = X^-_i \o K_i^{-1} + 1 \o X^-_i.
\end{gather*}
Its antipode and counit are given by:
\begin{gather*}
S(K_i) = K^{-1}_i, \, 
S(X^+_i)= - K^{-1}_i X^+_i, \, 
S(X^-_i)= - X^-_i K_i,
\\
\ep(K_i), \, \ep(X^\pm_i)=0.
\end{gather*}
As usual $q$-integers, $q$-factorials, and
$q$-binomial coefficients are denoted by 
\[
[n]_q = \frac{q^n - q^{-n}}{q - q^{-1}}, \, 
[n]_q ! = [1]_q \ldots [n]_q, \, 
\begin{bmatrix}
n \\ m  
\end{bmatrix}_{q}
= \frac{[n]_q}{[m]_q [n-m]_q},
\]
$n, m \in \Nset$, $m \leq n$. We refer to \cite[Ch. 4]{Ja} 
for more details.

Denote by $\UU_\pm$ the subalgebras of $\UU_q(\g)$
generated by $\{X^\pm_i\}_{i=1}^r$. Let $H$ be the group 
generated by $\{K_i^{\pm 1}\}_{i=1}^r$, i.e. the group
of all group-like elements of $\UU_q(\g)$.
\subsection{Type 1 modules and braid group action}
\label{2.2}
The sets of simple roots, simple coroots, and 
fundamental weights of $\g$ will be denoted by $\{\al_i\}_{i=1}^r,$ 
$\{ \al_i\spcheck \}_{i=1}^r$,
and $\{\om_i\}_{i=1}^r,$ respectively.
Denote by $P$ and $P^+$ the sets of integral and dominant integral 
weights of $\g$. For $\la = \sum_i n_i \om_i \in P$, let
\begin{equation}
\label{supp}
\Supp \la = \{ i \in \{ 1, \ldots, r \}  \mid n_i \neq 0 \}.
\end{equation}
Denote the root lattice $Q = \sum_i \Zset \al_i$
of $\g$. Set $Q^+ = \sum_i \Nset \al_i$.  Let $Q\spcheck$ 
be the coroot lattice of $\g$. We will use 
the following standard partial order on $P$:
\begin{equation}
\label{po}
\mbox{for} \; \;  \la_1, \la_2 \in P, \; 
\la_1 < \la_2, \; \; 
\mbox{if and only if} \; \; 
\la_2 - \la_1 \in Q^+ \backslash \{0\}.
\end{equation}
Denote by $\De^+$
and $\De^-$ the sets of positive and negative roots of $\g$.

Let $\lcor.,. \rcor$ be the nondegenerate bilinear form on 
$\Span_\Qset \{ \al_1, \ldots, \al_r \}$
defined by 
\begin{equation}
\label{inner}
\lcor \al_i, \al_j \rcor = d_i c_{ij}.
\end{equation}
 
The $q$-weight spaces of an $H$-module $V$ are defined by
\[
V_\mu = \{ v \in V \mid K_i v = q^{ \lcor \mu, \al_i \rcor} v, \; 
\forall i = 1, \ldots, r \}, \; \mu \in P.
\]
A $\UU_q(\g)$-module is called a type 1 module if it is the 
sum of its $q$-weight spaces, see \cite[Ch. 5]{Ja} for details. 
The irreducible finite dimensional type 1 $\UU_q(\g)$-modules 
are parametrized by $P^+$, \cite[Theorem 5.10]{Ja}. 
Let $V(\la)$ denote the irreducible 
type $1$ $\UU_q(\g)$-module with highest weight $\la \in P^+$.
Let $M(\la)$ denote the Verma module of $\UU_q(\g)$ with highest 
weight $\la$ and highest weight vector $u_\la$. 
For an arbitrary base field $\KK$, and $q \in \KK^*$ which is not a 
root of unity, the Weyl character formula holds for $V(\la)$ 
and $V(\la)$ is given as a quotient of $M(\la)$ by the standard
formula from the classical case (see 
\cite[Corollary 7.7]{APW} and \cite[p. 126]{Ja}):
\begin{equation}
\label{Vla}
V(\la) \cong M(\la) / \Big( \sum_{i=1}^r 
\UU_q(\g) (X_i^-)^{ \lcor \la, \al_i\spcheck \rcor +1} u_\la \Big). 
\end{equation}
All duals of finite dimensional $\UU_q(\g)$-modules
will be considered as left modules using the antipode of 
$\UU_q(\g)$. The category of finite dimensional type 1 
$\UU_q(\g)$-modules is semisimple \cite[Theorem 5.17]{Ja}
(cf. also the remark on p. 85 of \cite{Ja}) and is
closed under taking tensor products and duals.  

Denote by $W$ and $\B_\g$ the Weyl and braid groups 
associated to $\g$. The simple 
reflections of $W$ corresponding to $\al_1, \ldots, \al_r$ 
will be denoted by $s_1, \ldots, s_r$. The 
corresponding generators of $\B_\g$ will be 
denoted by $T_1, \ldots, T_r$.
For a Weyl group element $w$, 
$l(w)$ will denote its length. The Bruhat order 
on $W$ will be denoted by $\leq$.

Lusztig defined actions of $\B_\g$ on all finite dimensional type 1 
modules and $\UU_q(\g)$. On a finite dimensional type 1 module $V$ 
the generators $T_1, \ldots, T_r$ of $\B_\g$ act by
(see \cite[\S 8.6]{Ja} and \cite[\S 5.2]{L} for details):
\begin{equation}
\label{braid}
T_i(v) = \sum_{l, m, n}
(-1)^m q_i^{m- ln} (X^+_i)^{(l)} (X^-_i)^{(m)} (X^+_i)^{(n)} v, 
\quad v \in V_\mu, \mu \in P,
\end{equation}
where the sum is over $l,m, n \in \Nset$ such that 
$-l+m-n = \lcor \mu, \al_i\spcheck \rcor$ and
\[
(X^\pm_i)^{(l)}= \frac{(X^\pm_i)^l}{[l]_{q_i}} \cdot
\]
The action of the braid group $\B_\g$ satisfies
\[
T_w V(\la)_\mu = V(\la)_{w(\mu)}, \quad \forall 
\la \in P^+, \mu \in P, w \in W.
\]
This implies that $\dim V(\la)_{w \la} = 1$ for 
$\la \in P^+, w \in W$.
The braid group $\B_\g$ acts on $\UU_q(\g)$ by
\begin{align*}
& T_i(X_i^+) = - X_i^- K_i, \; T_i(X_i^-) = - K_i^{-1} X_i^+, \;
T_i(K_j) = K_j K_i^{-c_{ij}}, \\
& T_i(X_j^+) = \sum_{k=0}^{- c_{ij}} (- q_i)^{-k} (X_i^+)^{(-c_{ij}-k)}
X_j^+ (X_i^+)^{(k)}, \; j \neq i, \\
& T_i(X_j^-) = \sum_{k=0}^{- c_{ij}} (- q_i)^k (X_i^-)^{(k)}
X_j^- (X_i^-)^{(-c_{ij}-k)}, \; j \neq i,
\end{align*}
see \cite[\S 8.14]{Ja} and \cite[\S 37.1]{L} for details.
The two actions on type 1 finite dimensional 
modules and $\UU_q(\g)$ are compatible:
\begin{equation}
\label{compat}
T_w(x . v) = (T_w x). (T_w v),
\end{equation}
for all $w \in W$, $x \in \UU_q(\g)$, $v \in V(\la)$,
see \cite[eq. 8.14(1)]{Ja}.
\subsection{$H$-prime ideals of Quantum Groups}
\label{2.3}

Denote by $R_q[G]$ the Hopf subalgebra of the restricted 
dual of $\UU_q(\g)$ spanned by all matrix coefficients 
of the modules $V(\la)$. It is a noetherian domain, see 
\cite[Lemma 9.1.9 (i) and Proposition 9.2.2]{J} and 
\cite[Corollary 5.6]{BG2}. For $\xi \in V(\la)^*$, 
$v \in V(\la)$ define
\[
c^\la_{\xi, v} \in R_q[G] \; \; \mbox{by} \; \;  
c^\la_{\xi, v} (x) = \xi (xv), \; \forall x \in 
\UU_q(\g).
\]
There are two canonical left and right actions of $\UU_q(\g)$ 
on $R_q[G]$ given by
\begin{equation}
\label{action}
x \rha c = \sum c_{(2)}(x)c_{(1)}, \; 
c \lha x = \sum c_{(1)}(x)c_{(2)}, \; 
x \in \UU_q(\g), c \in R_q[G]
\end{equation}
and a corresponding $P\times P$-grading on $R_q[G]$
\begin{equation}
\label{Rgrade}
R_q[G]_{\nu, \mu} = \{ c^\la_{\xi, v} \mid
\la \in P^+, \xi \in (V(\la)^*)_\nu, v \in V(\la)_\mu \}, 
\quad \nu, \mu \in P.
\end{equation}

Define the subalgebras of $R_q[G]$
\begin{align*}
R^+ &= 
\Span \{c^\la_{\xi, v} \mid \la \in P^+, 
v \in V(\la)_\la, \xi \in V(\la)^* \},
\\
R^- & = \Span \{ c^\la_{\xi, v} \mid \la \in P^+, 
v \in V(\la)_{w_0 \la}, \xi \in V(\la)^* \}, 
\end{align*}
where $w_0$ denotes the longest element of $W$.
Joseph proved \cite[Proposition 9.2.2]{J} that 
$R_q[G]= R^+ R^- = R^- R^+$, see \S 3.4 below
for more details.

Throughout the paper we fix highest weight vectors 
$v_\la \in V(\la)_\la$, $\la \in P^+$. Denote the 
corresponding lowest weight vectors 
$v_{-\la} = T_{w_0} v_{- w_0 \la} \in V(-w_0 \la)_{-\la}$.
For $\xi \in V(\la)^*$ and $\xi' \in V(- w_0 \la)^*$ denote
\begin{equation}
\label{special-c}
c^\la_{\xi, \la} := c^\la_{\xi, v_\la } \; \; 
\mbox{and} \; \; 
c^{- w_0 \la}_{\xi', - \la} := c^{-w_0 \la}_{\xi', v_{-\la} }.
\end{equation}
As vector spaces $R^+$ and $R^-$ can be identified 
with $\oplus_{\la \in P^+} V(\la)^*$ by
\[
\xi \in V(\la)^* \mt c^\la_{\xi, \la} \; \; \mbox{and} \; \; 
\xi' \in V(- w_0 \la)^* \mt c^{- w_0 \la}_{\xi', - \la},
\]
respectively.
Then the multiplication in $R^\pm$ can be identified 
with the Cartan multiplication rule (see \cite[\S 9.1.6]{J}) 
\begin{equation}
\label{Cartan}
V(\la_1)^* V(\la_2)^* \to V(\la_1 + \la_2)^*,
\quad \xi_1. \xi_2 := (\xi_1 \otimes \xi_2) |_{V(\la_1 + \la_2)}, 
\end{equation}
where
$\la_i \in P^+$, $\xi_i \in V(\la_i)^*$, $i=1,2$. 
We normalize the embeddings
$V(\la_1 + \la_2) \hra V(\la_1) \otimes_\KK V(\la_2)$ so 
that $v_{\la_1 + \la_2} \mt v_{\la_1} \otimes v_{\la_2}$.
Then
$T_{w_0} v_{\la_1 + \la_2} \mt 
T_{w_0} v_{\la_1} \otimes T_{w_0} v_{\la_2}$,
see \eqref{RwT} below. Thus 
$v_{ -\la_1 - \la_2} \mt v_{-\la_1} \otimes v_{-\la_2}$ 
under 
$V(- w_0(\la_1 + \la_2) ) \hra V(- w_0 \la_1) \otimes_\KK V(- w_0 \la_2)$

Recall that for all $w \in W$ the weight spaces 
$V(\la)_{w \la}$ are one dimensional, see \S \ref{2.2}.
Define the Demazure modules
\begin{equation}
\label{Demazure}
V_w^+(\la) = \UU_+ V(\la)_{w \la} \subseteq V(\la), \quad 
V_w^-(\la) = \UU_- V(- w_0 \la)_{- w \la} \subseteq V(- w_0 \la),
\end{equation}
for $\la \in P^+$, $w \in W$, and the canonical projections
\begin{equation}
\label{g}
g^+_{w_+} \colon V(\la)^* \to (V_{w_+}^+(\la))^*
\; \; \mbox{and} \; \; 
g^-_{w_-} \colon V(- w_0 \la)^* \to (V_{w_-}^-(\la))^*. 
\end{equation}

Following Joseph \cite{J1, J}
and Hodges--Levasseur \cite{HL0,HL}, define 
\begin{align}
\label{Qw+}
I_w^+ &= \Span \{ c^\la_{\xi, v} \mid \la \in P^+, v \in V(\la)_\la, 
\xi \in V_w^+(\la)^\perp \} \subset R^+, \\
\label{Qw-}
I_w^- &= \Span \{ c^{-w_0 \la}_{\xi, v} \mid \la \in P^+, 
v \in V(- w_0 \la)_{-\la}, \xi \in V_w^-(\la)^\perp \} \subset R^-.
\end{align}
For ${\bfw}=(w_+, w_-) \in W \times W$ define 
\begin{equation}
\label{Qbfw}
I_{\bfw} = I_{w_+}^+ R^- + R^+ I_{w_-}^- \subset R_q[G].
\end{equation}

\bth{Hprim} (Joseph, \cite[Proposition 8.9]{J}, 
\cite[Proposition 10.1.8, Proposition 10.3.5]{J1}) 

(i) For each $w \in W$, 
$I_w^\pm$ is an $H$-invariant completely prime ideal
of $R^\pm$ with respect to the left action \eqref{action} of $H$.
All $H$-primes of $R^\pm$ are of this form.

(ii) For each ${\bfw} \in W \times W$,
$I_{\bfw}$ is a an $H$-invariant completely prime ideal
of $R_q[G]$ with respect to the left action of $H$.
All $H$-primes of $R_q[G]$ are of this form.
\eth

In \cite{J1} \thref{Hprim} was stated for $\KK = \Cset$, 
$q \in \Cset^*$ not a root of unity and in \cite{J} 
\thref{Hprim} was stated for $\KK= k(q)$ 
for a field $k$ of characteristic $0$.
It is well known that Joseph's proof works for an 
arbitrary field $\KK$, $q \in \KK^*$ not a root of 1, 
see \cite[\S 3.3]{Ko} for a related discussion.

The ideals $I_{\bfw}$ are also $H$-stable 
with respect to the right action \eqref{action}.
The left and right invariance property of the 
ideals $I_{\bfw}$ with respect to $H$ can be 
formulated in terms of invariance with respect 
to a torus action. See \S \ref{2.6} for details.
\subsection{Sets of normal elements}
\label{2.4}
Recall that for $w \in W$ the weight spaces 
$V(\la)_{w \la}$ are one dimensional.
For $\la \in P^+$, $w \in W$, denote by 
$\xi^+_{w, \la} \in (V(\la)^*)_{- w \la}$ and 
$\xi^-_{w, \la} \in (V(-w_0 \la)^*)_{w \la}$ 
the vectors normalized by
\begin{equation}
\label{normaliza}
\lcor \xi^+_{w, \la}, T_w v_\la \rcor = 1 \; \; 
\mbox{and} \; \; 
\lcor \xi^-_{w, \la}, T_{w^{-1}}^{-1} v_{-\la} \rcor = 1.
\end{equation}
Define
\begin{equation}
\label{cd-def}
c^+_{w, \la} = c^{\la}_{\xi^+_{w,\la} , v_\la}, \quad
c^-_{w, \la} = c^{- w_0 \la}_{ \xi^-_{w, -\la}, v_{-\la}}
\end{equation}
in terms of the highest and lowest weight vectors 
$v_\la$ and $v_{-\la}$, fixed in \S \ref{2.3}.
These normalizations are chosen to
match the Kogan--Zelevinsky
normalizations \cite{KZ}. This will ensure a proper alignment 
of the semiclassical and quantum pictures in 
Section \ref{Dixmier}, see \reref{special}.

One has: 
\begin{equation}
\label{normalization}
c^+_{w, \la_1} c^+_{w, \la_2} = c^+_{w, \la_1 + \la_2} \; \; 
\mbox{and} \; \; 
c^-_{w, \la_1} c^-_{w, \la_2} = c^-_{w, \la_1 + \la_2}, \quad 
\forall \la_1, \la_2 \in P^+, w \in W.
\end{equation}
This follows from the equalities 
\begin{align}
\label{RwT}
T_w( v_{\la_1} \otimes v_{\la_2} ) &= T_w( v_{\la_1}) \otimes T_w( v_{\la_2} ),
\\
T_{w^{-1}}^{-1}( v_{- \la_1} \otimes v_{- \la_2} ) 
&= T_{w^{-1}}^{-1}( v_{- \la_1}) \otimes T_{w^{-1}}^{-1}( v_{- \la_2} ),
\label{RwT2}
\end{align} 
$\forall \la_1, \la_2 \in P^+, w \in W$.
Eq. \eqref{RwT2} is a consequence of eq. \eqref{RwT}:
if $w_0 = w^{-1} w'$, then for all $\la \in P^+$ one has
$T_{w^{-1}}^{-1} v_{- \la} = 
T_{w^{-1}}^{-1} T_{w_0} v_{- w_0 \la} = 
T_{w^{-1}}^{-1} T_{w^{-1}} T_{w'} v_{ - w_0 \la} = T_{w'} v_{ - w_0 \la}$. 
If $w = s_{i_1} \ldots s_{i_l}$ 
is a reduced expression of $w$, then $T_{s_{j+1} \ldots s_{i_l}} v_\la$ is a 
highest weight vector for the $\UU_{q_{i_j}}({\mathfrak{sl_2}})$-subalgebra 
of $\UU_q(\g)$ generated by $X^\pm_{i_j}$, $K^{\pm 1}_{i_j}$, 
$\forall j= 1, \ldots, l$. Because of 
this it is sufficient to verify \eqref{RwT} for $\g = {\mathfrak{sl_2}}$. 
For $\g = {\mathfrak{sl_2}}$ one computes:
\begin{align*}
& T_1( v_{n \om_1} \otimes v_{m \om_1}) = 
\frac{(-1)^{n +m} q^{n +m}}{[n+ m]_q !} (X^-_1)^{n + m} 
( v_{n \om_1} \otimes v_{m \om_1}) = 
\\
= & 
\frac{(-1)^{n + m} q^{n + m}}{[n + m]_q !}
\left( \sum_{0 \leq k_1 < \ldots < k_n < n+m} 
q^{ k_1 + \ldots + k_n - nm} (X_1^-)^n v_{n \om_1} \otimes (X_1^-)^m v_{m \om_1} 
\right)
\\
= & \frac{(-1)^{n + m} q^{n + m}}{[n + m]_q !} 
\begin{bmatrix}
n + m \\ m  
\end{bmatrix}_{q}
(X_1^-)^n v_{n \om_1} \otimes (X_1^-)^m v_{m \om_1} =
T_1( v_{n \om_1} ) \otimes T_1 (v_{m \om_1}).
\end{align*}
\subsection{Localizations of quotients of $R_q[G]$ by its $H$-primes}
\label{2.5}
The algebra $\UU_q(\g)$ is $Q$-graded by 
\begin{equation}
\label{Uq-grading}
\deg X^\pm_i = \pm \al_i, \; \; \deg K_i = 0, \quad
i=1, \ldots, r. 
\end{equation}
The homogeneous component of $\UU_q(\g)$ corresponding to $\ga \in Q$ 
will be denoted by $(\UU_q(\g))_\ga$. 

For $\ga \in Q^+$, $\ga \neq 0$ denote 
$m(\ga) = \dim (\UU_+)_\ga= \dim (\UU_-)_{-\ga}$
and fix a pair of dual bases 
$\{u_{\ga, i} \}_{i=1}^{m(\ga)}$ and
$\{u_{-\ga, i} \}_{i=1}^{m(\ga)}$
of $(\UU_+)_\ga$ and $(\UU_-)_{-\ga}$ 
with respect to the Rosso--Tanisaki form,
see \cite[Ch. 6]{Ja} for a discussion 
of the properties of this form 
for arbitrary fields $\KK$.

The $R$-matrix commutation relations in 
$R_q[G]$ (see \cite[Theorem I.8.15]{BG}) imply:

\ble{comm} Let $\la_i \in P^+$, $\nu_i \in P$,
$i=1,2$ and $\xi_2 \in (V(\la_2)^*)_{-\nu_2}$.

(i) For all $\mu_1 \in P$, $v_1 \in V(\la_1)_{\mu_1}$
and $\xi_1 \in (V(\la_1)^*)_{-\nu_1}$:
\begin{multline*}
c_{\xi_1, v_1}^{\la_1} c_{\xi_2, \la_2}^{\la_2} =
q^{ \lcor \mu_1, \la_2 \rcor - \lcor \nu_1, \nu_2 \rcor}
c_{\xi_2, \la_2}^{\la_2} c_{\xi_1, v_1}^{\la_1} +
\\
\sum_{\ga \in Q^+, \ga \neq 0}
\sum_{i=1}^{m(\ga)}
q^{ \lcor \mu_1, \la_2 \rcor - \lcor \nu_1 + \ga , \nu_2 - \ga  \rcor} 
c_{S^{-1}(u_{\ga, i})\xi_2, \la_2}^{\la_2} 
c_{S^{-1}(u_{-\ga, i}) \xi_1, v_1}^{\la_1}.
\end{multline*}

(ii) For all $\mu_2 \in P$, $v_2 \in V(\la_2)_{\mu_2}$
and $\xi_1 \in (V( - w_0 \la_1)^*)_{-\nu_1}$:
\begin{multline*}
c_{\xi_1, - \la_1}^{- w_0 \la_1} c_{\xi_2, v_2}^{\la_2} =
q^{ - \lcor \la_1, \mu_2 \rcor - \lcor \nu_1, \nu_2 \rcor} 
c_{\xi_2, v_2 }^{\la_2} c_{\xi_1, -\la_1}^{- w_0 \la_1} +
\\
\sum_{\ga \in Q^+, \ga \neq 0}
\sum_{i=1}^{m(\ga)}
q^{ - \lcor \la_1, \mu_2 \rcor - \lcor \nu_1 + \ga, \nu_2 - \ga \rcor} 
c_{S^{-1}(u_{\ga, i})\xi_2, v_2}^{\la_2} 
c_{S^{-1}(u_{-\ga, i}) \xi_1, -\la_1}^{-w_0 \la_1}.
\end{multline*}
\ele

Thus for all $\la \in P^+$, $w \in W$,
$\nu, \mu \in P$ and $c \in R_q[G]_{-\nu, \mu}$
\begin{align}
\label{n1}
c^+_{w, \la} c &= q^{\lcor w \la, \nu \rcor - \lcor \la, \mu \rcor} 
c c^+_{w, \la} \mod I_{w}^+ R^-,  
\\
\label{n2}
c^-_{w, \la} c &= q^{ \lcor w \la, \nu \rcor - \lcor \la, \mu \rcor} 
c c^-_{w, \la} \mod R^+ I_{w}^-. 
\end{align}

By abuse of notation we will denote 
the images of $c_{\xi, v}^\la$ and $c^\pm_{w, \la}$ in 
$R^\pm/I^\pm_{w_\pm}$ and $R_q[G]/I_{\bfw}$ 
by the same symbols (recall \eqref{Qbfw}), as 
it is commonly done in \cite{J1,J,HLT}.
All  $c^\pm_{w, \la} \in R/I^\pm_{w_\pm}$ are nonzero 
normal elements, see \eqref{n1}--\eqref{n2}.
Their images in $R_q[G]/I_{\bfw}$ 
are also nonzero normal elements.  
Denote the multiplicative subsets 
of $R^\pm$, $R^\pm/I^\pm_{w^\pm}$ and $R_q[G]/I_{\bfw}$
\begin{equation}
\label{Ewpm}
E^\pm_{w_\pm} = \{c^\pm_{w_\pm, \la} \mid \la \in P^+\}.
\end{equation}
Denote the multiplicative subset of $R_q[G]$ and $R_q[G]/I_{\bfw}$  
\begin{equation}
\label{Ew}
E_{\bfw} = E^+_{w_+} E^-_{w_-},
\end{equation}
the localization
\begin{equation}
\label{Rw}
R_{\bfw} =(R_q[G]/I_{\bfw})[E_{\bfw}^{-1}], 
\end{equation}
and its center
\begin{equation}
\label{Zw}
Z_{\bfw} = Z(R_{\bfw}).
\end{equation}

Since the ideal $I_{\bfw}$ is homogeneous with respect to the 
$P \times P$-grading \eqref{Rgrade}
of $R_q[G]$, $R_q[G]/I_{\bfw}$ inherits a $P \times P$-grading.
Denote the corresponding components
\begin{equation}
\label{Rgr}
(R_q[G]/I_{\bfw})_{\nu, \mu} = (R_q[G]_{\nu, \mu} + I_{\bfw})/I_{\bfw}, 
\quad \nu, \mu \in P.
\end{equation}
The elements of $E_{\bfw}$ are $P\times P$-homogeneous. Thus 
$R_{\bfw}$ also inherits a $P\times P$-grading. Its 
components will be denoted by $(R_{\bfw})_{\nu, \mu}$.

Recall that 
$c_{w, \la_1}^\pm c_{w, \la_2}^\pm = 
c_{w, \la_1 + \la_2}^\pm$ for all $\la_1, \la_2 \in P^+$.
Write $\la \in P$ as $\la = \la_1 - \la_2$ for some 
$\la_1, \la_2 \in P^+$ and define
\begin{equation}
\label{cP}
c^\pm_{w, \la} = c^\pm_{w, \la_+} (c^\pm_{w, \la_-})^{-1} 
\in R_{\bfw}.
\end{equation}
This definition does not depend on the choice of $\la_1$ 
and $\la_2$, because of the above mentioned property of 
the elements $c_{w, \la}^\pm$. We have:
\begin{equation}
\label{mult}
c^\pm_{w, \la_1} c^\pm_{w, \la_2} = c^\pm_{w, \la_1 + \la_2}, \quad
\forall \la_1, \la_2 \in P.
\end{equation}

Eqs. \eqref{n1}--\eqref{n2} imply that
\begin{equation}
\label{n3}
c^\pm_{w, \la} c = q^{\lcor w \la, \nu \rcor - \lcor \la, \mu \rcor} 
c c^\pm_{w, \la},
\end{equation} 
for all $\la, \nu, \mu \in P$ and 
$c \in (R_{\bfw})_{-\nu, \mu}$.
\subsection{Spectral decomposition theorem for $R_q[G]$}
\label{2.6} Consider the torus $\Tset^r = (\KK^*)^{\times r}$ 
and define the characters
\begin{equation}
\label{tmu}
t \mt 
t^\mu = \prod_{i=1}^r t_i^{\lcor \mu, \al_i\spcheck \rcor}, \quad 
t = (t_1, \ldots, t_r) \in \Tset^r, \mu \in P.
\end{equation}
There are two commuting rational $\Tset^r$-actions on 
$R_q[G]$ by $\KK$-algebra automorphisms:
\begin{equation}
\label{Tract}
t \cdot c = t^{\mu} c, \quad
t \in \Tset^r, c \in R_q[G]_{-\nu, \mu}, \nu, \mu \in P
\end{equation}
and
\begin{equation}
\label{Tract0}
t \cdot c = t^{\nu} c, \quad
t \in \Tset^r, c \in R_q[G]_{-\nu, \mu}, \nu, \mu \in P. 
\end{equation}
These actions are extensions of the left and right actions 
\eqref{action} of $H$ on $R_q[G]$, respectively, 
under the embedding $H \hra \Tset^r$
given by $K_i \mt (1, \ldots, 1, q_i, 1, \ldots, 1)$, 
$i = 1, \ldots, r$, where $q_i= q^{d_i}$ is in position $i$. 

\bth{J-thm} (Joseph \cite{J1}, Hodges--Levasseur \cite{HL}) 
(i) For each prime ideal $J$ of $R_q[G]$, 
there exists a unique ${\bfw} \in W \times W$ such that 
$J \supseteq I_{\bfw}$ and $(J/I_{\bfw}) \cap E_{\bfw} = \emptyset.$ 

(ii) For each ${\bfw}=(w_+, w_-)$, the ring $Z_{\bfw}$ is isomorphic to a 
Laurent polynomial ring over $\KK$ of dimension $\dim \, \ker (w_+ - w_-)$. 
Moreover the stratum $\Spec_{\bfw} R_q[G] \subset \Spec R_q[G]$
of ideals corresponding to ${\bfw}$ by (i) is homeomorphic to $\Spec Z_{\bfw}$ 
via the map $\iota_{\bfw} \colon \Spec Z_{\bfw} \to \Spec_{\bfw} R_q[G]$
defined as follows. For each $J^0 \in \Spec Z_{\bfw}$, 
$\iota_{\bfw}(J^0)$ is the unique ideal of $R_q[G]$ containing $I_{\bfw}$
such that 
\[
\iota_{\bfw}(J^0)/I_{\bfw} = (R_{\bfw} J^0) \cap (R_q[G]/I_{\bfw}) .
\]

(iii) For each ${\bfw} \in W \times W$, the set of primitive ideals 
$\Prim_{\bfw} R_q[G]$ in the stratum 
$\Spec_{\bfw} R_q[G]$ is precisely $\iota_{\bfw}^{-1} (\Max Z_{\bfw})$.
If the base field $\KK$ is algebraically closed, then  $\Prim_{\bfw} R_q[G]$
is the $\Tset^r$-orbit of a single primitive ideal.
\eth

Hodges and Levasseur proved the theorem in the $A$ case 
in \cite{HL}. Joseph gave a proof in the general case \cite{J1}.
We refer the reader to Joseph's book \cite{J} for a detailed 
treatment of these and many other related results.
A multiparameter version of this result was obtained by 
Hodges, Levasseur, and Toro in \cite{HLT}. For part (i) 
see \cite[Corollary 6.4]{J1} and \cite[Theorem 4.4]{HLT}, 
and for part (iii) \cite[Theorem 9.2]{J1} and 
\cite[Theorem 4.16]{HLT}. Joseph states part (ii) 
of \thref{J-thm} in terms of orbits of $\Zset_2^{\times r}$, 
see \cite[Theorem 8.11]{J1}, \cite[Theorem 10.3.4]{J}.
In the above form it is stated in Hodges--Levasseur--Toro 
\cite[Theorem 4.15]{HLT}. Brown, Goodearl and Letzter
\cite{BG0,GL} observed that the strata of $\Spec R_q[G]$ 
can be also described by
\[
\Spec_{\bfw} R_q[G] = \{ J \in \Spec R_q[G] \mid
\cap_{t \in \Tset^r} t \cdot J = I_{\bfw} \}
\]
(with respect to either \eqref{Tract} or \eqref{Tract0})
and developed this point of view to a general stratification 
method for the spectra of algebras with torus actions 
\cite{GL,BG,GK}. In \cite{J1,HLT} \thref{J-thm} is stated 
for $\KK = \Cset$, $q \in \Cset^*$ not a root of unity
and in \cite{J} for $\KK= k(q)$ for a field $k$ of 
characteristic 0. It is well known that 
the proofs of Joseph and Hodges--Levasseur--Toro
of \thref{J-thm} work for all base fields $\KK$,  
$q \in \KK^*$ not a root of unity,
as was noted in a similar context for the results
in \thref{Hprim}.

Joseph \cite{J1,J} determined the centers $Z_{\bfw}$ up to a 
finite extension. The next section contains a detailed 
discussion of this and an explicit description of $Z_{\bfw}$. 
It follows from \thref{Hprim} (ii), as well as from
\thref{J-thm} (ii), that: (1) the ideals $I_{\bfw}$, 
${\bfw} \in W \times W$ are stable under both 
actions \eqref{Tract} and \eqref{Tract0} 
of $\Tset^r$ on $R_q[G]$, and (2) every 
prime ideal of $R_q[G]$ which is 
$\Tset^r$-stable 
under \eqref{Tract} or \eqref{Tract0} is of this form.

We also note that the algebras $R_{\bfw}$ 
play an important role in the work 
of Berenstein and Zelevinsky \cite{BZ} 
on quantum cluster algebras. They are 
quantizations of the coordinate 
rings of double Bruhat cells in simple Lie groups, 
which were proved to be upper cluster algebras
by Berenstein, Fomin and Zelevinsky \cite{BFZ}.
\subsection{The De Concini--Kac--Procesi algebras}
\label{2.7} 
Recall from \S \ref{2.2} that the braid group $\B_\g$ associated to $\g$ acts 
on $\UU_q(\g)$ by algebra automorphisms.

Fix $w \in W$. Let
\begin{equation}
\label{wdecomp}
w = s_{i_1} \ldots s_{i_l}
\end{equation}
be a reduced expression of $w$. Recall that the roots 
in $\De^+ \cap w(\De^-)$ are given by
\begin{equation}
\label{beta}
\beta_1 = \al_{i_1}, \beta_2 = s_{i_1} (\al_{i_2}), 
\ldots, \beta_l = s_{i_1} \ldots s_{i_{l-1}} (\al_{i_l}).
\end{equation}
Define Lusztig's root vectors
\begin{equation}
X^{\pm}_{\beta_1} = X^{\pm}_{i_1}, 
X^{\pm}_{\beta_2} = T_{i_1} (X^\pm_{i_2}), 
\ldots, X^\pm_{\beta_l} = T_{i_1} \ldots T_{i_{l-1}} (X^{\pm}_{i_l}),
\label{rootv}
\end{equation}
see \cite[\S 39.3]{L} for details. 
The elements $X^\pm_{\be_k}$ satisfy the
Levendorskii--Soibelman straightening rule \cite{LS}:
\begin{multline}
\label{LS0}
X^\pm_{\be_i} X^\pm_{\be_j} - q^{\lcor \be_i, \be_j \rcor }
X^\pm_{\be_j} X^\pm_{\be_j}  \\
= \sum_{ {\bf{k}} = (k_{i+1}, \ldots, k_{j-1}) \in \Nset^{\times (j-i-2)} }
p_{\bf{k}}^\pm (X^\pm_{\be_{j-1}})^{k_{j-1}} \ldots (X^\pm_{\be_{i+1}})^{k_{i+1}},
\; \; p_{\bf{k}}^\pm \in \KK,
\end{multline}
for $i < j$. We refer to \cite[Proposition I.6.10]{BG} 
for the plus case of \eqref{LS0} for the version of 
$\UU_q(\g)$ used in this paper. The minus case follows from  
it by applying the algebra automorphism $\om$ of $\UU_q(\g)$ 
defined by 
\[
\om(X^\pm_i) = X^\mp_i, \;  
\om(K_i) = K_i^{-1}, \; \; 
i = 1, \ldots, r
\]
on its generators. It satisfies 
\[
\om( T_i(x)) = 
(- q_i)^{ \lcor \al_i\spcheck, \ga \rcor }  T_i( \om (x)),
\; \; \forall \ga \in Q, x \in (\UU_q(\g))_\ga,
\] 
cf. \cite[eq. 8.14(9)]{Ja}.

De Concini, Kac and Procesi defined 
\cite{DKP} the subalgebras $\UU_\pm^w$ of $\UU_\pm$ generated by 
$X^{\pm}_{\beta_j}$, $j=1, \ldots, l$ and proved the following result:

\bth{DKP} (De Concini, Kac, Procesi) \cite[Proposition 2.2]{DKP}
The algebras $\UU_\pm^w$ do not depend on the choice of a 
reduced expression of $w$ and have the PBW basis
\begin{equation}
\label{vect}
(X^\pm_{\be_l})^{n_l} \ldots (X^\pm_{\be_1})^{n_1}, \; \; 
n_1, \ldots, n_l \in \Nset.
\end{equation}
\eth

Lusztig established independently \cite[Proposition 40.2.1]{L} that the 
space span\-ned by the monomials \eqref{vect} does not depend on the 
choice of a reduced expression of $w$.

 
In relation to \thref{DKP},
for ${\bf{n}}= (n_1, \ldots, n_l) \in \Nset^{\times l}$ 
denote the monomial 
\begin{equation}
\label{monomial}
(X^\pm)^{\bf{n}} = (X^\pm_{\be_l})^{n_l} \ldots (X^\pm_{\be_1})^{n_1}.
\end{equation}
These monomials form a $\KK$-basis of $\UU^w_\pm$.
We will say that $(X^\pm)^{\bf{n}}$ has {\em{degree}} ${\bf{n}}$.
Introduce the lexicographic order 
on $\Nset^{\times l}$:
\begin{multline}
\label{lexi}
{\bf{n}}= (n_1, \ldots, n_l) < {\bf{m}}= (m_1, \ldots, m_l), \; \; 
\mbox{if there exists} \; \; j \in \{1, \ldots, l\}
\\
\mbox{such that} \; \; n_j <m_j \; \; \mbox{and} \; \;  
n_{j+1} = m_{j+1}, \ldots, n_l = m_l.
\end{multline}

We will say that the {\em{highest term}} of 
a nonzero element
$u \in \UU^w_\pm$ is $p (X^\pm)^{\bf{n}}$, 
where ${\bf{n}} \in \Nset^{\times l}$ and 
$p \in \KK^*$, if
\[
u - p (X^\pm)^{\bf{n}} \in \Span 
\{ (X^\pm)^{{\bf{n}}'} \mid {\bf{n}}'  \in \Nset^{ \times l}, 
{\bf{n}'} < {\bf{n}} \}. 
\]
The Levendorskii--Soibelman straightening rule implies that 
one obtains an $\Nset^{ \times l}$-filtration on $\UU^w_\pm$ by
collecting the elements with highest terms of degree 
$\leq {\bf{n}}$ for ${\bf{n}} \in \Nset^{\times l}$:

\ble{LSmult} For all ${\bf{n}}, {\bf{n}}' \in \Nset^{\times l}$
the highest term of the product $(X^\pm)^{\bf{n}} (X^\pm)^{ {\bf{n}}' }$
is $q^{ m_{ {\bf{n}}, {\bf{n}}'} } (X^\pm)^{ {\bf{n}} + {\bf{n}}' }$, 
for some $m_{ {\bf{n}}, {\bf{n}}'} \in \Zset$.
\ele
\subsection{A second presentation of $\UU^w_\pm$}
\label{2.8} 
The algebras $\UU^w_\pm$ are (anti)isomorphic 
(see \thref{isom} below) to the algebras 
$S^\mp_w$ 
which play an important role in Joseph's work 
\cite{J1,J}. The latter algebras are defined as follows. 
Let $w \in W$.
The quotients $R^\pm/I^\pm_w$ can be canonically identified as vector 
spaces with $\oplus_{\la \in P^+} V_w^\pm(\la)^*$ by
$c^{\la_1}_{\xi_1, \la_1} \mapsto g^+_w(\xi_1)$ and
$c^{-w_0 \la_2}_{\xi_2, - \la_2} \mapsto g^-_w(\xi_2)$,
for $\la_1, \la_2 \in P^+$, $\xi_1 \in V(\la_1)^*$,
$\xi_2 \in V(-w_0 \la_2)^*$,
where we used the projections \eqref{g}.
Recall that $R^\pm_w = (R^\pm/I^\pm_w)[E^\pm_w]$.
The invariant subalgebras of $R^\pm_w$ with respect to the 
left action \eqref{action} of $H$ will be denoted by $S^\pm_w$. 
In terms of the above vector space identifications
\begin{equation}
\label{ident}
S^\pm_w = \underrightarrow{\lim}_{\la \in P^+} 
(c^\pm_{w, \la})^{-1} V^\pm_w(\la)^*.
\end{equation}
For $\la_1, \la_2 \in P^+$ the embedding 
\[
(c^\pm_{w, \la_2})^{-1} V^\pm_w(\la_2)^* \hra 
(c^\pm_{w, \la_1+\la_2})^{-1} V^\pm_w(\la_1+\la_2)^*
\]
is given by $(c^\pm_{w, \la_2})^{-1} \xi  
\mt (c^\pm_{w, \la_1+\la_2})^{-1} 
(g^\pm_w(\xi_{w, \la_1}^\pm) . \xi)$,
where $\xi \in V_w^\pm(\la_2)^*$. The product in the right hand 
side is the Cartan multiplication \eqref{Cartan} and 
$\xi^\pm_{w, \la_1}$ are the weight vectors,
defined in \S \ref{2.4}.
The $P \times P$-grading of $R_q[G]$ induces 
$P \times P$-gradings 
on $R^\pm/I^\pm_w$, $R^\pm_w$, 
and $S^\pm_w$, analogously to \eqref{Rgr}. 
Denote the graded components of the algebra $S^\pm_w$ by 
$(S^\pm_w)_{\nu, \mu}$, $\nu, \mu \in P$.
It is clear that $(S^\pm_w)_{\nu, \mu}=0$, 
if $\nu \notin Q$ or $\mu \neq 0$. Thus, effectively 
we have a $Q$-grading on $S^\pm_w$. 
Eq. \eqref{gradS+-} below describes the nonzero 
components of this grading.

The $Q$-grading \eqref{Uq-grading} of $\UU_q(\g)$ induces
a $Q$-grading of the algebras $\UU^w_\pm$, explicitly 
given by 
\begin{equation}
\label{U-grading}
\deg X^\pm_{\beta_j} = \pm \beta_j, \quad j=1, \ldots, l.
\end{equation}
It is clear that 
\begin{equation}
\label{non0U}
(\UU^w_\pm)_\ga \neq 0 \; \; 
\mbox{if and only if} \; \; 
\pm \ga \in \sum_{j=1}^l 
\Nset \beta_j.
\end{equation}
The group $H$ acts on $\UU_q(\g)$ by conjugation. The subalgebras
$\UU^w_\pm$ are stable under this action. The eigenspaces for the 
action are precisely the graded components with respect to
the grading \eqref{U-grading}.

For $\ga \in Q^+$, $\ga \neq 0$ denote
$m_w(\ga) = \dim (\UU^w_+)_\ga= \dim (\UU^w_-)_{-\ga}$.
Fix a pair of dual bases $\{u_{\ga, i} \}_{i=1}^{m_w(\ga)}$ and
$\{u_{-\ga, i} \}_{i=1}^{m_w(\ga)}$ of $(\UU^w_+)_\ga$ and 
$(\UU^w_-)_{-\ga}$ with respect to the Rosso--Tanisaki form,
see \cite[Ch. 6]{Ja}. The quantum $R$ matrix corresponding 
to $w$ is given by 
\[
\RR^w = \sum_{\ga \in Q_+} \sum_{i=1}^{m_w(\ga)} 
u_{\ga, i} \otimes u_{- \ga, i} \in \UU_+^w \wh{\otimes} \UU_-^w.
\]
Here $\UU_+^w \wh{\otimes} \UU_-^w$ denotes
the completion of $\UU_+^w \otimes_\KK \UU_-^w$ with respect to the
descending filtration \cite[\S 4.1.1]{L}. More explicitly, 
for a reduced expression of $w$ as in \eqref{wdecomp}, 
$\RR^w$ is given by 
\begin{multline}
\label{Rwe}
\RR^w = \sum_{n_1, \ldots, n_l \in \Nset}
\left( \prod_{j=1}^l (-1)^{n_j} q_{i_j}^{- n_j (n_j -1)/2} 
\frac{ (q_{i_j} - q_{i_j}^{-1})^{n_j}}{[n_j]_{q_{i_j}}!} \right)   
\times
\\
(X_{\be_l}^+)^{n_l} \ldots (X_{\be_1}^+)^{n_1} 
\otimes
(X_{\be_l}^-)^{n_l} \ldots (X_{\be_1}^-)^{n_1}
\end{multline}
in terms of the notation \eqref{beta}--\eqref{rootv}, 
see \cite[eq. 8.30(2)]{Ja}. This implies that 
\begin{equation}
\label{Rw-prod}
\RR^{w_0} = \left( T_w (\RR^{w^{-1} w_0}) \right) \RR^w.
\end{equation}

Recall that there is a unique graded algebra antiautomorphism 
$\tau$ of $\UU_q(\g)$ such that 
\begin{equation}
\label{tau}
\tau(X_i^\pm) = X_i^\pm, 
\, 
\tau(K_i) = K_i^{-1}, \; \; 
i = 1, \ldots, r,
\end{equation}
see \cite[Lemma 4.6(b)]{Ja}. It satisfies
\begin{equation}
\label{tau-ident}
\tau (T_w x) = T_{w^{-1}}^{-1} ( \tau (x)), \; \; 
\forall x \in \UU_q(\g), w \in W,
\end{equation}
cf. \cite[ eq. 8.18(6)]{Ja}. 

For $w \in W$, we define the maps 
\[
\varphi_w^\pm \colon S^\pm_w \to \UU^w_\mp 
\]
by
\begin{align}
\varphi_w^+( c^\la_{\xi, \la} ( c^+_{w, \la} )^{-1} )
&=( c^\la_{\xi, T_w v_\la} \otimes \id) (\tau \otimes \id) \RR^w \; \; 
\mbox{and}
\label{map+}
\\
\varphi_w^-( (c^-_{w, \la})^{-1} c^{-w_0 \la}_{\xi', - \la} )
&=(\id \otimes c^{- w_0 \la}_{\xi', T_{w^{-1}}^{-1} v_{-\la}}) (\id \otimes \tau) \RR^w,
\label{map-}
\end{align}
for $\la \in P^+$, $\xi \in V(\la)^*$, $\xi' \in V(- w_0 \la)^*$. 
In the right hand sides the elements of $R_q[G]$ are 
viewed as functionals on $\UU_q(\g)$. 
The choice of $T_{w^{-1}}^{-1}$ in \eqref{map-} 
instead of $T_w$
matches the second normalization in \eqref{cd-def}
and the Poisson side of the picture discussed 
in \S \ref{4.3}--\ref{4.4}.

\bth{isom} The maps 
$\varphi_w^+ \colon S^+_w \to \UU^w_-$
are well defined antiisomorphisms
of $Q$-graded algebras. The maps 
$\varphi_w^- \colon S^-_w \to \UU^w_+$
are well defined isomorphisms
of $Q$-graded algebras.
\eth

\thref{isom} is an analog of \cite[Theorem 3.7]{Y}.
In \cite{Y} we used a version of $\UU_q(\g)$ equipped with 
the opposite coproduct, a different braid group action and 
Lusztig's root vectors. As a result of this
the map in \cite[Theorem 3.7]{Y} is an isomorphism.
In \cite{Y} we also formulated the result for a base field
$\KK$ of characteristic $0$ and $q \in \KK$ transcendental
over $\Qset$. Because of this we will give a proof of 
\thref{isom}. We will need the following simple lemma.

\ble{hom} \cite[Lemma 3.2]{Y}
Let $\HH$ be a Hopf algebra and $A$ be an 
$\HH$-module algebra equipped with a right $\HH$-action.
If $\epsilon \colon A \to \KK$ is an algebra homomorphism, 
where $\KK$ is the ground field, then
the map $\phi \colon A \to \HH^*$ given by
\[
\phi(a) (h) = \epsilon (a . h)
\]
is an algebra homomorphism. If, in addition the action of $\HH$ on $A$
is locally finite, then the image of $\phi$ is contained in
the restricted dual $\HH^\circ$ of $\HH$.
\ele
\noindent
{\em{Proof of \thref{isom}.}} We will prove the plus case. 
The minus case is analogous and is left to the reader. 
In the definition of the maps $\varphi^\pm_w$ the inverses of the 
elements $c^\pm_{w, \la}$ appear on different sides 
because of the differences between the coproducts of 
$X^\pm_i$.

It follows from \eqref{RwT} that for all $w \in W$ the map
\[
\epsilon_w \colon R^+ \to \KK \quad 
\mbox{defined by} \quad
\epsilon_w ( c_{\xi, \la}^\la ) = \xi (T_w v_\la), \; \; 
\la \in P^+, \xi \in V(\la)^*  
\]
is an algebra homomorphism. Denote by $\UU_q(\b_+)$ the 
(Hopf) subalgebra of $\UU_q(\g)$ generated by 
$X_i^+$, $K_i$, $i = 1, \ldots, r$.
We apply \leref{hom} to
$A = R^+$, $\HH = \UU_q(\b_+)$, $\epsilon = \epsilon_w$ 
and the restriction of the right action \eqref{action}
of $\UU_q(\g)$ on $R^+$ to $\UU_q(\b_+)$. This action is locally finite.
Denote the corresponding homomorphism from \leref{hom} 
by $\phi_w^+ \colon R^+ \to (\UU_q(\b_+))^\circ$.
We will identify $(\UU_+)_0$ with $\KK$ 
via $t. 1 \mt t$, $t \in \KK$. 
For all $x \in (\UU_+)_\ga$, $\ga \in Q^+$, 
$n_1, \ldots, n_r \in \Zset$
\[ 
\lcor \phi_w^+( c_{w,\la}^+), x \prod_{i=1}^r K_i^{n_i} \rcor = 
\delta_{\ga, 0} \prod_{i=1}^r q^{ n_i \lcor \al_i, w \la \rcor} x.
\]
One easily deduces from this that 
$\phi_w^+( c_{w,\la}^+) \in \UU_q(\b_+)$, $\la \in P^+$
are not zero divisors. More generally for all 
$x \in (\UU_+)_\ga$, $\ga \in Q^+$, 
$n_1, \ldots, n_r \in \Zset$
\[ 
\lcor \phi_w^+( c_{\xi, \la}^\la), x \prod_{i=1}^r K_i^{n_i} \rcor = 
\prod_{i=1}^r q^{ n_i \lcor \al_i, w \la \rcor} \xi(x T_w v_\la).
\]
Thus $I^+_w \subset \ker \phi^+_w$ and $\phi^+_w$ induces an
algebra homomorphism from $R^+/I_w^+$ to $(\UU_q(\b_+))^\circ$,  
which by abuse of notation will be denoted by the same symbol.
For $\la \in P^+$, $\xi \in (V(\la))^*$ define the 
elements $a_{w, \la} \in (\UU_q(\b_+))^\circ$ by 
\begin{equation}
\label{aaa1}
\lcor a_{w, \la}, x \prod_{i=1}^r K_i^{n_i} \rcor
: = \xi(x T_w v_\la), \; \; 
\forall x \in (\UU_+)_\ga, \ga \in Q^+,
n_1, \ldots, n_r \in \Zset.
\end{equation}
The formulas for the coproduct
of $\UU_q(\g)$ imply that  
\begin{equation}
\label{a-elem}
a_{w, \la} \phi_w^+( c_{w,\la}^+)
=
\phi_w^+( c_{\xi, \la}^\la).
\end{equation}
The fact that $\phi_w^+( c_{w,\la}^+) \in (\UU_q(\b_+))^\circ$
are not zero divisors imply that the assignment
\begin{equation}
\label{aaa2}
\phi_w^+( c_{\xi, \la}^\la (c_{w,\la}^+)^{-1} ) 
:= a_{w, \la}, \; \; \la \in P^+, \xi \in V(\la)^*
\end{equation}
is a well defined algebra homomorphism from 
$S_w^+$ to $(\UU_q(\b_+))^\ci$. (It will be denoted by 
the same symbol $\phi_w^+$ as the previous homomorphism.) 
We have the embedding of algebras $\UU_- \hra (\UU_q(\b_+))^\ci$ 
via the Rosso--Tanisaki form. It follows from \eqref{aaa1}
and \eqref{aaa2} that $\phi_w^+(S_w^+) \subseteq \UU_-$ 
and that $\phi_w^+ \colon S_w^+ \to \UU_-$ is
injective.
Moreover, \eqref{aaa1} and \eqref{aaa2} imply that 
$\phi_w^+ \colon S_w^+ \to \UU_-$ is given by 
\[
\phi_w^+ (  c_{\xi, \la}^\la (c_{w,\la}^+)^{-1} ) = 
( c^\la_{\xi, T_w v_\la} \otimes \id ) \RR^{w_0}, \; \; 
\forall \la \in P^+, \xi \in V(\la)^*.
\]
Define the injective algebra antihomomorphism
\[
\varphi_w^+ := \tau \phi_w^+ \colon S_w^+ \to \UU_-.
\] 
The fact that $(\tau \otimes \tau) \RR^{w_0} = \RR^{w_0}$,
see \cite[eq. 7.1(2)]{Ja}, implies that $\varphi_w^+$ is given by 
\begin{equation}
\label{vpp}
\varphi_w^+ (  c_{\xi, \la}^\la (c_{w,\la}^+)^{-1} ) = 
( c^\la_{\xi, T_w v_\la} \tau \otimes \id ) \RR^{w_0}, \; \; 
\forall \la \in P^+, \xi \in V(\la)^*.
\end{equation}  
It follows from this formula that $\varphi_w^+ \colon S^+_w \to \UU_-$ is an 
antihomomorphism of $Q$-graded algebras. 

Fix a reduced expression of $w$ as in \eqref{wdecomp} and extend it 
to a reduced expression $w_0 = s_{i_1} \ldots s_{i_l} s_{i_{l+1}} \ldots s_{i_N}$ 
of the longest element of $W$. We claim that $\varphi_w^+(S^+_w) \subseteq \UU^w_-$ 
and that $\varphi_w^+$ is given by \eqref{map+}. Both statements follow 
from \eqref{Rw-prod} and the fact that 
\begin{equation}
\label{int1}
[\tau(T_{i_1} \ldots T_{i_{j-1}} (X_{i_j}^+))] V(\la)_{w \la} = 0, \; \; 
\forall j = l+1, \ldots, N.
\end{equation}
We have
\begin{align*}
&[\tau(T_{i_1} \ldots T_{i_{j-1}} (X_{i_j}^+))] T_{w^{-1}}^{-1} v_\la = 
(T_{i_1}^{-1} \ldots T_{i_{j-1}}^{-1} (X_{i_j}^+) ) (T_{i_1}^{-1} \ldots T_{i_l}^{-1} v_\la) 
\\
=&(T_{i_1}^{-1} \ldots T_{i_l}^{-1}) 
[ (T_{i_{l+1}}^{-1} \ldots T_{i_{j-1}}^{-1} (X_{i_j}^+)) v_\la ]=0, \; \; 
\forall j= l+1, \ldots, N,
\end{align*}
because $(T_{i_{l+1}}^{-1} \ldots T_{i_{j-1}}^{-1} (X_{i_j}^+)) \in \UU^+$.
This implies \eqref{int1} since $V(\la)_{w \la} = \KK T_{w^{-1}}^{-1} v_\la$. 
We have proved that $\varphi_w^+(S^+_w) \subseteq \UU^w_-$ and all that remains to be 
shown now is that $\varphi_w^+(S^+_w) = \UU^w_-$. Assume the opposite 
that $\varphi_w^+(S^+_w) \subsetneq \UU^w_-$. Since the Rosso--Tanisaki 
form restricts to a nondegenerate pairing between $\UU^w_+$ and $\UU^w_-$, 
$\varphi_w^+(S^+_w) \subsetneq \UU^w_-$ implies that 
there exists $\ga \in Q^+$ and $x \in (\UU^w_+)_\ga$, $x \neq 0$ 
such that $c^\la_{\xi, T_w v_\la} (x) = 0$ for all $\la \in P^+$, $\xi \in V(\la)^*$. 
This means that $x T_w v_\la = 0$, $\forall \la \in P^+$, $\xi \in V(\la)^*$. 
Set $x_1 = T_w^{-1}(x)$, $\ga_1 = - w^{-1}(\ga)$. Then 
$\ga_1 \in Q^+$, $x_1 \in (\UU_-)_{- \ga_1}$, $x_1 \neq 0$ and 
$x_1 v_\la = 0$ for all $\la \in P_+$. If we choose 
$\la \in P^+$ such that 
$\lcor \la, \al_i\spcheck \rcor > \lcor \ga_1, \om_1 + \ldots + \om_r \rcor$,
$\forall i =1, \ldots, r$, then the equality $x_1 v_\la = 0$ contradicts with 
\eqref{Vla}. This completes the proof of the theorem.
\qed
\sectionnew{A description of the centers of Joseph's localizations}
\lb{cent}
\subsection{Statement of the main result}
\label{3.1} 
In this section we obtain an explicit description of the centers 
$Z_{\bfw}$ of Joseph's localizations $R_{\bfw}$. This is done
in \thref{bfw-center}. It is the building block of the paper. On the one 
hand, it leads to a more explicit description of the prime ideals of 
$R_q[G]$, which in particular allows to compute the stabilizers 
of those ideals under the actions \eqref{Tract} and \eqref{Tract0}
of $\Tset^r$ and to construct a torus equivariant Dixmier type map in the 
next section. This description of prime ideals eventually 
leads to a classification of the maximal spectrum of $R_q[G]$,
which allows us to settle a question of Goodearl and Zhang \cite{GZ}, 
by proving that all maximal ideals of $R_q[G]$ have finite 
codimension. On the other hand, \thref{bfw-center} and the methods 
developed in its proof play a key role in two 
freeness theorems which we prove 
in Sections \ref{free1} and \ref{Module} for the De Concini--Kac--Procesi
algebras and Joseph's localizations $R_{\bfw}$. The first 
is a freeness result for $\UU^w_\pm$ as a module over its 
subalgebra generated by homogeneous normal elements, and the second is a freeness 
result for $R_{\bfw}$ over its subalgebra generated by 
Joseph's set of normal elements $(E_{\bfw})^{ \pm 1}$. 
The latter supplies the second key ingredient in the classification 
of $\Max R_q[G]$ in Section \ref{Max}.

For a subset $I \subset \{ 1, \ldots, r \}$ denote
\begin{equation}
\label{PQi}
P_I = \bigoplus_{i \in I} \Zset \om_i, \quad
P^+_I = \bigoplus_{i \in I} \Nset \om_i, \quad
Q_I = \bigoplus_{i \in I} \Zset \al_i, \quad 
Q_I\spcheck = \bigoplus_{i \in I} \Zset \al_i\spcheck.
\end{equation} 
For $w \in W$ set 
\begin{equation}
\label{IS}
\II(w) = \{ i = 1, \ldots, r \mid w(\om_i) = \om_i \}
\; \; 
\mbox{and} 
\; \; 
\SS(w) = \{1, \ldots, r \} \backslash \II(w).
\end{equation}
For ${\bfw} = (w_+, w_-) \in W \times W$ set 
\begin{equation}
\label{Ibfw}
\II({\bfw}) = \II(w_+) \cap \II(w_-) \; \; 
\end{equation}
and
\begin{equation}
\label{Sbfw}
\SS({\bfw}) = \SS(w_+) \cup \SS(w_-)= 
\{1, \ldots, r\} \backslash \II({\bfw}). 
\end{equation}
The intersection
\begin{equation}
\wt{\LL}({\bfw}) = \ker (w_+ - w_-) \cap P
\label{L}
\end{equation}
is a lattice of rank $\dim \, \ker (w_+ - w_-)$. Its reduced version 
\begin{equation}
\wt{\LL}_{\red}({\bfw}) = \ker (w_+ - w_-) \cap P_{\SS({\bfw})}
\label{Lred}
\end{equation}
is a lattice of rank
\begin{equation}
\label{k}
k = \dim \, \ker (w_+ - w_-) - |\II({\bfw})|,
\end{equation}
because $P_{\II({\bfw})} \subset \ker(w_+ - w_-)$ and thus
\begin{equation}
\label{LredL}
\wt{\LL}({\bfw}) = P_{\II({\bfw})} \oplus \wt{\LL}_{\red}({\bfw}).
\end{equation}

Choose a basis $\la^{(1)}, \la^{(2)}, \ldots, \la^{(k)}$
of $\wt{\LL}_{\red}({\bfw})$. For each $j = 1, \ldots, k$ denote
\begin{equation}
\label{a}
a_j = c_{w_+, \la^{(j)}}^+ (c_{w_-, \la^{(j)}}^-)^{-1},
\end{equation}
recall \eqref{cP}.

\bth{bfw-center} Assume that $\KK$ is an arbitrary base field, 
and $q \in \KK^*$ is not a root of unity. Then 
for each ${\bfw} =(w_+, w_-) \in W \times W$ the center $Z_{\bfw}$ 
of the algebra $R_{\bfw}$ coincides with the Laurent polynomial 
algebra over $\KK$ of dimension $\dim \ker (w_+ - w_-)$ with generators
\begin{equation}
\label{Laurent-cent}
\{ c^+_{w_+, \om_i} \mid i \in \II({\bfw}) \} \sqcup
\{a_1, \ldots, a_k \}.
\end{equation}
Here $k$ and $a_1, \ldots, a_k$ are given by \eqref{k} and \eqref{a}.
\eth

Kogan and Zelevinsky \cite{KZ} proved that similar equations are cutting the 
symplectic leaves of the standard Poisson structure on the corresponding 
connected, simply connected, complex, simple Lie group within a double 
Bruhat cell. Section \ref{Dixmier} will establish a connection between the 
two results. 

The cases of $\g= {\mathfrak{sl}}_2$ and $\g = {\mathfrak{sl}}_3$ 
of \thref{bfw-center} were obtained 
by Hodges--Levasseur \cite{HL0} and Goodearl--Lenagan \cite{GLen0},
respectively. Their methods are very different from ours and use in an 
essential way the low rank of the underlying Lie algebra.
\subsection{Associated root and weight spaces}
\label{3.2} 
Next, we gather some simple facts for the sets $\II(w)$ and $\SS(w)$,
$w \in W$.

\ble{Iw} Fix $w \in W$.

(i) Then $\SS(w) = \{ i = 1, \ldots, r \mid s_i \leq w \}$ 
with respect to the Bruhat order $\leq$ on $W$, i.e. 
for each reduced expression $w = s_{i_1} \ldots s_{i_l}$
\[
\SS(w) = \cup_{j=1}^l \{ i_j \}.
\]

(ii) We have
\begin{equation}
\label{wQ}
\sum_{\be \in \De^+ \cap w (\De^-)} \Zset \be = Q_{\SS(w)}, \; \; 
\sum_{\be \in \De^+ \cap w (\De^-)} \Zset \be\spcheck = Q_{\SS(w)}\spcheck
\end{equation}
and
\begin{equation}
\label{wP}
(\De^+ \cap w(\De^-))^\perp \cap P = 
(Q_{\SS(w)})^\perp \cap P = P_{\II(w)}. 
\end{equation}
\ele
\begin{proof} For the reduced expression in (i) denote
$S= \cup_{j=1}^l \{ i_j \}$ and $I = \{1, \ldots, r \} \backslash S$.
One has
\[
\De^+ \cap w(\De^-) = \{\beta_j =s_{i_1} \ldots s_{i_{j-1}} (\al_{i_j})
\mid j= 1, \ldots, l \},
\]
cf. \S \ref{2.7}. Since 
\[
\beta_j - \al_{i_j} \in \sum_{n=1}^{j-1} \Zset \al_{i_n}, \quad
\forall j = 1, \ldots, l,
\]
we have
\begin{equation}
\label{summ}
\sum_{\be \in \De^+ \cap w (\De^-)} \Zset \be = 
\bigoplus_{i \in S} \Zset \al_i = Q_S.
\end{equation}
Analogously 
\begin{equation}
\label{summ2}
\sum_{\be \in \De^+ \cap w (\De^-)} \Zset \be\spcheck = Q_{S}\spcheck.
\end{equation}

Obviously $I \subseteq \II(w)$. If $i \in \II(w)$, then for all 
$\be \in \De^+ \cap w(\De^-)$,
\[
0 \leq \lcor \om_i, \be \rcor = \lcor w^{-1} (\om_i), w^{-1} (\be) \rcor = 
\lcor \om_i, w^{-1} (\be) \rcor \leq 0, 
\]
thus $\lcor \om_i, \be \rcor = 0$. Taking \eqref{summ} into account, 
we obtain that $i \in \II(w)$ implies $\om_i \in (Q_S)^\perp \cap P = P_I$, 
i.e. $i \in I$. Therefore $I = \II(w)$ and $S = \SS(w)$. 
Now the second part follows from \eqref{summ} and \eqref{summ2}.
\end{proof}
\noindent
\subsection{One side inclusion in \thref{bfw-center}}
\label{3.3}
Joseph proved \cite{J1} that 
\begin{equation}
\label{cent1}
c_{w_+, \la}^+ (c_{w_-, \la}^-)^{-1} \in Z_{\bfw}, \quad
\mbox{for all} \; \; \la \in \LL({\bfw}).
\end{equation}
This follows from \eqref{n3}. In particular,
in the setting of \S \ref{3.1},
\[
a_j \in Z_{\bfw}, \quad \forall j =1, \ldots, k.
\] 
The following proposition provides the rest needed to claim that $Z_{\bfw}$ 
contains all elements in \eqref{Laurent-cent}.

\bpr{cent} For all ${\bfw} = (w_+, w_-) \in W \times W$ and $i \in \II({\bfw})$,
\[
c_{w_+, \om_i}^+ \in Z_{\bfw}.
\]
\epr
\begin{proof}
Fix $i \in \II({\bfw})$. Since 
$R_q[G] = R^+ R^-$, it is sufficient to prove that $c_{w_+, \om_i}^+$
commutes with the images of $R^+$ and $R^-$ in $R_{\bfw}$. We will prove
the former. The latter is analogous and is left to the reader. 
Let $\la \in P^+$. Recall that $(V_{w_+}(\la))_\nu \neq 0$ implies that 
$\nu = w_+ (\la) + \ga$ for some 
\begin{equation}
\label{gamma}
\ga \in \sum_{\be \in \De^+ \cap w_+ (\De^-)} \Nset \be \subset Q_{\SS(w_+)},
\end{equation}
cf. \leref{Iw}.
The definition of $I_{w_+}^+$ implies that,
if the image of $c^\la_{\xi, \la}$ in $R_{\bfw}$ is nonzero for some 
$\xi \in (V(\la)^*)_{- \nu}$, then $\nu = w_+(\la) + \ga$ 
with $\ga$ as in \eqref{gamma}, in particular
$\ga \in Q_{\SS(w_+)}$.
Using \eqref{n3}, we obtain that
\[
c_{w_+, \om_i}^+ c^\la_{\xi, \la} =
q^{ \lcor w \om_i, w \la + \ga \rcor - \lcor \om_i, \la \rcor } 
c^\la_{\xi, \la} c_{w_+, \om_i}^+ = q^{\lcor \om_i, \ga \rcor} 
c^\la_{\xi, \la} c_{w, \om_i}^+ = c^\la_{\xi, \la} c_{w, \om_i}^+
\]  
in $R_{\bfw}$, since 
$i \in \II({\bfw}) \subseteq \II(w_+)$
implies $q^{\lcor \om_i, \ga \rcor} =1$, 
$\forall \ga \in Q_{\SS(w_+)}$, see \leref{Iw} (ii).  
This completes the proof of \prref{cent}.
\end{proof}
\subsection{Joseph's description of $R_{\bfw}$}
\label{3.4}
Our treatment of $R_{\bfw}$ and its center uses 
a model of Joseph of $R_{\bfw}$, which represents 
it as a kind of ``bicrossed product'' of the algebras 
$S^\pm_{w_\pm}$, modulo a simple additional localization 
and a smash product by a Laurent polynomial ring. 
We refer the reader to \cite[ \S 9.1-9.2 and \S 10.3]{J} 
for details. 
This model 
and \thref{isom} allow the simultaneous application 
of techniques from quantum function
algebras (for the algebra $R_q[G]$, its quotients and localizations)  
and quantized universal enveloping algebras of nilpotent 
Lie algebras (for the algebras $\UU^w_\pm$).

First, denote by $R^+ \ast R^-$ the free product of the 
$\KK$-algebras $R^+$ and $R^-$. 
Define $R^+ \circledast R^-$ as the quotient of $R^+ \ast R^-$
by the following relations (which are 
analogous to the ones in \leref{comm}):
\begin{multline}
\label{commR}
c_{\xi_1, - \la_1 }^{- w_0 \la_1} c_{\xi_2, \la_2}^{\la_2} =
q^{ - \lcor \la_1, \la_2 \rcor - \lcor \nu_1, \nu_2 \rcor} 
c_{\xi_2, \la_2}^{\la_2} c_{\xi_1, -\la_1 }^{-w_0 \la_1} +
\\
\sum_{\ga \in Q^+, \ga \neq 0}
\sum_{i=1}^{m(\ga)}
q^{ - \lcor \la_1, \la_2 \rcor - \lcor \nu_1 + \ga , \nu_2 - \ga \rcor} 
c_{S^{-1}(u_{\ga, i})\xi_2, \la_2}^{\la_2} 
c_{S^{-1}(u_{-\ga, i}) \xi_1, -\la_1 }^{ -w_0 \la_1},
\end{multline}
for all $\la_i \in P^+$, $\nu_i \in P$, $i=1,2$,
$\xi_1 \in (V(- w_0 \la_1)^*)_{-\nu_1}$, 
$\xi_2 \in (V(\la_2)^*)_{-\nu_2}$.
Joseph proved \cite[Lemma 9.1.8]{J}
that the multiplication map 
in $R^+ \circledast R^-$
induces the vector space isomorphism
\begin{equation}
\label{R+-}
R^+ \otimes_\KK R^- \congto R^+ \circledast R^-
\end{equation}
and that $R^+ \circledast R^-$ is a noetherian domain, 
\cite[Lemma 9.1.9 (ii) and Proposition 9.1.11]{J}.
He also proved that the multiplication 
map $R^+ \otimes_\KK R^- \to R_q[G]$
induces a surjective $\KK$-algebra homomorphism
$R^+ \circledast R^- \to R_q[G]$ and described
its kernel in \cite[Corollary 9.2.4]{J}.

For the remainder of this section we fix 
${\bfw} = (w_+, w_-) \in W \times W$.
By \cite[Corollary 10.1.10]{J} 
\[
\wh{I}_{\bfw} = I^+_{w_+} R^- + R^+ I^-_{w_-}
\]
is a completely prime ideal of $R^+ \circledast R^-$.
The embeddings $R^\pm \hra R^+ \circledast R^-$ 
induce \cite[\S 10.3.1]{J} embeddings 
$R^\pm/I^\pm_{w_\pm} \hra (R^+ \circledast R^-)/I_{\bfw}$.
The images of $c^\pm_{w_\pm, \la}$ are nonzero 
normal elements in $R^\pm/I^\pm_{w_\pm}$ and 
$(R^+ \circledast R^-)/I_{\bfw}$. These 
images will be denoted by the same symbols. 
Recall the definition \eqref{Ewpm} of
the multiplicative subsets $E^\pm_{w_\pm}$ of
$R^\pm/I^\pm_{w_\pm}$. Define the 
multiplicative subset
\[
\wh{E}_{\bfw} = E^+_{w_+} E^-_{w_-}
\]
of $(R^+ \circledast R^-)/I_{\bfw}$ and
denote the localization
\[
\wh{R}_{\bfw} = ((R^+ \circledast R^-)/I_{\bfw})[\wh{E}_{\bfw}^{-1}].
\] 
Recall the definition of the subalgebras $S^\pm_{w_\pm}$ of $R^\pm_{w_\pm}$
from \S \ref{2.8}.
The embeddings $R^\pm/I^\pm_{w_\pm} \hra (R^+ \circledast R^-)/I_{\bfw}$
induce embeddings $R^\pm_{w_\pm} \hra \wh{R}_{\bfw}$.
We denote the images of 
$S^\pm_{w_\pm}$ in $\wh{R}_{\bfw}$ by the same symbols. 
Following Joseph \cite[\S 10.3.2]{J}, define
\[
S_{\bfw} = S_{w_+}^+ S_{w_-}^-.
\]
By \eqref{commR}, $S_{w_+}^+ S_{w_-}^- = S_{w_-}^- S_{w_+}^+$. 
More precisely, \eqref{n3} and \eqref{commR} imply the 
following commutation relation between 
the elements of $S^+_{w_+}$ and $S^-_{w_-}$. 
In terms of the identifications \eqref{ident} and the 
projections $g^\pm_{w_\pm}$ from \eqref{g}, we have
\begin{align}
\label{commRR}
&\big[ (c^-_{w_-, \la_1})^{-1} g_{w_-}^-(\xi_1) \big] 
\big[  g_{w_+}^+(\xi_2) (c^+_{w_+, \la_2})^{-1} \big]
\\
\nn
= &q^{ - \lcor \nu_1 + w_- \la_1 , \nu_2 - w_+ \la_2 \rcor} 
\big[  g_{w_+}^+(\xi_2) (c^+_{w_+, \la_2})^{-1} \big]
\big[ (c^-_{w_-, \la_1})^{-1} g_{w_-}^-(\xi_1) \big]
\\
\nn
+ & \sum_{\ga \in Q^+, \ga \neq 0}
\sum_{i=1}^{m(\ga)}
q^{ - \lcor \nu_1 + \ga + w_- \la_1, \nu_2 - \ga - w_+ \la_2 \rcor } 
\big[ g_{w_+}^+( S^{-1}(u_{\ga, i}) \xi_2) (c^+_{w_+, \la_2})^{-1} \big] .
\\
\nn
& \hspace{6.7cm}
\big[ (c^-_{w_-, \la_1})^{-1} g_{w_-}^-( S^{-1}(u_{-\ga, i}) \xi_1) \big],
\end{align}
for all 
$\la_i \in P^+$, $\nu_i \in P$, $i=1,2$,
$\xi_1 \in (V(- w_0 \la_1)^*)_{-\nu_1}$, 
$\xi_2 \in (V(\la_2)^*)_{-\nu_2}$.
It follows from \eqref{R+-} that the multiplication in 
$S_{\bfw}$ induces the vector space isomorphism
\begin{equation}
\label{Sisom}
S_{w_+}^+ \otimes_\KK S_{w_-}^- \congto S_{\bfw}.
\end{equation}

The algebra $R^+ \circledast R^-$ inherits a canonical 
$P \times P$-grading from the 
$P \times P$-gradings \eqref{Rgrade} of $R^\pm$. 
This induces a $P \times P$-grading on $R_{\bfw}$ and $S_{\bfw}$. 
For $\ga \in P$, there exists $\la \in P^+$ such that
$(V_{w_\pm}^\pm(\la))_{\pm w_\pm(\la) + \ga} \neq 0$, if and only if 
\[
\pm \ga \in \sum_{\be \in \De^+ \cap w_\pm (\De^-)} \Nset \be. 
\]
For $\xi^\pm \in ((V_{w_\pm}^\pm(\la))^*)_{\mp w_\pm(\la) - \ga}$, 
\begin{equation}
\label{Sdeg}
(c_{w_\pm, \la}^\pm)^{-1} \xi^\pm \in 
(S_{w_\pm}^\pm)_{- \ga, 0} 
\end{equation}
in terms of the identifications \eqref{ident}. Therefore:
\begin{equation}
\label{gradS+-}
(S^{\pm}_{w_\pm})_{\ga, 0} \neq 0, \; \;
\forall \ga \in \mp 
\sum_{\be \in \De^+ \cap w_\pm (\De^-)} \Nset \be
\; \; \; 
\mbox{and} \; \; \; 
(S^{\pm}_{w_\pm})_{\nu, \mu}=0, \; \; 
\mbox{otherwise}. 
\end{equation}
This also follows from \eqref{non0U} and the 
(anti)isomorphisms in \thref{isom}.
Eq. \eqref{gradS+-} implies that
\begin{multline}
\label{gradS}
(S_{\bfw})_{\ga, 0} \neq 0, \; \;
\forall \ga \in 
- \sum_{\be \in \De^+ \cap w_+ (\De^-)} \Nset \be
+ \sum_{\be \in \De^+ \cap w_- (\De^-)} \Nset \be
\\
\mbox{and} \; \; \; 
(S_{\bfw})_{\nu, \mu}=0, \; \; 
\mbox{otherwise}. 
\end{multline}
In \cite[Theorem 3.6]{Y2} we proved that the algebras $S_{w_\pm}^\pm$ 
play the role of quantum Schubert cells in relation to the 
$H$-spectra of quantum partial flag varieties. In a forthcoming
publication we will prove that the algebras 
$R^+ \circledast R^-$ and $S_{\bfw}$ are 
closely related to the quantizations of the standard Poisson structure 
on the double flag variety \cite{WY} and its restrictions to 
double Schubert cells, and will study the spectra of 
related double versions of the De Concini--Kac--Procesi algebras.

Denote by $\wh{L}^\pm_{w_\pm}$ and $\wh{L}_{\bfw}$ the subalgebras of 
$\wh{R}_{\bfw}$ generated by $(E^\pm_{w_\pm})^{\pm1}$ 
and $(\wh{E}_{\bfw})^{ \pm 1}$, respectively.
The algebras $\wh{L}^+_{w_+}$ and $\wh{L}^-_{w_-}$ 
are $r$ dimensional Laurent polynomial 
algebras over $\KK$ with generators 
$(c^+_{w_+, \om_i})^{\pm 1}$, $i=1, \ldots, r$
and 
$(c^-_{w_-, \om_i})^{\pm 1}$, $i=1, \ldots, r$.
We have the algebra isomorphism
\cite[10.3.2(2)]{J},
\begin{equation}
S_{\bfw} \# \wh{L}_{\bfw} \congto \wh{R}_{\bfw},
\label{SLR}
\end{equation}
where the smash product is computed using the actions
\begin{equation}
\label{smash}
c^\pm_{w_\pm, \la} \cdot u = q^{\lcor w_\pm \la, \nu \rcor} x, 
\quad \mbox{for} \; \;  
u \in (S_{\bfw})_{-\nu, 0}, \nu \in Q,
\end{equation}
because of \eqref{n3} and \eqref{commR}.

For $\la \in P^+$ choose an identification $V(\la)^* \cong V(- w_0 \la)$ 
normalized so that $\xi_{1, \la}^+ \mt v_{- \la}$ 
in terms of the lowest weight vectors fixed in \S \ref{2.3}
and the vectors $\xi^+_{w, \la}$ defined in \S \ref{2.4}. 
Let $\{ \xi_i \}$ 
and $\{ v_i \}$ be two sets of dual weight vectors of 
$V(\la)^*$ and $V(\la)$.
Define
\begin{equation}
\label{x}
x_\la = \sum_i c^\la_{\xi_i, \la} c^{-w_0 \la}_{v_i, - \la},
\end{equation}
where in the second term we used the identification 
$V(-w_0 \la)^* \cong V(\la)^{**} \cong V(\la)$. Then  
$x_\la \in Z(R^+ \circledast R^-)$, see \cite[Lemma 9.1.12]{J}. The images 
of $x_\la$ in $\wh{R}_{\bfw}$ will be denoted by 
the same symbols. Denote by $\wh{E}$ the multiplicative 
subset of $\wh{R}_{\bfw}$ generated by $x_{\om_i}$, 
$i=1, \ldots, r$ and by 
$\wh{L}$ the $\KK$-subalgebra of 
$\wh{R}_{\bfw}[\wh{E}^{-1}]$ generated by $x_{\om_i}^{\pm 1}$, 
$i=1, \ldots, r$. Denote 
\begin{equation}
\label{y}
y_{\om_i} = (c_{w_+, \om_i}^+)^{-1} (c_{w_-, \om_i}^-)^{-1} 
x_{\om_i} \in S_{\bfw}, \; \; i =1, \ldots, r.
\end{equation}
Continuing \eqref{SLR}, we have \cite[10.3.2(4)]{J}, 
\[
\wh{L} \otimes_\KK (S_{\bfw}[y^{-1}_{\om_i},  i=1, \ldots, r] \# \wh{L}^-_{w_-}) \congto
\wh{R}_{\bfw} [\wh{E}^{-1}].
\]
Joseph proved \cite[\S 9.2.4]{J} that 
the evaluation map $x_{\om_i} \mt 1$, 
$i=1, \ldots, r$ (i.e. 
$y_{\om_i} \mt (c_{w, \om_i}^+)^{-1} (c_{w, \om_i}^-)^{-1}$)
induces a surjective homomorphism $\wh{R}_{\bfw} [{\wh{E}}^{-1}] \to R_{\bfw}$,
from which he obtained the algebra isomorphism \cite[10.3.2(5)]{J},
\begin{equation}
\label{ySR}
\psi_{\bfw} \colon 
S_{\bfw}[y^{-1}_{\om_i},  i=1, \ldots, r] \# 
\wh{L}^-_{w_-} \congto R_{\bfw}.
\end{equation}
\subsection{Homogeneous $P$-normal elements of the algebras $S^\pm_{w_\pm}$}
\label{3.5}
Our proof of \thref{bfw-center} is based upon a study of 
a special kind of normal elements of the 
algebras $S^\pm_{w_\pm}$ and $S_{\bfw}$. These normal 
elements commute with the elements of the algebras 
$S^\pm_{w_\pm}$ and $S_{\bfw}$ up to an automorphism 
coming from the action \eqref{Tract0} of $\Tset^r$, 
restricted to a subgroup of $\Tset^r$ isomorphic
to the weight lattice $P$.

\bde{qnormal} We say that $z_\pm \in S^\pm_{w_\pm}$ is a $P$-normal 
element if there exists $\de_\pm \in P$ such that 
\[
z_\pm s = q^{\lcor \de_\pm, \ga \rcor } s z_\pm, \quad
\forall s \in (S^\pm_{w_\pm})_{-\ga,0}, \ga \in Q.
\] 
Analogously, we say that $z \in S_{\bfw}$ is a 
$P$-normal element if there exists $\de \in P$ such that 
\[
z s = q^{\lcor \de, \ga \rcor } s z, \quad
\forall s \in (S_{\bfw})_{-\ga,0}, \ga \in Q.
\] 
\ede

The motivation for the above definition is as follows. 
The abelian group  $P$ acts on $S_{\bfw}$ by 
\[
\mu \cdot s = q^{ \lcor \mu, \ga \rcor} s, \; \; 
\mbox{for} \; \; s \in (S_{\bfw})_{-\ga,0}, \mu \in P
\]
and preserves its subalgebras $S^\pm_{w_\pm}$. (It
is easy to see that the action \eqref{Tract0} of 
$\Tset^r$ on $R_q[G]$ induces an action on $S_{\bfw}$.
The above action is a restriction of this action to a subgroup
of $\Tset^r$ isomorphic to $P$.) An element
$z \in S_{\bfw}$ is $P$-normal, if it is a normal element and it 
commutes with the elements of $S_{\bfw}$ via an algebra 
automorphism coming from the $P$-action:
\[
z s = (\de \cdot s) z, \quad \forall s \in S_{\bfw},
\]
for some $\de \in P$.
The same applies to the subalgebras $S^\pm_{w_\pm}$.

\bre{qnormal1}
\leref{Iw} and \eqref{gradS+-} imply that the $\Zset$-span 
of all roots $\ga \in Q$ such that $(S^\pm_{w_\pm})_{-\ga,0} \neq 0$ is 
$Q_{\SS(w_\pm)}$. Thus in \deref{qnormal} one can assume 
that $\de_\pm \in P_{\SS(w_\pm)}$. Analogously 
\eqref{gradS} and \leref{Iw} imply that the $\Zset$-span 
of all $\ga \in Q$ such that $(S_{\bfw})_{-\ga,0} \neq 0$ is 
$Q_{\SS({\bfw})}$. Therefore in \deref{qnormal} one can assume 
that $\de \in P_{\SS({\bfw})}$.
\ere

For $\la \in P^+$ denote
\begin{equation}
\label{dd}
d_{w_\pm, \la}^\pm = (c_{w_\pm, \la}^\pm)^{-1} c_{1, \la}^\pm
\in (S^\pm_{w_\pm})_{ \pm (w_\pm- 1)\la, 0}. 
\end{equation}
These elements are $P$-normal; applying \leref{comm}, we 
obtain
\begin{equation}
\label{ddcomm}
d_{w_\pm, \la}^\pm s = q^{- \lcor (w_\pm + 1 ) \la , \ga \rcor } 
s d_{w_\pm, \la}^\pm, \quad \forall s \in (S_{w_\pm}^\pm)_{-\ga, 0}.
\end{equation}
For all $\la_1, \la_2 \in P^+$, 
\begin{equation}
\label{multd}
d_{w_\pm, \la_1}^\pm d_{w_\pm, \la_2}^\pm = 
q^{\pm \lcor \la_1, (w_\pm - 1) \la_2 \rcor }
d_{w_\pm, \la_1 + \la_2}^\pm.
\end{equation}
One verifies this using \eqref{n3} 
and \eqref{normalization}:
\begin{align*}
&d_{w_\pm, \la_1}^\pm d_{w_\pm, \la_2}^\pm = 
(c_{w_\pm, \la_1}^\pm)^{-1} c_{1, \la_1}^\pm
(c_{w_\pm, \la_2}^\pm)^{-1} c_{1, \la_2}^\pm 
\\
= &q^{\pm \lcor \la_1, (w_\pm -1) \la_2 \rcor }
(c_{w_\pm, \la_1}^\pm)^{-1} (c_{w_\pm, \la_2}^\pm)^{-1} 
c_{1, \la_1}^\pm c_{1, \la_2}^\pm \\
= & q^{\lcor \pm \la_1, (w_\pm -1) \la_2 \rcor }
(c_{w_\pm, \la_1 + \la_2 }^\pm)^{-1} c_{1, \la_1+ \la_2}^\pm = 
q^{\lcor \pm \la_1, (w_\pm -1) \la_2 \rcor }
d_{w_\pm, \la_1 + \la_2}^\pm. 
\end{align*}

The following result relates the degrees of the homogeneous 
$P$-normal elements of the algebras $S_{w_\pm}^\pm$ 
and the weights $\delta_\pm$ in \deref{qnormal}.

\bth{nS1} Assume that $\KK$ is an arbitrary base field and 
$q \in \KK^*$ is not a root of unity.
Let $z_\pm \in (S^\pm_{w_\pm})_{ \nu_\pm, 0}$ be a 
homogeneous $P$-normal element. Then there exists 
$\eta_\pm \in P_{\SS(w_\pm)}$ such that 
$\nu_\pm = \pm (w_\pm - 1) \eta_\pm$ and 
\[
z_\pm s = q^{- \lcor (w_\pm +1) \eta_\pm, \ga \rcor } s z_\pm, 
\quad \forall s \in (S^\pm_{w_\pm})_{-\ga,0}, \ga \in Q.
\]
\eth

Caldero determined \cite{Ca} the set of normal elements of $\UU_+$
with very different methods, using the Joseph--Letzter results \cite{JL}.
In the special 
case of $w_\pm = w_0$ (where $w_0$ is the longest element of $W$),
\thref{nS1} follows from \cite{Ca}. In Section \ref{normal}, building 
upon \thref{nS1} and other results, we will prove
that every homogeneous normal element of $S^\pm_{w_\pm}$ is $P$-normal
and eventually show that all homogeneous normal elements 
of $S^\pm_{w_\pm}$ are scalar multiples of $d_{w_\pm, \la}^\pm$
for $\la \in P^+$.
Those results are postponed to a later section, since they 
require various intermediate steps. 

For the proof of \thref{nS1} we will need the following lemma.

\ble{nS2} Assume that $z_\pm \in (S^\pm_{w_\pm})_{ \nu_\pm, 0}$ is a 
homogeneous $P$-normal element such that 
\begin{equation}
\label{zz}
z_\pm s = q^{\lcor \de_\pm, \ga \rcor } s z_\pm, \quad
\forall s \in (S^\pm_{w_\pm})_{-\ga,0}, \ga \in Q,
\end{equation}
for some $\de_\pm \in P$. Then for all $i \in \II({\bfw})$,
\[
\lcor \nu_\pm + \de_\pm, \al_i\spcheck \rcor \; \;  
\mbox{and} \; \; 
\lcor \nu_\pm - \de_\pm, \al_i\spcheck \rcor
\]
are even integers.
\ele 
\begin{proof} 
Fix a reduced expression $w_\pm= s_{j_1} \ldots s_{j_l}$.
Denote by $\be_1, \ldots, \be_l$ the roots \eqref{beta} and 
by $X^\pm_{\be_1}, \ldots, X^\pm_{\be_l}$ the root 
vectors \eqref{rootv}. Recall 
the graded (anti)isomor\-phisms $\varphi_w^\pm \colon S^\pm_w \to \UU^w_\mp$
from \thref{isom}.
We have 
\begin{equation}
\label{Xcomm}
\varphi_{w_\pm}^\pm(z_\pm) X_{\be_j}^\mp = 
q^{- \lcor \de_\pm, \be_j \rcor} X_{\be_j}^\mp
\varphi_{w_\pm}^\pm(z_\pm), \quad \forall j = 1, \ldots l.
\end{equation}
Recall the notation \eqref{monomial} and the 
notion of highest term of a nonzero element of 
$\UU^w_\pm$, defined in \S 2.7.  
Let $p (X^\mp)^{\bf{n}}$, $p \in \KK^*$, 
${\bf{n}} \in \Nset^{ \times l}$
be the highest term 
of $\varphi_{w_\pm}^\pm(z_\pm)$. Since $\UU^w_\pm$ 
are $Q$-graded algebras, we have
\begin{equation}
\label{nu}
\nu_\pm = \mp \sum_{i=1}^l n_i \be_i.
\end{equation}
Applying \eqref{LS0} we obtain that 
for $j =1, \ldots, l$
\[
\varphi_{w_\pm}^\pm(z_\pm) X_{\be_j}^\mp  - 
p q^{ \lcor \sum_{i=1}^{j-1} n_i \be_i, \be_j \rcor } 
(X^\mp)^{(n_1, \ldots, n_j +1, \ldots n_l)}
\]
and
\[
X_{\be_j}^\mp \varphi_{w_\pm}^\pm(z_\pm)
- p q^{ \lcor \sum_{i=j+1}^l n_i \be_i, \be_j \rcor } 
(X^\mp)^{(n_1, \ldots, n_j +1, \ldots n_l)}
\]
belong to $\Span \{ (X^\pm)^{{\bf{n}}'} \mid {\bf{n}}' < 
(n_1, \ldots, n_j +1, \ldots n_l) \}$.
Comparing this with \eqref{Xcomm} leads to 
\[
- \lcor \de_\pm , \be_j \rcor = \sum_{i=1}^{j-1} n_i \lcor \be_i, \be_j \rcor
- \sum_{i=j+1}^l n_i \lcor \be_i, \be_j \rcor, \quad j =1, \ldots, l.
\]
Now \eqref{nu} implies 
\[
- \lcor \de_\pm \pm \nu_\pm , \be_j \rcor =  
2 \sum_{i=1}^{j-1} n_i \lcor \be_i, \be_j \rcor
+ \lcor \be_j , \be_j \rcor.
\]
Hence
\begin{align*}
- \lcor \de_\pm \pm \nu_\pm, \be_j\spcheck \rcor &= 
2 \sum_{i=1}^{j-1} n_i \lcor \be_i, \be_j\spcheck \rcor
+ \lcor \be_j , \be_j \spcheck \rcor \\
&=  
2 \sum_{i=1}^{j-1} n_i \lcor \be_i, \be_j\spcheck \rcor
+ 2.
\end{align*}
is even for $j=1, \ldots, l$. Part (ii) of \leref{Iw}
implies that $\lcor \de_\pm \pm \nu_\pm, \al_i\spcheck \rcor$  
is even for all $i \in \SS(w_\pm)$. Therefore
\[
\lcor \de_\pm \mp \nu_\pm, \al_i\spcheck \rcor =  
\lcor \de_\pm \pm \nu_\pm, \al_i\spcheck \rcor \mp
2 \lcor \nu_\pm, \al_i\spcheck \rcor
\]
is also even for all $i \in \SS(w_\pm)$.
\end{proof}
\noindent
{\em{Proof of \thref{nS1}.}} Assume that $\de_\pm \in P_{\SS(w_\pm)}$ 
is such that \eqref{zz} holds, recall \reref{qnormal1}. Using 
\eqref{zz} and \eqref{ddcomm} we obtain
\begin{align*}
z_\pm d_{w_\pm, \la}^\pm &= 
q^{- \lcor \de_\pm, \pm(w_\pm-1) \la \rcor } 
d_{w_\pm, \la}^\pm z_\pm
\\
&=q^{ - \lcor \nu_\pm, (w_\pm + 1 ) \la \rcor }
d_{w_\pm, \la}^\pm z_\pm, \quad \forall \la \in P^+.
\end{align*}
Since $S^\pm_{w_\pm}$ is a domain and 
$z_\pm, d_{w_\pm, \la}^\pm \neq 0$,
\[
\lcor \de_\pm, \pm (w_\pm - 1) \la \rcor =
\lcor \nu_\pm, (w_\pm  + 1 ) \la \rcor, \quad \forall \la \in P^+,
\]
i.e.
\[
\lcor \nu_\pm \mp \de_\pm, w_\pm (\la) \rcor 
= - \lcor \nu_\pm \pm \de_\pm, \la \rcor
= - \lcor w_\pm( \nu_\pm \pm \de_\pm), w_\pm(\la) \rcor,
\quad \forall \la \in P^+.
\]
Therefore
\[
w_\pm(\nu_\pm \pm \de_\pm) + (\nu_\pm \mp \de_\pm) = 0.
\]
So,
\begin{equation}
\label{lin}
(w_\pm + 1) \nu_\pm = \mp (w_\pm -1) \de_\pm.
\end{equation}
Decompose 
\begin{equation}
\h = \h_\pm^{(1)} \oplus \h_\pm^{(-1)} \oplus \h_\pm^{(c)},
\label{hdecomp}
\end{equation}
where $\h_\pm^{(1)}$, $\h_\pm^{(-1)}$ are 
the eigenspaces of $w_\pm$ with eigenvalues $1$, $-1$,
and $\h_\pm^{(c)}$ is the direct sum of the other eigenspaces 
of $w_\pm$. Denote by $\nu_\pm^{(1)}$, $\nu_\pm^{(-1)}$, 
$\nu_\pm^{(c)}$ and $\de_\pm^{(1)}$, $\de_\pm^{(-1)}$, 
$\de_\pm^{(c)}$ the components of $\nu_\pm$ and $\de_\pm$ 
in the decomposition \eqref{hdecomp}. Then \eqref{lin} 
implies that $\nu_\pm^{(1)} = 0$, $\de_\pm^{(-1)}= 0$ and 
$(w_\pm + 1) \nu_\pm^{(c)} = \mp (w_\pm -1) \de_\pm^{(c)}$.
Therefore
\[
\wt{\eta}_\pm = - [ \de_\pm^{(1)}/2 \pm \nu_\pm^{(-1)}/2 + 
(w_\pm+1)^{-1} \de_\pm^{(c)} ]  
\]
satisfies
\begin{align*}
\nu_\pm &= \pm(w_\pm-1) \wt{\eta}_\pm,
\\
\de_\pm &= - (w_\pm+1) \wt{\eta}_\pm.
\end{align*}
We have $\wt{\eta}_\pm= - (\de_\pm \pm \nu_\pm)/2 \in (1/2)P$. 
Let
\[
\wt{\eta}_\pm = \eta_\pm + \ol{\eta}_\pm, 
\quad \mbox{where} \; \; 
\eta_\pm \in (1/2) P_{\SS(w_\pm)}, 
\ol{\eta}_\pm \in (1/2) P_{\II(w_\pm)}.
\] 
\leref{nS2} implies that $\eta_\pm \in P_{\SS(w_\pm)}$.
Since $\ol{\eta}_\pm \in \ker (w_\pm -1)$, 
\begin{equation}
\label{nupm}
\nu_\pm = \pm(w_\pm-1) \eta_\pm.
\end{equation}
Moreover 
\[
\de_\pm + (w_\pm+1) \eta_\pm = -(w_\pm +1) \ol{\eta}_\pm = 
- 2 \ol{\eta}_\pm 
\]
belongs to $P_{\II(w_\pm)}$ and is thus orthogonal to all 
$\ga \in Q$ such that $(S^\pm_{w_\pm})_{\ga, 0} \neq 0$, 
because of \eqref{Sdeg} and \leref{Iw}. Hence \eqref{zz} 
implies 
\[
z_\pm s = q^{\lcor - w_\pm(\eta_\pm) - \eta_\pm, \ga \rcor } s z_\pm, 
\quad \forall s \in (S^\pm_{w_\pm})_{-\ga,0}, \ga \in Q.
\]
This equation and \eqref{nupm} establish the
statement of \thref{nS1}.  
\qed
\subsection{Homogeneous $P$-normal elements of the algebras $S_{\bfw}$}
\label{3.6}
We proceed with establishing certain properties 
of the homogeneous $P$-normal elements of $S_{\bfw}$,
which are similar to the ones in \thref{nS1} for 
the algebras $S^\pm_{w_\pm}$.

\bth{nS3} Assume that $\KK$ is an arbitrary base field, 
and $q \in \KK^*$ is not a root of unity.
Let $z \in (S_{\bfw})_{\nu, 0}$ be a 
homogeneous $P$-normal element. Then there exists 
$\eta \in P_{\SS({\bfw})}$ such that 
$\nu = (w_+ - w_-)\eta$ and 
\begin{equation}
\label{dede}
z s = q^{\lcor - (w_+ + w_-)\eta, \ga \rcor } s z, \quad
\forall s \in (S_{\bfw})_{-\ga,0}, \ga \in Q_{\SS({\bfw})}.
\end{equation}    
\eth
\begin{proof} 
Let $\de \in P$ be such that 
\begin{equation}
\label{zs}
z s = q^{\lcor \de, \ga \rcor } s z, \quad
\forall s \in (S_{w_\pm})_{-\ga,0}, \ga \in Q.
\end{equation}
For $\tau \in Q^+$ denote 
\begin{equation}
\label{S+ord}
(S^+_{w_+})_{>- \tau,0} = \bigoplus_{\tau' \in Q^+, \tau' < \tau}
(S^+_{w_+})_{-\tau',0}
\end{equation}
and
\begin{equation}
\label{S-ord}
(S^-_{w_-})_{< \tau,0} = \bigoplus_{\tau' \in Q^+, \tau' < \tau}
(S^-_{w_-})_{\tau',0}
\end{equation}
in terms of the partial order \eqref{po}.
Eq. \eqref{commRR} implies that for all $\tau \in Q^+$, 
$S^+_{w_+}(S^-_{w_-})_{< \tau,0} = (S^-_{w_-})_{< \tau,0}S^+_{w_+}$
and 
$(S^+_{w_+})_{>- \tau,0} S^-_{w_-} = S^-_{w_-} (S^+_{w_+})_{>- \tau,0}$.
We have
\begin{equation}
\label{Sdec}
(S_{\bfw})_{\nu, 0} = \bigoplus_{\tau \in Q^+} (S^+_{w_+})_{-\tau, 0}
(S^-_{w_-})_{\nu+ \tau, 0}.
\end{equation}
Denote by $\tau_1, \ldots, \tau_m$ the set of maximal elements of 
the set consisting of those $\tau \in Q^+$  
for which $z$ has a nontrivial component in $(S^+_{w_+})_{-\tau, 0}
(S^-_{w_-})_{\nu+ \tau, 0}$, recall \eqref{Sisom}.
Denote the component of $z$ in 
$(S^+_{w_+})_{-\tau_i, 0} (S^-_{w_-})_{\nu+ \tau_i, 0}$ by
\[
z_i = \sum_{j=1}^{h(i)} z^+_{ij} z^-_{ij}, 
\]
where $z^+_{ij} \in (S^+_{w_+})_{-\tau_i, 0}$, 
$z^-_{ij} \in (S^-_{w_-})_{\nu+\tau_i, 0}$, for all 
$1 \leq i \leq m$, $1 \leq j \leq h(i)$ and for each 
$1 \leq i \leq m$
\begin{equation}
\label{lin-ind}
z^+_{i1}, \ldots, z^+_{i h(i)} \; \; \: 
\mbox{are linearly independent}. 
\end{equation}
Fix $s_- \in (S^-_{w_-})_{-\ga_-,0}$, for some 
\begin{equation} 
\label{ga-} 
\ga_- \in \sum_{\be \in \De^+ \cap w_-(\De^-)} 
\Nset \be,
\end{equation}
recall \eqref{gradS+-}. From \eqref{commRR} we obtain
\[
s_- z =
\sum_{ij} s_- z^+_{ij} z^-_{ij} = 
\sum_{ij} q^{- \lcor \tau_i, \ga_- \rcor } z^+_{ij} s_- z^-_{ij} 
\mod  \Big( \sum_i (S^+_{w_+})_{> - \tau_i, 0} S^-_{w_-} \Big),
\]
while \eqref{zs} implies
\[
s_- z = q^{- \lcor \de, \ga_- \rcor} z s_- = 
q^{- \lcor \de, \ga_- \rcor} \sum_{ij} z^+_{ij} z^-_{ij} s_- 
\mod  \Big( \sum_i (S^+_{w_+})_{> - \tau_i, 0} S^-_{w_-} \Big).
\] 
Applying \eqref{Sdec}, \eqref{lin-ind} and the fact that the 
multiplication map $S^+_{w_+} \otimes_\KK S^-_{w_-} \to S_{\bfw}$ 
is a vector space isomorphism, leads to 
\begin{equation}
\label{zi}
z^-_{ij} s_- = q^{\lcor \de - \tau_i, \ga_- \rcor} s_- z^-_{ij}, 
\quad \forall i, j.
\end{equation}
Therefore all $z^-_{ij}$ are homogeneous 
$P$-normal elements of $S^-_{w_-}$. 
\thref{nS1} implies that there exists $\eta_- \in P_{\SS(w_-)}$
such that $z^-_{11} \in (S^-_{w_-})_{- (w_- -1) \eta_-}$, i.e. 
\begin{equation}
\label{nueta1}
\nu + \tau_1 = - (w_- - 1) \eta_-
\end{equation}
and 
\[
q^{\lcor \de - \tau_i, \ga_- \rcor} =  
q^{- \lcor (w_- + 1 )\eta_-, \ga_- \rcor},
\]
for all $\ga_-$ as in \eqref{ga-}, recall \eqref{gradS+-}.
Taking into account \leref{Iw}, we obtain
\begin{equation}
\label{w1}
\de - \tau_1 + (w_- + 1) \eta_- \in P_{\II(w_-)}.
\end{equation}
Interchanging the roles of $S^+_{w_+}$ and $S^-_{w_-}$, 
we represent $z_i$ as in \eqref{zi} with $z_{ij}^\pm$ such that
\[
z^-_{i1}, \ldots, z^-_{i h(i)} \; \; \: 
\mbox{are linearly independent}, 
\]
instead of \eqref{lin-ind}. For all  $s_+ \in (S^+_{w_+})_{-\ga_+, 0}$
we obtain
\[
z s_+ =
\sum_{ij} z^+_{ij} z^-_{ij} s_+ = 
\sum_{ij} q^{\lcor \nu + \tau_i, \ga_+ \rcor } z^+_{ij} s_+ z^-_{ij} 
\mod  \Big( \sum_i S^+_{w_+} (S^-_{w_-})_{< \nu+\tau_i, 0} \Big)
\]
and 
\[
z s_+ = q^{\lcor \de, \ga_+ \rcor} s_+ z = 
q^{\lcor \de, \ga_+ \rcor} \sum_{ij} s_+ z^+_{ij} z^-_{ij} 
\mod  \Big( \sum_i S^+_{w_+} (S^-_{w_-})_{ <\nu+\tau_i, 0} \Big),
\] 
from \eqref{commRR} and \eqref{zs}, respectively. Therefore all
$z^-_{ij}$ are homogeneous $P$-normal elements of $S^-_{w_-}$ and
\[
z^+_{ij} s_+ =
q^{\lcor \de - \nu - \tau_i, \ga_+ \rcor } s_+ z^+_{ij}.
\]
Applying \thref{nS1}, we obtain that there exists $\eta_+ \in P_{\SS(w_+)}$
such that $z^+_{11} \in (S^+_{w_+})_{ (w_+ - 1) \eta_+}$, i.e. 
\begin{equation} 
\label{nueta2}
- \tau_1 = ( w_+ - 1 ) \eta_+
\end{equation}
and 
\[
q^{\lcor \de - \nu - \tau_1, \ga_+ \rcor} =  
q^{- \lcor ( w_+ + 1 ) \eta_+, \ga_+ \rcor}
\]
for all $\ga_+ \in Q_{\SS(w_+)}$, recall \leref{Iw} and \eqref{gradS+-}.
The latter is equivalent to
\begin{equation}
\label{w2}
\de - \nu - \tau_1 + ( w_+ + 1 ) \eta_+ \in P_{\II(w_+)}.
\end{equation}

Adding \eqref{nueta1} and \eqref{nueta2} gives
\begin{equation}
\label{conc1}
\nu = (w_+ -1)\eta_+ - (w_- - 1)\eta_-.
\end{equation}
Combining \eqref{w1} and \eqref{nueta2} leads to
\begin{equation}
\label{conc2}
\de + (w_+ -1) \eta_+ + (w_- + 1) \eta_- \in P_{\II(w_-)}.
\end{equation}
Similarly \eqref{w2} and \eqref{nueta1} imply
\begin{equation}
\label{conc3}
\de + (w_+ + 1)\eta_+ + (w_- - 1) \eta_- \in P_{\II(w_+)}.
\end{equation}
Decompose
\[
\eta_\pm = \wh{\eta}_\pm + \wt{\eta}_\pm,
\]
so that 
\begin{align}
\wh{\eta}_+ \in P_{\SS(w_+) \cap \SS(w_-)}, \; \; 
\wt{\eta}_+ \in P_{\SS(w_+) \backslash (\SS(w_+) \cap \SS(w_-))} = 
P_{\SS(w_+) \cap \II(w_-)}, \\
\wh{\eta}_- \in P_{\SS(w_+) \cap \SS(w_-)}, \; \;
\wt{\eta}_- \in P_{\SS(w_-) \backslash (\SS(w_+) \cap \SS(w_-))} =
P_{\SS(w_-) \cap \II(w_+)}.
\end{align}
In particular,
\begin{equation}
\label{kerr}
(w_\pm-1) \wt{\eta}_\mp = 0.
\end{equation}
Subtracting the left hand sides of \eqref{conc2} and \eqref{conc3}, 
shows that
\[
2(\eta_+ - \eta_-) \perp Q_{\SS(w_+) \cap \SS(w_-)}.
\]
Therefore 
$\wh{\eta}_+ = \wh{\eta}_-$. Denote
\[
\eta = \wh{\eta}_+  +  \wt{\eta}_+  +  \wt{\eta}_-
= \eta_+ + \wt{\eta}_- = \eta_- + \wt{\eta}_+.
\]
From \eqref{conc1} we have
\begin{multline}
\label{nuX}
\nu = (w_+ -1)(\eta - \wt{\eta}_-) - 
(w_- - 1)(\eta - \wt{\eta}_+) \\ 
= (w_+ -1) \eta - (w_- -1) \eta = (w_+ - w_-) \eta,
\end{multline}
because of \eqref{kerr}. Eqs. \eqref{conc2} and \eqref{kerr} 
imply
\[
\de + (w_+ - 1)(\eta - \wt{\eta}_-) + 
(w_- + 1)(\eta - \wt{\eta}_+) = 
\de + (w_+ + w_-) \eta - 2 \wt{\eta}_+ 
\in P_{\II(w_-)},
\]
so
\[ 
\de + (w_+ + w_-) \eta \in P_{\II(w_-)}.
\]
Analogously \eqref{conc3} and \eqref{kerr} imply
\[
\de + (w_+ + w_-) \eta \in P_{\II(w_+)},
\]
i.e.
\[
\de + (w_+ + w_-) \eta \in P_{\II(w_+)} \cap P_{\II(w_-)} = 
P_{\II(\bfw)}.
\]

From \eqref{zs} we obtain that $\eta$ satisfies \eqref{dede}. Since 
it also satisfies \eqref{nuX}, it provides the needed weight for
the theorem.
\end{proof}
\subsection{Proof of \thref{bfw-center}}
\label{3.7}
Denote 
\begin{equation}
(R_{\bfw})_\mu = \bigoplus_{\nu \in P} (R_{\bfw})_{\nu, \mu}. 
\end{equation}
Recall that $Z_{\bfw} = Z (R_{\bfw})$ and denote
\[
(Z_{\bfw})_{\nu, \mu} = Z_{\bfw} \cap (R_{\bfw})_{\nu, \mu},
\quad
(Z_{\bfw})_\mu = Z_{\bfw} \cap (R_{\bfw})_\mu, \quad
\forall \nu, \mu \in P.
\]
Obviously 
\[
Z_{\bfw} = \bigoplus_{\nu, \mu \in P} (Z_{\bfw})_{\nu, \mu}.
\]
We will need the following theorem of Joseph and Hodges--Levasseur--Toro. 

\bth{Joseph2} (Joseph \cite[Theorem 8.11]{J1}, 
Hodges-Levasseur--Toro \cite[Theorem 4.14 (3)]{HLT}) 
For all $\mu \in P$, 
\[
\dim (Z_{\bfw})_\mu = 0 \; \; \mbox{or} \; \; 1.
\]
\eth
\noindent
Similarly to \thref{J-thm} the proof of this result 
in \cite{J1,HLT} only uses the assumption that 
$q \in \KK^*$ is not a root of unity, without 
restrictions on the characteristic of $\KK$.
 
Denote by $A_{\bfw}$ the subalgebra of $R_{\bfw}$ generated by
\[
\{ (c^+_{w, \om_i})^{\pm 1} \mid i \in \II({\bfw}) \} \cup
\{a_1^{\pm 1}, \ldots, (a_k)^{\pm 1} \},
\]
recall \S \ref{3.1}.
Since each of the generators of $A_{\bfw}$ is 
$P \times P$ homogeneous, 
\[
A_{\bfw} = \bigoplus_{\nu, \mu \in P} 
(A_{\bfw})_{\nu, \mu},
\quad \mbox{where} \; \; 
(A_{\bfw})_{\nu, \mu} = A_{\bfw} \cap 
(R_{\bfw})_{\nu, \mu}.
\]
Define
\[
(A_{\bfw})_\mu = \bigoplus_{\nu \in P} (A_{\bfw})_{\nu, \mu}.
\]
Because $\{\la^{(1)}, \ldots, \la^{(k)} \} \cup 
\{ \om_i \mid i \in \II({\bfw}) \}$ 
is a linearly independent set (recall \S \ref{3.1}), 
the monomials 
\[
\prod_{ i \in \II({\bfw})} (c^+_{w_+, \om_i})^{n_i} 
\prod_{j=1}^k a_j^{m_j} \in (A_{\bfw})_\nu,
\quad
\nu = \sum_{i \in \II({\bfw})} n_i \om_i + 2 \sum_{j=1}^k m_j \la^{(j)}
\]
are linearly independent for different
$(n_i \mid i \in \II({\bfw})) \in \Zset^{ \times |\II({\bfw})|}$,
$(m_1, \ldots, m_k) \in \Zset^{\times k}$.
Therefore
\begin{equation}
\label{Aw}
A_{\bfw} \cong \KK[ (c^+_{w, \om_i})^{\pm 1}, a_j^{\pm 1}, 
i \in \II({\bfw}), j = 1, \ldots k].
\end{equation}

Recall \eqref{Ibfw}, \eqref{Sbfw}, \eqref{L} and \eqref{Lred},
and denote
\begin{equation}
\label{LL}
\LL({\bfw}) = 
2 \wt{\LL}_{\red}({\bfw}) \bigoplus 
\left( \bigoplus_{i \in \II({\bfw})} \Zset \om_i \right) = 
2 \wt{\LL}({\bfw}) + P_{\II({\bfw})}.
\end{equation}
Since $P_{\II({\bfw})} \subseteq \LL({\bfw})$,
\[
2 \wt{\LL}({\bfw}) \subset \LL({\bfw}) 
\subset \wt{\LL}({\bfw})
\]
and
\[
\LL({\bfw})/2 \wt{\LL}({\bfw}) \cong \Zset_2^{\times |\II({\bfw})| }.
\]

We have
\[
a_j \in (A_{\bfw})_{-(w_+ - w_-)\la^{(j)}, 2 \la^{(j)}}, \quad j =1, \ldots, k,
\]
cf. \eqref{a} and
\[
c_{w, \om_j}^+ \in (A_{\bfw})_{\om_j, \om_j}, \quad
\forall j \in \II({\bfw}),
\]
which leads to:
\ble{A} For all $\mu \in \LL({\bfw})$,  
$\dim (A_{\bfw})_\mu =1$ and for all $\mu \notin \LL({\bfw})$, 
$(A_{\bfw})_\mu = 0$.
\ele
\bre{set}
Joseph \cite{J1} and Hodges--Levasseur--Toro \cite{HLT}
that the set of all $\mu \in P$ such that 
$(Z_{\bfw})_{\mu} \neq 0$ contains $2 \wt{\LL}({\bfw})$ and is contained 
in $\wt{\LL}({\bfw})$. \thref{bfw-center} determines explicitly this set;
it is equal to $\LL({\bfw})$.
\ere

\noindent
{\em{Proof of \thref{bfw-center}.}} By \eqref{cent1}
and \prref{cent}, $A_{\bfw}$ is a subalgebra of 
$Z_{\bfw}$. We need to prove that $Z_{\bfw} = A_{\bfw}$.
Let $\nu', \mu \in P$. We will prove that 
\begin{equation}
\label{d}
(Z_{\bfw})_{\nu', \mu} \neq 0
\end{equation}
forces 
\begin{equation}
\label{inc}
\mu \in \LL({\bfw}). 
\end{equation}
Then we can apply \thref{Joseph2} and \leref{A} to deduce
that $(Z_{\bfw})_\mu = (A_{\bfw})_\mu$, $\forall \mu \in P$.  
Therefore $Z_{\bfw} = A_{\bfw}$. 
  
We are left with showing that \eqref{d} implies 
\eqref{inc}. Fix $\nu', \mu \in P$ and 
\[
d \in (Z_{\bfw})_{\nu', \mu}, \; \; d \neq 0.
\]
The isomorphism \eqref{ySR} and eq. \eqref{gradS} imply 
that
\[
d = \psi_{\bfw}( u \# (c^-_{w_-, \mu})^{-1}), \quad
\mbox{for some} \; \; u \in 
(S_{\bfw}[y^{-1}_{\om_i},  i=1, \ldots, r])_{\nu'+w_-(\mu), 0}.
\]
For $\la = \sum_{i=1}^r n_i \om_i \in P^+$ write
\begin{equation}
\label{yla}
y_\la = (y_{\om_1})^{n_1} \ldots (y_{\om_r})^{n_r}.
\end{equation}
From \eqref{x} and \eqref{y} we have 
\begin{equation}
\label{yc}
y_\la \in (S_{\bfw})_{(w_+ - w_-)\la, 0}.
\end{equation}
Let $u = z y_\la^{-1}$ for some $\la \in P^+$ and 
\begin{equation}
\label{z}
z \in (S_{\bfw})_{\nu, 0}, \; \; z \neq 0,
\end{equation}
where $\nu = \nu' - (w_+- w_-)(\la) + w_-(\mu)$. 
Thus
\begin{equation}
\label{ce}
\psi_{\bfw} \big( (z y_\la^{-1}) \# (c^-_{w_-, \mu})^{-1} \big) \in 
Z(R_{\bfw})_{\nu + (w_+ -w_-)(\la) - w_-(\mu), \mu} \; \; 
\mbox{with} \; \; z \neq 0. 
\end{equation}
In particular, it commutes with 
$c^\pm_{w_\pm, \mu'}$, for all $\mu' \in P$.
Using \eqref{n3}, we obtain
\[
\lcor w_\pm \mu', \nu + (w_+ -w_-)\la - w_-(\mu) \rcor +
\lcor w_\pm \mu', w_\pm \mu \rcor = 0, \quad
\forall \mu' \in P.
\]
Therefore 
\[
\nu + (w_+ -w_-)\la - w_-(\mu) = - w_+ (\mu) = - w_- (\mu),
\]
i.e.
\begin{equation}
\label{eq1}
\nu = -(w_+ - w_-) \la \; \; 
\mbox{and} \; \; \mu \in \ker(w_+ - w_-).
\end{equation}

Since $x_{\om_i} \in Z (R^+ \circledast R^-)$, \eqref{yla} and 
\eqref{n3} imply
\begin{equation}
\label{yla-norm}
y_\la s' = q^{\lcor -(w_+ + w_-)\la, \ga' \rcor} s' y_\la, 
\quad \forall s' \in (S_{\bfw})_{-\ga', 0}. 
\end{equation}
Because of \eqref{smash} one has
\begin{equation}
\label{sm2}
(1 \# c^-_{w_-, \mu}) (s' \# 1) = 
q^{\lcor w_-(\mu), \ga' \rcor} (s' \# 1) (1 \# c^-_{w_-, \mu} ), 
\; \forall s' \in (S_{\bfw})_{-\ga', 0}, 
\ga' \in Q_{\SS({\bfw})}.
\end{equation}
From \eqref{ce}, \eqref{yla-norm} and \eqref{sm2} it follows that
\begin{equation}
\label{zcommut}
z s' = q^{\lcor - (w_+ + w_-)\la + w_-(\mu), \ga' \rcor} 
s' z, \quad \forall s' \in (S_{\bfw})_{-\ga', 0}, \ga' \in Q_{\SS({\bfw})},
\end{equation}
recall \eqref{gradS} and \leref{Iw}.
In particular, $z \in (S_{\bfw})_{\nu, 0}$ is a homogeneous 
$P$-normal element.
\thref{nS3} implies that there exists $\eta \in P_{\SS({\bfw})}$ such that 
\begin{equation}
\label{ker2}
\nu = (w_+ - w_-)\eta
\end{equation}
and 
\begin{equation}
\label{zcommut2}
z s' = q^{\lcor - (w_+ + w_-)\eta, \ga' \rcor}
s' z, \quad \forall s' \in (S_{\bfw})_{-\ga', 0},
\ga' \in Q_{\SS({\bfw})}.
\end{equation}
Comparing \eqref{eq1} and \eqref{ker2}, gives that 
$\la - \eta \in \ker (w_+ - w_-)$. Therefore
\begin{equation}
\label{la-eta}
\la - \eta \in \wt{\LL}({\bfw}).
\end{equation}
Combining \eqref{zcommut} and \eqref{zcommut2} implies
that
\[
w_-(\mu) - (w_+ + w_-)(\la - \eta) = w_-(\mu) - 2 w_-(\la - \eta) 
\in P_{\II({\bfw})}.
\]
Thus 
\[
\mu - 2 (\la - \eta) \in P_{\II({\bfw})},
\]
because each element of $P_{\II({\bfw})}$ is fixed under $w_-^{-1}$. Finally 
this, together with \eqref{la-eta}, leads to
\[
\mu \in \wt{\LL}({\bfw}) + P_{\II({\bfw})} = \LL({\bfw}).
\] 
Therefore \eqref{d} implies \eqref{inc}, which completes 
the proof of \thref{bfw-center}
\qed
\\

\thref{bfw-center} makes Joseph's description of prime ideals 
of $R_q[G]$ more explicit. In parts (ii) and (iii) of \thref{J-thm}
one can replace $Z_{\bfw}$ with the explicit Laurent 
polynomial ring $A_{\bfw}$ given by \eqref{Aw}. 

\bco{Zco1} Assume that $\KK$ is an arbitrary base field,
and $q \in \KK^*$ is not a root of unity. For 
${\bfw} \in W \times W$ and  
$J^0 \in \Spec A_{\bfw}$ define
\[
\iota_{\bfw}(J^0)=\{r \in R_q[G] \mid (r + I_{\bfw}) \in R_{\bfw} J^0 \}.  
\]
Then $\iota_{\bfw}(J^0) \in \Spec_{\bfw} R_q[G]$ and 
\[
\iota_{\bfw} \colon \Spec A_{\bfw} \to \Spec_{\bfw} R_q[G]
\]
is a homeomorphism for all ${\bfw} \in W \times W$.
Moreover $\iota_{\bfw}$ restricts to a 
homeomorphism from $\Max A_{\bfw}$ to 
$\Prim_{\bfw} R_q[G]$. 
\eco

The application of this result to the primitive spectrum
of $R_q[G]$, described in \thref{prim}, is the starting point 
for explicitly relating $\Prim R_q[G]$ to the symplectic foliation 
of the underlying Poisson Lie group, discussed 
in the next section.
\sectionnew{Primitive ideals of $R_q[G]$ and a Dixmier map for $R_q[G]$}
\lb{Dixmier}
\subsection{A formula for the primitive ideals of $R_q[G]$}
\label{4.1}
When the base field $\KK$ is algebraically closed, the results from 
the previous section lead to an explicit parametrization 
of $\Prim R_q[G]$ and to a more explicit formula for the 
primitive ideals of $R_q[G]$ than the previously known ones, 
which is in turn used in Section \ref{Max} to classify $\Max R_q[G]$.
Based on this formula, we explicitly determine the stabilizers of the 
primitive ideals of $R_q[G]$ under the $\Tset^r \times \Tset^r$-action 
obtained by combining the actions \eqref{Tract} and \eqref{Tract0}.
This was not possible with the previously known formulas. 
In light of \thref{J-thm} (iii), we obtain the exact structure 
of $\Prim_{\bfw} R_q[G]$ is a $\Tset^r \times \Tset^r$-homogeneous space.
For $\KK= \Cset$, we combine this with the Kogan--Zelevinsky results \cite{KZ} 
to construct a $\Tset^r \times \Tset^r$-equivariant map from the symplectic 
foliation of the corresponding Poisson Lie group to $\Prim R_q[G]$. In this paper 
we use the term Dixmier type map in the wide sense, referring to a map from 
the topological space of 
the symplectic foliation associated with the semiclassical limit of an 
algebra $R$ to $\Prim R$, which is expected to be a homeomorphism.

Throughout the section the base field $\KK$ will be assumed to be algebraically 
closed. Recall the setting of \S \ref{3.1}. For ${\bfw}=(w_+, w_-)$, we fix a 
basis $\la^{(1)}, \ldots, \la^{(k)}$ of $\wt{\LL}_{\red}({\bfw})$, 
where $k = \dim \ker(w_+ - w_-) - |\II({\bfw})|$, recall \eqref{Lred}.
Represent
\[
\la^{(j)} = \la^{(j)}_+ - \la^{(j)}_-, 
\]
for some $\la^{(j)}_+$ and $\la^{(j)}_-$, 
which belong to $P^+$ and have disjoint support,
cf. \eqref{supp}. For $\zeta_j \in \KK$ define
\begin{equation}
\label{bzeta}
b_j(\zeta_j) = c^+_{w_+, \la^{(j)}_+} c^-_{w_-, \la^{(j)}_-} - 
\zeta_j c^+_{w_+, \la^{(j)}_-} c^-_{w_-, \la^{(j)}_+},
\quad j = 1, \ldots, k.
\end{equation}
Then
\begin{align}
\label{ab}
a_j - \zeta_j &= c^+_{w_+, \la^{(j)}} (c^-_{w_-, \la^{(j)}})^{-1} - \zeta_j
\\
\nn
&= (c^+_{w_+, \la^{(j)}_-})^{-1} c^+_{w_+, \la^{(j)}_+} 
             c^-_{w_-, \la^{(j)}_-} (c^-_{w_-, \la^{(j)}_+})^{-1} - \zeta_j
\\
\nn
&= (c^+_{w_+, \la^{(j)}_-})^{-1} b_j(\zeta_j) (c^-_{w_-, \la^{(j)}_+})^{-1},
\end{align}  
recall \eqref{a}. Thus 
$b_j(\zeta_j) = c^+_{w_+, \la^{(j)}_-} (a_j - \zeta_j) c^-_{w_-, \la^{(j)}_+}$.
Using \eqref{n3} and the fact that $a_j \in R_{\bfw}$ 
are central elements, we 
obtain that $b_j(\zeta_j) \in R_q[G]/I_w$ are normal:
\begin{equation}
\label{bcomm}
b_j(\zeta_j) c = q^{\lcor w_+(\la^{(j)}_-) +w_-(\la^{(j)}_+), \mu \rcor - 
\lcor \la^{(j)}_+ + \la^{(j)}_-, \nu \rcor }c b_j(\zeta_j), \; \; 
\forall c \in (R_q[G]/I_{\bfw})_{-\nu, \mu}, \nu, \mu \in P.
\end{equation}
For $\zeta = (\zeta_1, \ldots, \zeta_k) \in (\KK^*)^{\times k}$
and 
$\theta = \{ \theta_i \}_{i \in \II({\bfw}) } \in 
(\KK^*)^{\times |\II({\bfw})|}$
denote
\begin{equation}
\label{Jw0}
J_{\bfw, \zeta, \theta} = \iota_{\bfw} \Big( \sum_{j=1}^k R_{\bfw} (a_j - \zeta_j) + 
\sum_{i \in \II({\bfw})} R_{\bfw}(c^+_{w_+, \om_i} - \theta_i) \Big).
\end{equation}
Eq. \eqref{bcomm} implies that
\begin{multline}
\label{Jw}
J_{\bfw, \zeta, \theta} =
\big\{ r \in R_q[G] \; \big| \; c r \in  
\sum_{j=1}^k R_q[G] b_j(\zeta_j) 
\\
+ \sum_{i \in \II({\bfw})} R_q[G](c^+_{w_+, \om_i} - \theta_i) 
+ I_{\bfw} \; \mbox{for some} \; c \in E_{\bfw} \big\},
\end{multline}
recall \eqref{Ew}.
\thref{J-thm} (iii) and \thref{bfw-center} lead to the following result,
cf. \coref{Zco1}.

\bth{prim} Assume that $\KK$ is an algebraically closed field and 
$q \in \KK^*$ is not a root of unity. Then for all 
${\bfw} = (w_+, w_-) \in W \times W$, the stratum of primitive 
ideals $\Prim_{\bfw} R_q[G]$ consists of the ideals 
$J_{\bfw, \zeta, \theta}$ given by \eqref{Jw}, where 
$\zeta = (\zeta_1, \ldots, \zeta_k) \in (\KK^*)^{\times k}$,
$\theta = \{ \theta_i \}_{i \in \II({\bfw}) } \in 
(\KK^*)^{\times |\II({\bfw})|}$ and 
$k = \dim \ker(w_+ - w_-) - |\II({\bfw})|$.
\eth

This result plays a key role in our classification 
of the maximal ideals of $R_q[G]$ in Section \ref{Max}.

The cases of $\g= {\mathfrak{sl}}_2$ and $\g = \mathfrak{sl}_3$ 
of \thref{prim} were obtained
by Hodges--Levasseur \cite{HL0} and Goodearl--Lenagan \cite{GLen0}, 
respectively, who also proved a stron\-ger result without the $E_{\bfw}$
localization in \eqref{Jw}. 
Their methods are very different from ours and use in an 
essential way the low rank of the underlying Lie algebra.
\subsection{Structure of $\Prim_{\bfw} R_q[G]$ as a 
$\Tset^r \times \Tset^r$-homogeneous space}
\label{4.2}

The commuting $\Tset^r$-actions \eqref{Tract} and \eqref{Tract0} 
on $R_q[G]$ can be combined to the following rational 
$\Tset^r \times \Tset^r$-action by $\KK$-algebra automorphisms:
\begin{equation}
\label{Tract2}
(t',t) \cdot c = (t')^\nu t^{\mu} c, \quad
t', t \in \Tset^r, c \in R_q[G]_{-\nu, \mu}, \nu, \mu \in P. 
\end{equation}
We obtain induced $\Tset^r \times \Tset^r$-actions 
on $R_{\bfw}$, $Z_{\bfw} = A_{\bfw}$, $\Spec A_{\bfw}$, 
$\Spec_{\bfw} R_q[G]$. 

Denote by $\Stab_{\Tset^r}(.)$ and $\Stab_{\Tset^r \times \Tset^r} (.)$ 
the stabilizers with respect to the 
actions \eqref{Tract} and \eqref{Tract2}, 
respectively. The map 
$\iota_{\bfw} \colon \Spec A_{\bfw} \to \Spec_{\bfw} R_q[G]$
is $\Tset^r \times \Tset^r$-equivariant. In particular,
\[
\Stab_{\Tset^r \times \Tset^r} \iota_{\bfw} (J^0) = 
\Stab_{\Tset^r \times \Tset^r} (J^0), \quad
\forall J^0 \in \Spec A_{\bfw}.
\]
The equivariance of $\iota_{\bfw}$ and \eqref{Jw0} imply that
the $\Tset^r \times \Tset^r$-action on $\Prim_{\bfw} R_q[G]$ 
is given by
\begin{equation}
\label{tt1}
(t',t) \cdot J_{\bfw, \zeta, \theta} = 
J_{\bfw, (t', t) \cdot \zeta, (t',t) \cdot \theta}, 
\end{equation}
where
\begin{align}
\label{tt2}
(t', t) \cdot \{ \zeta_j \}_{j=1}^k &= 
\{ (t')^{-(w_+-w_-)\la^{(j)}  } t^{- 2 \la^{(j)} } \zeta_j \}_{j=1}^k
\; \; \mbox{and}
\\
\label{tt3}
(t', t) \cdot \{ \theta_i \}_{i \in \II({\bfw})} &= 
\{ (t_i t'_i)^{-1} \theta_i \}_{i \in \II({\bfw})},
\end{align}  
because
\[
(t', t) \cdot a_j = (t')^{(w_+-w_-)\la^{(j)}} t^{2 \la^{(j)}}, 
j=1, \ldots, k, 
\; \; 
(t', t) \cdot c^+_{w_+, \om_i}= t'_i t_i c^+_{w_+, \om_i}, 
i \in \II({\bfw}).
\]
This implies the following result describing 
the stabilizers of $J_{\bfw, \zeta, \theta}$ under 
the action \eqref{Tract2} of $\Tset^r \times \Tset^r$ and in particular 
under the action \eqref{Tract} of $\Tset^r$.

\bpr{st} If $\KK$ is algebraically closed and $q \in \KK^*$ is not a root 
of unity, then for all ${\bfw} = (w_+, w_-) \in W \times W$, 
$\zeta = (\zeta_1, \ldots, \zeta_k) \in (\KK^*)^{\times k}$,
$\theta = \{ \theta_i \}_{i \in \II({\bfw}) } \in 
(\KK^*)^{\times |\II({\bfw})|}$:
\begin{multline}
\label{stabb}
\Stab_{\Tset^r \times \Tset^r} J_{\bfw, \zeta, \theta} = \{ 
(t', t ) \in \Tset^r \times \Tset^r \mid 
\\
t^{2 \la} = (t')^{-(w_+-w_-)\la}, \;  \forall 
\la \in \wt{\LL}_{\red}({\bfw}), \; \; t_i = (t'_i)^{-1}, \; 
\forall i \in \II({\bfw}) \},
\end{multline}
recall \eqref{Lred}. In particular, we have:
\begin{align}
\Stab_{\Tset^r} J_{\bfw, \zeta, \theta} 
\label{stab}
&= \{ t \in \Tset^r \mid t_i = 1, \forall i \in \II({\bfw}), 
\; t^{2 \la^{(j)}} = 1, \forall j = 1, \ldots, k \}
\\
\label{stab'}
&= \{ t \in \Tset^r \mid 
t^\la = 1, \; \forall \la \in \LL({\bfw}) \},
\end{align}
cf. \eqref{LL}.
\epr
\begin{proof}
Eq. \eqref{stabb} follows directly from \eqref{tt1}, \eqref{tt2}, 
and \eqref{tt3}. Eq. \eqref{stab} is the restriction of 
\eqref{stabb}. Eq. \eqref{stab'} is a consequence of
\eqref{LL} and \eqref{stab}.
\end{proof}
\subsection{The standard Poisson Lie structure on $G$ and its symplectic leaves}
\label{4.3}
In the remaining part of this section we assume that the base field 
is $\KK= \Cset$. The assumption on the deformation parameter $q \in \Cset^*$ 
will be that it is not a root of unity, as before. Thus $\g$ will be a complex 
simple Lie algebra. We will denote by $G$ the connected, simply connected 
algebraic group with Lie algebra $\g$. Let $B_\pm$ be a pair of opposite Borel 
subgroups of $G$, and $T = B_+ \cap B_-$ be the corresponding maximal torus of 
$G$. One has the isomorphism of complex tori:
\begin{equation}
\label{TrT}
T \cong \Tset^r, \quad \exp( \zeta_1 \al_1\spcheck + \ldots + \zeta_r \al_r\spcheck ) =
(\exp \zeta_1, \ldots, \exp \zeta_r).  
\end{equation}
Denote $\h = \Lie \, T$. Let $\lcor.,.\rcor$ be the nondegenerate invariant 
bilinear on $\g$ which matches the form \eqref{inner} on $\h^*$. For 
$\mu \in P$ define the characters $t^\mu$ of $T$ by 
\[
\exp(h)^\mu = \exp ( \lcor \mu, h \rcor).
\]
This matches \eqref{tmu} under the isomorphism \eqref{TrT}.
Additionally, denote $t^w = w^{-1} t w$ 
for $w \in W$, $t \in T$. For ${\bfw}=(w_+, w_-) \in W \times W $, 
set \cite[\S 2.4]{KZ}
\[
T^{\bfw} = \{ (t^{w_+})^{-1} t^{w_-} \mid t \in T \}.
\]

Let $\{e_\al\}$ and $\{f_\al\}$, $\al \in \De^+$ be sets of 
dual root vectors of $\g$, 
normalized by $\lcor e_\al, f_\al \rcor = 1$.  For $x \in \g$ denote by 
$L(x)$ and $R(x)$ the left and right invariant vector fields on $G$. The 
standard Poisson structure on $G$ is given by 
\[
\pi_G = \sum_{\al \in \De^+} L(e_\al) \wedge L(f_\al) - 
\sum_{\al \in \De^+} R(e_\al) \wedge R(f_\al).
\]

For $j=1, \ldots, r$ choose the representative 
\begin{equation}
\label{ols}
\ol{s}_j = \exp (e_{\al_j} ) \exp (-f_{\al_j}) \exp(e_{\al_j}) \in N_G(T)
\end{equation}
of $s_j \in W$,
where $N_G(T)$ denotes the normalizer of $T$ in $G$. This choice is 
slightly different from the one of Kogan and Zelevinsky \cite{KZ}, 
but we need it to match it to the braid group action \eqref{braid}.
For $w \in W$ choose a reduced expression 
$w= s_{j_1} \ldots s_{j_l}$ and define 
\[
\ol{w} = \ol{s}_{j_1} \ldots \ol{s}_{j_l} \in N_G(T).
\]
This choice of representative of a Weyl group element 
in the normalizer of the torus $T$ does not depend on the 
choice of the reduced expression, because the elements
\eqref{ols} satisfy the braid relations analogously to \cite{KZ}.

Denote the unipotent radicals of $B_\pm$ by $U_\pm$. 
We have $U_- T U_+ \cong U_- \times T \times U_+$ 
under the group product. For $g \in U_- T U_+$ denote 
its components in $U_-$, $T$ and $U_+$ by 
$[g]_-$, $[g]_0$ and $[g]_+$, respectively. 

The left and right regular actions of $T$ on $G$,
preserve $\pi_G$.
The $T$-orbits of symplectic leaves of $\pi_G$ 
(under any of those actions) are 
\cite{HL0} the double Bruhat cells 
$G^{\bfw}= G^{w_+, w_-} = B_+ w_+ B_+ \cap B_- w_- B_-$ of $G$, 
${\bfw}= (w_+, w_-) \in W \times W$. 
The symplectic leaves of $(G^{\bfw}, \pi_G)$
were determined by Kogan and Zelevinsky \cite[Theorem 2.3]{KZ}.
\bth{KZ} (Kogan--Zelevinsky \cite{KZ}) 
For every ${\bfw} = (w_+, w_-) \in W \times W$, the set 
\begin{multline}
\label{Sww}
\SLl_{\bfw} = \big\{ g \in G^{\bfw} \mid 
\big[ \ol{w_+}\, {}^{-1} g \big]_0 .
\big( \big[ g \ol{w_-^{-1}} \big] \big)^{w_-} \in 
T^{\bfw}, 
\\
\big[ \ol{w_+}\, {}^{-1} g \big]_0^{\om_i} = 1, \; 
\forall i \in \II({\bfw})
\big\}   
\end{multline}
is a symplectic leaf of $(G, \pi_G)$. All symplectic 
leaves of $G$ have the form $\SLl_{\bfw} . t$ for some 
$t \in T$, ${\bfw} \in W \times W$. 
\eth

In particular, the double Bruhat cell $G^{\bfw}$ is 
the $T$-orbit of the symplectic leaf $\SLl_{\bfw}$ under 
both the left and right $T$-actions.
\subsection{Equations for the symplectic leaves of $(G^{\bfw}, \pi_G)$}
\label{4.4}
Next, we make a minor reformulation of \thref{KZ} to match 
it to \thref{prim}. We have the following description 
of $T^{\bfw}$.

\ble{Tw} The torus $T^{\bfw}$ is given by
\[
T^{\bfw} = \{ t \in T \mid t^\mu = 0, \forall 
\mu \in \wt{\LL}({\bfw}) \},
\]
recall \eqref{L}.
\ele
\noindent
{\em{Sketch of the proof.}} For $\mu \in P$, $w \in W$ one has 
$(t^w)^\la = t^{w \la}$. Therefore
\[
T^{\bfw} \subseteq \{ t \in T \mid t^\mu = 0, \forall 
\mu \in \wt{\LL}({\bfw}) \},
\]
because $((t^{w_+})^{-1} t^{w_-})^\la= t^{(w_- - w_+)\la} = 1$ 
for all $\la \in \wt{\LL}({\bfw})$. It is clear that both sides
of the above inclusion are algebraic subgroups of $T$ of codimension 
$\dim \ker(w_+ - w_-)$. One easily checks that they are both 
connected, thus they coincide.
\qed
\\ \hfill \\
For $\la \in P^+$ denote by $\wt{V}(\la)$ the irreducible 
finite dimensional module of $G$ with highest weight $\la$. 
For $v \in \wt{V}(\la)$ and $\xi \in \wt{V}(\la)^*$ denote 
the matrix coefficient 
\[
\wt{c}^\la_{\xi, v} \in \Cset[G], \quad
\wt{c}^\la_{\xi, v} (g) = \xi( g v), \; g \in G.
\]
Let $v_\la \in V(\la)_\la$ and $\xi_\la \in V(\la)^*_{-\la}$,
be such that $\xi_\la(v_\la)=1$. Similarly let 
$v_{-\la} \in V(-w_0 \la)_{-\la}$ and 
$\xi_{-\la} \in V(- w_0 \la)^*_{\la}$, be such that $\xi_{-\la} (v_{-\la})=1$.
Analogously to the quantum case for $\la \in P^+$ and $w \in W$ define
\begin{equation}
\label{cwtdef}
\wt{c}^+_{w, \la} = \wt{c}^\la_{\ol{w} \xi_\la, v_\la}, \quad
\wt{c}^-_{w, \la} = \wt{c}^{-w_0 \la}_{ \left( \ol{w^{-1}} \right)^{-1} \xi_{- \la}, v_{-\la}}.
\end{equation}
Their key property is that
\begin{equation}
\label{c+-}
\wt{c}^+_{w, \la}(g) =
\big( \big[ \ol{w_+}\, {}^{-1} g \big]_0 \big)^\la,  \quad
\wt{c}^-_{w, \la}(g) = 
\big( \big[ g \ol{w_-^{-1}} \big] \big)^{- w_- \la},
\end{equation}
which is verified by a direct computation. This property is the 
reason for the above normalization of $\wt{c}^-_{w, \la}$.

We also have 
\begin{equation}
\label{prodc}
\wt{c}^+_{w, \la_1} \wt{c}^+_{w, \la_2} = \wt{c}^+_{w, \la_1 + \la_2}, 
\quad \forall \la_1, \la_2 \in P^+.
\end{equation}
For all $\la \in P^+$,
$\wt{c}^\pm_{w_\pm, \la}$ are regular functions on $G$ which are
nowhere vanishing on $G^{\bfw}$.
Fix $\la \in P$, represent it as $\la = \la_1 - \la_2$ for some 
$\la_1, \la_2 \in P^+$, and define 
\begin{equation}
\label{wtcP}
\wt{c}^\pm_{w_\pm, \la} = \wt{c}^\pm_{w_\pm, \la_1} (\wt{c}^\pm_{w_\pm, \la_2})^{-1},
\end{equation} 
considered as a rational function on $G$ and a regular function on $G^{\bfw}$. 
The definition \eqref{wtcP} does not depend on the choice of $\la_1$ and $\la_2$,
because of \eqref{prodc}. Eq. \eqref{c+-} holds for all $\la \in P$.

For $j = 1, \ldots, k$ denote
\[
\wt{a}_j = \wt{c}^+_{w_+, \la^{(j)}} (\wt{c}^-_{w_-, \la^{(j)}})^{-1}.
\]
\bco{sleaf} Let ${\bfw} = (w_+, w_-) \in W \times W$. 
Then the symplectic leaves 
of $(G, \pi_G)$ inside the double Bruhat cell $G^{\bfw}$ are 
parametrized by $(\Cset^*)^{ \times \dim \ker(w_+- w_-)}$. They are 
exactly the sets 
\begin{equation}
\label{Swzt}
\SLl_{\bfw, \zeta, \theta} = \{ 
g \in G^{\bfw} \mid \wt{a}_j(g) = \zeta_j, j=1, \ldots, k, \; 
\wt{c}^+_{w_+, \om_i}(g) = \theta_i, i \in \II({\bfw}) \},
\end{equation}
for $\zeta = (\zeta_1, \ldots, \zeta_k) \in (\Cset^*)^{\times k}$,
$\theta = \{ \theta_i \}_{i \in \II({\bfw}) } \in 
(\Cset^*)^{\times |\II({\bfw})|}$ and 
$k = \dim \ker(w_+ - w_-) - |\II({\bfw})|$.
\eco
\begin{proof} \leref{Tw} and \eqref{c+-} imply that 
$\SLl_{\bfw} = \SLl_{\bfw, \zeta, \theta}$ for 
$\zeta_j = \theta_i=1$, $\forall j=1, \ldots, k$, $i \in \II({\bfw})$.  
\thref{KZ} now implies the statement using the right 
regular action of $T$.
\end{proof}
\subsection{A $\Tset^r \times \Tset^r$-equivariant 
Dixmier map for $R_q[G]$}
\label{4.5}
Denote by $\Sympl (G, \pi_G)$ the symplectic foliation 
space of the Poisson structure $\pi_G$ 
(i.e the set of symplectic 
leaves with the induced topology from the Zariski 
topology on $G$). Define the Dixmier type map 
\[
D_G \colon \Sympl (G, \pi_G) \to \Prim R_q[G], \quad
D_G(\SLl_{\bfw, \zeta, \theta}) = J_{\bfw, \zeta, \theta},
\]
${\bfw} \in W \times W$, 
$\theta = \{ \theta_i \}_{i \in \II({\bfw}) } \in 
(\Cset^*)^{\times |\II({\bfw})|}$, 
$\zeta = (\zeta_1, \ldots, \zeta_k) \in (\Cset^*)^{\times k}$, 
where $k = \dim \ker(w_+ - w_-) - |\II({\bfw})|$.

Consider the $T\times T$-action on $G$ coming from the 
left and right regular actions 
\[
(t, t') \cdot g = (t')^{-1} g t^{-1}
\]
and transfer it to a $\Tset^r \times \Tset^r$-action
on $G$ via \eqref{TrT}. This action preserves $\pi_G$ 
and thus induces an action on $\Sympl (G, \pi_G)$.
(The choice of the inverses is made to match 
this action with the actions \eqref{Tract2} and \eqref{action} 
in the quantum situation.)
Analogously to \eqref{tt1} one shows that
\begin{equation}
\label{TG}
(t',t) \cdot \SLl_{\bfw, \zeta, \theta} = 
\SLl_{\bfw, (t', t) \cdot \zeta, (t',t) \cdot \theta}, 
\end{equation}
in terms of \eqref{tt2} and \eqref{tt3}. Combining
\thref{prim}, \thref{KZ}, \eqref{tt1} 
and \eqref{TG}, we obtain:
\bth{equivD} Assume that the base field $\KK$ is $\Cset$ and 
$q \in \Cset^*$ is not a root of unity. Then  
the Dixmier type map 
$D_G \colon \Sympl (G, \pi_G) \to \Prim R_q[G]$
is a $\Tset^r \times \Tset^r$-equivariant 
bijection. 
\eth

The original orbit method conjecture \cite{HL0} of 
Hodges and Levasseur for $R_q[G]$ can be formulated 
more precisely as follows:

\bcj{conj} Under the above assumptions, 
the Dixmier map
\[
D_G \colon \Sympl (G, \pi_G) \to \Prim R_q[G]
\]
is a homeomorphism.
\ecj

\bre{special} In the special case when the base field 
is $\KK= \Cset$ and $q$ is transcendental over $\Qset$,
one can prove that 
the elements $c^\pm_{w, \la}$  defined in \eqref{cd-def} 
specialize to the elements $\wt{c}^\pm_{w, \la}$
defined in \eqref{cwtdef} for all $\la \in P^+$, $w \in W$,
when $q$ is specialized to $1$. 
The elements $\wt{c}^\pm_{w, \la}$ are in turn related to the 
setting of Kogan and Zelevinsky \cite{KZ} via \eqref{c+-}. 
The normalization of the elements $c^\pm_{w, \la}$ in 
\S \ref{2.4} was made so that our setting matches the 
latter whenever specialization can be defined (i.e. 
when $\KK$ has characteristic $0$ 
and $q$ is transcendental over $\Qset$).

The special case of \thref{bfw-center} for base fields $\KK$ 
of characteristic $0$ and $q \in \KK$ transcendental over 
$\Qset$ can be proved in a simpler way using specialization
and the Kogan--Zelevinsky result \cite{KZ}. One should point 
out though that the results on $P$-normal elements of the algebras
$S^\pm_{w_\pm}$ and $S_{\bfw}$ which are the building blocks 
of the proof of \thref{bfw-center} play an important role 
throughout the rest of the paper.
\ere

A result of \cite{RY} proves that 
in the complex case the Haar functional on $R_q[G]$
is an integral of the traces of the irreducible $*$-representations
of $R_q[G]$ (with respect to the $*$-involution associated
to the compact form of $G$)
classified in \cite{LS}, see \cite[Theorem 5.2]{RY}
for details. Those representations correspond to particular 
primitive ideals in $\Prim_{(w,w)} R_q[G]$ for $w \in W$. 
We finish with raising the question whether irreducible representations 
corresponding to the other primitive ideals of $R_q[G]$ play 
any (noncommutative) differential geometric role.  
\sectionnew{Separation of variables for the algebras $S^\pm_w$}
\lb{free1}
\subsection{Statement of the freeness result}
\label{5a.1}
Recall that Joseph's isomorphism \eqref{ySR} represents 
the localizations $R_{\bfw}$, ${\bfw}=(w_+, w_-) \in W \times W$
in terms of the algebras $S^\pm_{w_\pm}$. In this and the next sections
we prove a number of results for the algebras $S^\pm_{w_\pm}$
which will play a key role in our study of $R_{\bfw}$ 
and $\Max R_q[G]$ in the following sections. These results 
also establish important properties of the 
De Concini--Kac--Procesi algebras via the (anti)isomorphisms 
from \thref{isom}.

Throughout this section we fix a Weyl group element $w \in W$. 
Denote by $N^\pm_w$ the subalgebras of $S^\pm_w$ generated 
by the normal elements $d^\pm_{w, \om_i}$, 
$i \in \SS(w)$, recall \eqref{dd}.
In this section we describe the structure
of the algebras $S^\pm_w$, considered as $N^\pm_w$-modules. 
We apply these results in several directions. 
In Section \ref{normal} we use them to classify 
all normal elements of the algebras $S^\pm_w$
and equivalently the De Concini--Kac--Procesi algebras
$\UU^w_\pm$. In fact, we prove that all homogeneous normal 
elements of the algebras $S^\pm_w$ are scalar multiples 
of $d^\pm_{w, \la}$, $\la \in P^+_{\SS(w)}$, and 
all normal elements of $S^\pm_w$ are equal to 
(certain) linear combinations of 
$d^\pm_{w, \la}$, $\la \in P^+_{\SS(w)}$.
As another
application in Section \ref{normal} we classify all 
prime elements of the algebras $S^\pm_w$. 
In Section \ref{Module} the results of this section are used 
to describe the structure of $R_{\bfw}$ as a module
over its subalgebra generated by the sets of normal elements 
$E_{\bfw}^{\pm 1}$, recall \eqref{Ew}. This is then applied 
to classify the maximal ideals of $R_q[G]$ in 
Section \ref{Max}.

We start by noting that \eqref{dd} implies  
$d^\pm_{w, \om_i} \in \KK^*$, for $i \in \II(w)$. 
Because of this, one only needs to consider $d^\pm_{w, \om_i}$
for $i \in \SS(w)$. It follows from \eqref{ddcomm}
that 
\[
d^\pm_{w, \la_1} d^\pm_{w, \la_1}
=
q^{ \pm( \lcor w_\pm(\la_1), \la_2 \rcor - 
\lcor \la_1, w_\pm(\la_2) \rcor ) }  
d^\pm_{w, \la_2} d^\pm_{w, \la_1}, 
\quad \forall \la_1, \la_2 \in P^+,
\]
and in particular
\begin{equation}
\label{dij}
d^\pm_{w, \om_i} d^\pm_{w, \om_j}
=
q^{ \pm( \lcor w_\pm(\om_i), \om_j \rcor -
\lcor \om_i, w_\pm(\om_j) \rcor ) }
d^\pm_{w, \om_j} d^\pm_{w, \om_i}, \quad
\forall i,j \in \SS(w).
\end{equation} 

The main result of the section is:

\bth{freeS} Let $\KK$ be an arbitrary base field, $q$ be an
element of $\KK^*$ which is not a root of unity, and
$w \in W$. Then:

(i) The algebra $N^\pm_w$ is isomorphic to 
the quantum affine space algebra over $\KK$ 
of dimension $|\SS(w)|$ with generators 
$d^\pm_{w, \om_i}$, $i \in \SS(w)$ and relations 
\eqref{dij}.

(ii) The algebra $S^\pm_w$ is a free left and right 
$N^\pm_w$-module in which $N^\pm_w$ is a direct summand, 
viewed as a module over itself.
\eth

An explicit form of the freeness result in the second
part of the theorem is obtained in \thref{freeS2} below.
The special case of $\g = {\mathfrak{sl}}_{r+1}$ and
$w = w_0$ in Theorems \ref{tfreeS} (ii) and \ref{tfreeS2}
is due to Lopes \cite{Lo}. 

In the next section we classify 
the normal elements of $S^\pm_w$. A
consequence of this result is that $N^\pm_w$ 
coincides with the subalgebra of $S^\pm_w$ generated by all 
of its homogeneous normal elements. In particular, 
$Z(S^\pm_w) \subset N^\pm_w$. Theorems of the above kind 
are motivated by the desire to extend the 
theorems for separation of variables of Kostant \cite{K} 
and Joseph--Letzter \cite{JL}
to quantized universal enveloping algebras of nilpotent 
Lie algebras. Kostant, and Joseph and Letzter
proved that $\UU(\g)$ and 
$\UU_q(\g)$ are free as modules over their centers
and deduced further important properties of the 
corresponding bases. In our case 
the centers of $S^\pm_w$ are in general too small compared
to the centers of $S^\pm_w[E^\pm_w]$, see \leref{Z}.
Thus one would obtain weaker results by considering 
the module structure of $S^\pm_w$ 
over their centers $Z(S^\pm_w)$ as opposite to the subalgebras
generated by the ``numerators'' and ``denominators''
of the central elements of $S^\pm_w[E^\pm_w]$.
It is the structure 
of $S^\pm_w$ as a module over the ``normal subalgebra''
$N^\pm_w$ that has applications to the structure
of $\Spec R_q[G]$ and $\Spec S^\pm_w$.
Two additional freeness results  will be obtained 
in Section \ref{Module} for the 
algebras $S_{\bfw}$ and $R_{\bfw}$.

We recall that a quantum affine space algebra is an algebra 
over $\KK$ with generators $X_1, \ldots, X_m$ 
and relations
\begin{equation}
\label{q-plane}
X_i X_j = p_{ij} X_j X_i, \quad i, j = 1, \ldots, m,
\end{equation}
for some $p_{ij} \in \KK^*$ such that 
$p_{ij} p_{ji} =1$, for all $i \neq j \in \{ 1, \ldots, m \}$,
$p_{ii} =1$, $i \in \{ 1, \ldots, m \}$.
Such an algebra has Gelfand--Kirillov dimension equal 
to $m$. It has a $\KK$-basis, consisting of the 
monomials
\begin{equation}
\label{mo}
(X_1)^{n_1} \ldots (X_m)^{n_m}, \quad 
n_1, \ldots, n_m \in \Nset.
\end{equation}
On the other hand, if a $\KK$-algebra is generated 
by some elements $X_1, \ldots, X_m$, which satisfy 
\eqref{q-plane} and the monomials \eqref{mo} are 
linearly independent, then the algebra is 
isomorphic to the above quantum affine space algebra.
We also recall that the localization of this 
algebra by the multiplicative subset generated by 
$X_1^{-1}, \ldots, X_m^{-1}$ is called 
quantum torus algebra.

Because of \eqref{multd}, the first part of \thref{freeS} 
essentially claims that the elements $d^\pm_{w, \la}$ 
are linearly independent
over $\KK$ for different $\la \in P^+_{\SS(w)}$. 
\subsection{Leading terms of the normal elements $\varphi_w^\pm(d^\pm_{w, \la})$}
\label{5a.2}
For the rest of this section we fix a reduced expression of $w$
\begin{equation}
\label{ii}
w = s_{i_1} \ldots s_{i_l},
\end{equation}
where $l = l(w)$ is the length of $w$. 
Denote this reduced expression by $\vec{w}$.
For $j \in \SS(w)$, let
\[
\Supp_j(\vec{w}) = \{ k = 1, \ldots, l \mid i_k = j \}.
\]

Recall the definition of the roots $\be_k$ (see \eqref{beta}) 
and the root vectors $X^\pm_{\be_k}$ (see \eqref{rootv}), 
$k = 1, \ldots, l$, associated to the reduced expression 
$\vec{w}$. Recall the definition \eqref{monomial} 
of the monomials $(X^\pm)^{\bf{n}}$, ${\bf{n}} \in \Nset^{\times l}$, 
the notions of leading term of an element of 
$\UU^w_\pm$ and degree of a monomial from \S 2.7.

For $j \in \SS(w)$ denote
\[
{\bf{e}}(\vec{w})_j = (n_{j1}, \ldots, n_{jl}) \in \Nset^{\times l},
\]
where 
\begin{equation}
\label{ew}
n_{jk} = 1 \; \; \mbox{if} 
\; \; k \in \Supp_j(\vec{w}), \; \; 
n_{jk} = 0 \; \; \mbox{if} 
\; \; k \notin \Supp_j(\vec{w}).
\end{equation}

Recall the (anti)isomorphisms $\varphi^\pm_w \colon S^\pm_w \to \UU^w_\mp$ 
from \thref{isom}. We will need the following fact for proof of \thref{freeS}.

\bpr{lead} Let $\KK$ be an arbitrary base field, $q$ be an 
element of $\KK^*$ which is not a root of unity, and 
$\vec{w}$ be a reduced expression of $w \in W$. Then for 
all $\la \in P^+$ the leading term of
$\varphi^\pm_w(d^\pm_{w, \la})$ 
has degree 
\[
(\lcor \la, \al_{i_1}\spcheck \rcor, \ldots, 
\lcor \la, \al_{i_l}\spcheck \rcor).
\]
In particular, for all $j \in \SS(w)$, the leading term of 
$\varphi^\pm_w(d^\pm_{w, \om_j})$ has degree 
${\bf{e}}(\vec{w})_j$.
\epr
\begin{proof}
We will prove the statement 
in the plus case. The minus case is 
analogous and is left to the reader.
Assume that the reduced expression $\vec{w}$ is given by \eqref{ii}.
Recall from \S \ref{2.3} that for $\la \in P^+$, 
$v_\la$ denotes a fixed highest weight vector
of $V(\la)$. Recall the definition of the vectors 
$\xi^+_{1, \la} \in V(\la)^*_{- \la}$ from \S \ref{2.4}.  
Taking into account the definition \eqref{map+} of the 
antiisomorphism $\varphi_w^+ \colon S^+_w \to \UU^w_-$,
and eqs. \eqref{n1}-\eqref{n2}, \eqref{Rwe} and \eqref{dd}, 
we see that the statement 
of the proposition is equivalent to:
\begin{equation}
\label{cond1}
\lcor \xi^+_{1, \la}, 
(\tau(X^+_{\be_1}))^{\lcor \la, \al_{i_1}\spcheck \rcor}
\ldots
(\tau(X^+_{\be_l}))^{\lcor \la, \al_{i_l}\spcheck \rcor} T_w v_\la \rcor \neq 0
\end{equation}
and
\begin{equation}
\label{cond2}
\lcor \xi^+_{1, \la}, (\tau((X^+)^{\bf{n}})) T_w v_\la \rcor \neq 0, 
\; {\bf{n}} \in \Nset^{ \times l} \; \; 
\Rightarrow \; \;  {\bf{n}} \leq (\lcor \la, \al_{i_1}\spcheck \rcor, \ldots,
\lcor \la, \al_{i_l}\spcheck \rcor)
\end{equation}
in the lexicographic order from \eqref{lexi}.

Since $\dim V(\la)_{s_i \la} = 1$, $i=1, \ldots, l$,
we have that $T_i^{-1} v_\la = p_i 
(X^-_{i_k})^{\lcor \la, \al_{i_k}\spcheck \rcor} v_\la$ 
for some $p_i \in \KK^*$. Because for all $i =1, \ldots, l$,
$v_\la$ is a highest weight vector for the 
$\UU_{q_i}({\mathfrak{sl}}_2)$-subalgebra of $\UU_q(\g)$ 
generated by $X^\pm_i$, $K^{\pm 1}_i$
with highest weight $\lcor \la, \al_i\spcheck \rcor \om_i$,
we have that 
$(X^+_i)^{\lcor \la, \al_i\spcheck \rcor}
(X^-_i)^{\lcor \la, \al_i\spcheck \rcor}  v_\la= 
p_i' v_\la$ for some $p_i' \in \KK^*$.
  
For $k=0, 1, \ldots, l$ denote 
$w_{(k)} = s_{i_1} \ldots s_{i_{k}}$. 
Using the above facts and eqs. \eqref{compat}, \eqref{tau-ident}, 
we obtain 
\begin{align}
\label{i1}
&(\tau(X^+_{\be_k}))^{\lcor \la, \al_{i_k}\spcheck \rcor}  
\big( T_{w(k)^{-1}}^{-1} v_\la \big)
\\
\nn
=& p_{i_k} \Big( T_{w(k-1)^{-1}}^{-1}(X^+_{i_k})^{\lcor \la, \al_{i_k}\spcheck \rcor} \Big)
\Big( T_{w(k-1)^{-1}}^{-1}\big((X^-_{i_k})^{\lcor \la, \al_{i_k}\spcheck \rcor} 
v_\la \big) \Big)
\\
\nn
= & p_{i_k}
T_{w(k-1)^{-1}}^{-1} \big( (X^+_{i_k})^{\lcor \la, \al_{i_k}\spcheck \rcor}
(X^-_{i_k})^{\lcor \la, \al_{i_k}\spcheck \rcor}  v_\la \big)
=  p_{i_k} p_{i_k}' T_{w(k-1)^{-1}}^{-1} v_\la.
\end{align}
Analogously one proves that for $m> 0$ 
\begin{equation}
\label{i2}
(\tau(X^+_{\be_k}))^{\lcor \la, \al_{i_k}\spcheck \rcor +m }  
\big( T_{w(k)^{-1}}^{-1} (v_\la) \big) = 0.
\end{equation}
 
Because $\dim V(\la)_{w \la} = 1$, $T_w v_\la = p T_{w^{-1}}^{-1} v_\la$
for some $p \in \KK^*$. Recursively applying \eqref{i1}, one obtains
\begin{align*}
&\lcor \xi^+_{1, \la}, 
(\tau(X^+_{\be_1}))^{\lcor \la, \al_{i_1}\spcheck \rcor}
\ldots
(\tau(X^+_{\be_l}))^{\lcor \la, \al_{i_l}\spcheck \rcor} T_w v_\la \rcor =
\\
= &p \lcor \xi^+_{1, \la}, 
(\tau(X^+_{\be_1}))^{\lcor \la, \al_{i_1}\spcheck \rcor}
\ldots
(\tau(X^+_{\be_l}))^{\lcor \la, \al_{i_l}\spcheck \rcor} T_{w^{-1}}^{-1} v_\la \rcor = \ldots
\\
= &p p_{i_l} p_{i_l}' \ldots p_{i_{k+1}} p_{i_{k+1}}'
\lcor \xi^+_{1, \la}, 
(\tau(X^+_{\be_1}))^{\lcor \la, \al_{i_1}\spcheck \rcor}
\ldots
(\tau(X^+_{\be_k}))^{\lcor \la, \al_{i_k}\spcheck \rcor} T_{w(k)^{-1}}^{-1} v_\la \rcor =
\ldots 
\\
= & p p_{i_l} p_{i_l}' \ldots p_1 p_1' \lcor  \xi^+_{1, \la}, v_\la \rcor \neq 0.
\end{align*}
This proves \eqref{cond1}. Assume that ${\bf{n}} \in \Nset^{ \times l}$ 
and ${\bf{n}} > (\lcor \la, \al_{i_1}\spcheck \rcor, \ldots,
\lcor \la, \al_{i_l}\spcheck \rcor)$. Then there exists 
$k \in [1,l]$ such that 
$n_j = \lcor \la, \al_{i_j}\spcheck \rcor$ for $j=k+1, \ldots, l$
and $n_k > \lcor \la, \al_{i_k}\spcheck \rcor$. 
Using \eqref{i1} and \eqref{i2}, one obtains
\begin{multline*}
\lcor \xi^+_{1, \la}, 
( \tau (X^+)^{\bf{n}} )  T_w v_\la \rcor = 
\lcor \xi^+_{1, \la},
p ( \tau (X^+)^{\bf{n}} )  T_{w^{-1}}^{-1} v_\la \rcor = 
\ldots 
\\
= p p_{i_l} p_{i_l}' \ldots p_{i_{k+1}} p_{i_{k+1}}' 
\lcor \xi^+_{1, \la}, (\tau(X^+_{\be_1}))^{n_1}
\ldots
(\tau(X^+_{\be_k}))^{n_k} T_{w(k)^{-1}}^{-1} v_\la \rcor=0.
\end{multline*}
This proves \eqref{cond2} and completes the proof 
of the Proposition.
\end{proof}
\subsection{Proof of \thref{freeS}}
\label{5a.3}
We begin with the proof of the first part of \thref{freeS}.
The second part of the theorem requires some additional 
facts. It is given at the end of the subsection.
\\ \hfill \\
\noindent
{\em{Proof of part (i) of \thref{freeS}.}} If $\la_1, \la_2 \in P^+_{\SS(w)}$
and $\la_1 \neq \la_2$, then there exists $j \in \SS(w)$ such that 
$\lcor \la_1, \al_j\spcheck \rcor \neq \lcor \la_2, \al_j\spcheck \rcor$.
\prref{lead} implies that all elements 
$\{ d^\pm_{w, \la} \}_{\la \in P^+_{\SS(w)} }$ have leading terms 
of different degrees. Therefore they are linearly independent
because of \thref{DKP}, which proves part (i) of \thref{freeS}. 
\qed

Denote the following two subsets of $\Nset^{\times l}$:
\[
\Sig(\vec{w}) = \bigoplus_{j \in \SS(w)} 
\Nset {\bf{e}}(\vec{w})_j
\]
and
\begin{equation}
\label{De}
\De(\vec{w}) = \{ (n_1, \ldots, n_l) \in \Nset^{\times l} \mid
\forall j \in \SS(w), \; \exists k \in \Supp_j(\vec{w}) \; \; 
\mbox{such that} \; \; n_{k} = 0 \}. 
\end{equation}
According to \prref{lead} the first subset consists of the degrees
of the leading terms of the elements 
$d^\pm_{w, \la} \in N_w^\pm$, $\la \in P_{\SS(w)}^+$. 
The following fact shows 
that the second subset is complementary to the first one. 
Its proof is left to the reader.

\ble{complement} Each element of $\Nset^{\times l}$ is
representable in a unique way as the sum of an element of 
$\Sig(\vec{w})$ and an element of $\De(\vec{w})$.
\ele

The second part of \thref{freeS} follows directly from the following 
theorem which provides an explicit presentation of $S^\pm_w$ as a
free $N^\pm_w$-module. 

\bth{freeS2} For an arbitrary base field $\KK$, $q \in \KK^*$ not 
a root of unity and a reduced expression $\vec{w}$ of 
$w \in W$:
\[
S^\pm_w = 
\bigoplus_{ {\bf{n}} \in \De(\vec{w}) } 
N^\pm_w \cdot (\varphi^\pm_w)^{-1} \big( (X^\mp)^{\bf{n}} \big) 
= 
\bigoplus_{ {\bf{n}} \in \De(\vec{w}) } 
(\varphi^\pm_w)^{-1} \big( (X^\mp)^{\bf{n}} \big)
\cdot N^\pm_w. 
\]
\eth

Note that ${\bf{0}} \in \De(\vec{w})$ and 
\[
N^\pm_w \cdot (\varphi^\pm_w)^{-1} \big( (X^\mp)^{\bf{0}} \big)
= N^\pm_w.
\]
\hfill \\
\noindent
{\em{Proof of \thref{freeS2}.}} It is sufficient to prove the first 
equality since the algebra $N^\pm_w$ is spanned by $P$-normal 
elements, which $q$-commute with 
$(\varphi^\pm_w)^{-1} \big( (X^\pm)^{\bf{n}} \big)$,
for ${\bf{n}} \in \Nset^{ \times l }$.

The theorem follows from the fact that the associated graded 
of $\UU^w_\pm$ with respect to the filtration induced from the ordering 
\eqref{lexi} is free over the associated graded of 
$\varphi_w^\pm (N^\pm_w)$, because of \prref{lead} and \leref{complement}. 
Here are the details. Recall the (anti)isomorphisms $\varphi^\pm_w \colon S^\pm_w \to \UU^w_\mp$
from \thref{isom}. To show
\begin{equation}
\label{sum}
S^\pm_w =
\sum_{ {\bf{n}} \in \De(\vec{w}) }
N^\pm_w \cdot (\varphi^\pm_w)^{-1} \big( (X^\mp)^{\bf{n}} \big),
\end{equation}
fix $s \in S^\pm_w$, $s \neq 0$. \prref{lead}, and 
Lemmas \ref{lLSmult} and \ref{lcomplement} 
imply that there exist $\la \in P^+_{\SS(w)}$,
${\bf{n}} \in \De(\vec{w})$, and $p \in \KK^*$, 
such that either  
$\varphi^\pm_w(s) - (X^\mp)^{\bf{n}} \varphi^\pm_w(d^\pm_{w, \la}) =0$,
or it is a nonzero element whose leading term has degree
strictly less than that of the leading term of $\varphi^\pm_w(s)$, recall 
\eqref{lexi}. Iterating this, gives \eqref{sum}.

Finally, the set
\[
\{ (X^\mp)^{\bf{n}} \varphi^\pm_w(d^\pm_{w, \la})  \mid 
\la \in P^+_{\SS(w)}, {\bf{n}} \in \De(\vec{w}) 
\}
\]
is linearly independent over $\KK$, because the 
elements of this set have leading terms of different degrees.
This follows from \prref{lead} and \leref{complement}.
\qed
\sectionnew{A classification of the normal  
and prime elements of the De Concini--Kac--Procesi algebras}
\lb{normal}
\subsection{Statement of the classification result}
\label{5.1}
In this section we develop further the line of argument of 
\S 3.5 and obtain a classification
of the sets of homogeneous normal elements of 
all De Concini--Kac--Procesi algebras 
$\UU^w_\pm$. Equivalently, 
this gives a classification of the homogeneous normal 
elements of the algebras $S^\pm_w$. 
We combine these results with the results
from the previous section to obtain an explicit 
description of the primitive ideals in 
the Goodearl--Letzter stratum \cite{GL} of 
$\Prim S^\pm_w$ over the $\{0 \}$ ideal.
These results are then applied 
to obtain a classification of all prime elements 
of the algebras $S^\pm_w$. At the end of the section, 
in \thref{allnorma} we prove that all normal elements of 
the algebras $S^\pm_w$ are equal to (certain) linear 
combinations of homogeneous normal elements. This 
produces an explicit classification of all normal elements of 
the algebras $S^\pm_w$ and $\UU^w_\pm$.

Our approach to the classification 
problem for the homogeneous normal elements of the algebras 
$S^\pm_w$ (which is the key step in the arguments 
in this section), is to prove first that each such element 
is $P$-normal, recall \deref{qnormal}. We then obtain the 
classification by an argument, which combines \thref{bfw-center}, 
\thref{freeS} on separation of variables for the algebras $S^\pm_w$, 
and a strong rationality result for $H$-primes of iterated skew 
polynomial extensions of Goodearl \cite{BG}.  

Throughout this section $w$ will denote a fixed 
element of the Weyl group $W$. Recall that the algebras
$\UU^w_\pm$ and $S^\pm_w$ are $Q_{\SS(w)}$-graded, 
by \eqref{non0U} and \eqref{gradS+-}. We call an element of these
algebras homogeneous, if it is homogeneous with 
respect to the corresponding grading.

Recall that an element $u$ of a noetherian domain 
$R$ is called prime if it is normal and $R u$ is a
height one prime ideal, which is completely prime. 
Recall the definition \eqref{dd} of the normal 
elements $d^\pm_{w, \la} \in (S^\pm_w)_{\pm(w-1)\la, 0}$, 
cf. \eqref{ddcomm}. The following theorem contains 
our classification result for homogeneous prime 
and homogeneous normal elements of the algebras 
$S^\pm_w$. In \thref{inhompr} below we obtain 
a classification of the inhomogeneous 
prime elements of the algebras $S^\pm_w$.

\bth{normal1} Assume that $\KK$ is an arbitrary base field and 
$q \in \KK^*$ is not a root of unity.
Let $w \in W$. Then:

(i) Every nonzero homogeneous normal element of 
$S^\pm_w$ is equal to an element of the form
\begin{equation}
\label{normla}
p d^\pm_{w, \la} \in (S^\pm_w)_{\pm(w-1) \la, 0}
\end{equation}  
for some $p \in \KK^*$, $\la \in P^+_{\SS(w)}$.
All such elements are distinct and even 
more the elements $d^\pm_{w, \la}$, $\la \in P^+_{\SS(w)}$
are linearly independent, cf. \thref{freeS}.

(ii) For all $i \in \SS(w)$, 
$d^\pm_{w, \om_i} \in (S^\pm_w)_{\pm(w-1)\om_i, 0}$
are pairwise nonproportional prime elements of 
$S^\pm_w$ and all homogeneous
prime elements of $S^\pm_w$ are 
nonzero scalar multiples of them.
\eth

Here and below, ``pairwise nonproportional elements'' means not 
a scalar multiple of each other. In \thref{allnorma} we obtain 
a further classification of all normal elements of the algebras 
$S^w_\pm$. Although that is an interesting extension 
of \thref{normal1}, it is the explicit classification 
of all homogeneous normal elements of $S^\pm_w$ that is needed for 
all applications. 

In view of \eqref{multd}, another 
way to formulate \thref{normal1} is to say that 
every nonzero homogeneous normal element of
$S^\pm_w$ is equal to an element of the form
\begin{equation}
\label{dprod}
p \prod_{i \in \SS(w)} (d^\pm_{w, \om_i})^{n_i}
\end{equation}
for some $n_1, \ldots, n_r \in \Nset$, 
$p \in \KK^*$. 
The elements $d^\pm_{w, \om_i}$ do not commute. In
\eqref{dprod} we take the product over $i$ in any 
fixed order. (Recall that the elements $d^\pm_{w, \om_i}$
$q$-commute.)

Recall the graded (anti)isomorphisms 
$\varphi^\pm_w \colon S^\pm_w \to \UU^w_\mp$
from \thref{isom}. 
We have the following reformulation of \thref{normal1},
which provides a classification of the 
sets of homogeneous prime elements and 
homogeneous normal elements 
of the De Concini--Kac--Procesi algebras
$\UU^w_\pm$.

\bth{normal2} (i) In the setting of \thref{normal1},
every nonzero homogeneous normal element of
$\UU^w_\pm$ is equal to an element
of the form
\[
p \varphi_w^\mp(d^\mp_{w, \la}) \in (\UU^w_\pm)_{\mp(w-1)\la}
\]
for some $p \in \KK^*$, $\la \in P^+_{\SS(w)}$.
All such elements are distinct and even linearly independent 
for different $\la$'s.

(ii) For all $i \in \SS(w)$,
$\varphi_w^\mp(d^\mp_{w, \om_i}) \in (\UU^w_\pm)_{\mp(w-1)\om_i}$
are pairwise nonproportional prime elements of $\UU^w_\pm$
and all homogeneous prime elements of $\UU^w_\pm$
are nonzero scalar multiples of them.
\eth

The case of $w= w_0$ of the first part of the theorem 
was proved by Caldero \cite{Ca}, using very different 
methods from ours, based on the Joseph--Letzter results \cite{JL}.
The case of the second part of the theorem for the algebras 
of quantum matrices is due to Launois, Lenagan and 
Rigal \cite[Proposition 4.2]{LLR}.

In the case when the characteristic of $\KK$ is 0 and $q$ 
is transcendental over $\Qset$, one can deduce part (ii) 
of \thref{normal2} from \cite[Theorem 1.1 (c)]{Y}. 
\bre{centers} The structure of the centers of the universal enveloping 
algebras $\UU(\n_\pm)$ was described by Joseph in \cite{J0} and Kostant 
in an unpublished work and \cite{K2}. They proved that these centers 
are polynomial algebras, described their generators, and obtained various 
other structure results. Theorems \ref{tnormal2} and \ref{tinhompr} can be 
considered 
as quantum counterparts of these results. These theorems imply that the 
subalgebras generated by all normal elements of the algebras 
$\UU^w_\pm$ are quantum affine space algebras in the generators
$\varphi_w^\mp(d^\mp_{w, \om_i})$, $i \in \SS(w)$
for all $w \in W$. We note that 
noncommutative rings which are not universal enveloping
algebras can have small centers but much bigger subalgebras generated 
by all normal elements. This is why the latter subalgebras exhibit 
richer structure. Secondly, quantum affine space algebras 
are the simplest analogs of polynomial rings in the class 
of noncommutative rings.   

It is very interesting that the subalgebras generated by all normal elements of the 
quantized universal enveloping algebras $\UU^w_\pm$ are better organized 
than the centers of the universal enveloping algebras $\UU(\n_\pm \cap w(\n_\mp))$.
In the former case we always have quantum affine space algebras 
by Theorems \ref{tnormal2} and \ref{tinhompr}. In the latter case Joseph and 
Hersant \cite[\S 8.5]{J0} showed that the centers of the universal enveloping
algebras of the nilradicals of certain parabolic subalgebras of $\g$ of
even type $A$ are not polynomial algebras.
In another respect a construction of Lipsman--Wolf \cite{LW} and 
Kostant \cite{K1} recovers a part of the center of $\UU(\n_\pm)$ in terms 
of matrix coefficients of finite dimensional $\g$-modules. At the same time 
all homogeneous normal elements of $\UU^w_\pm$ (even for an arbitrary
Weyl group element $w$) are given in terms of $R$-matrices and matrix 
coefficients of finite dimensional $\UU_q(\g)$-modules by \thref{normal2} (i).
This raises the question of the relation between the 
two constructions in the specialization $q=1$.
\ere
\subsection{Homogeneous normal and $P$-normal elements of $S^\pm_w$}
\label{5.2}
First, we show that 
each homogeneous normal element of the algebras 
$S^\pm_w$ is $P$-normal in the sense 
of \deref{qnormal}.

\bpr{normalp} All homogeneous normal elements of
$S^\pm_w$ are $P$-normal. 
\epr
The proof of this proposition will be given in \S \ref{5.4}.
For this proof we will need two results. The first concerns the
number of pairwise nonproportional prime elements of the algebras $\UU^w_\pm$
(proved in this subsection) and the second concerns a special 
kind of ``diagonal'' automorphisms of the algebras $\UU^w_\pm$
(proved in \S \ref{5.3}). 

A noetherian domain $R$ is said \cite{Ch} to be a 
unique factorization domain, if $R$ has at least 
one height one prime ideal, and every height one 
prime ideal is generated by a prime element.
Torsion free CGL extensions (for Cauchon--Goodearl--Letzter) 
are skew polynomial algebras  with a rational action of a torus, 
satisfying certain general conditions, see 
\cite[Definition 3.1]{LLR}. Launois, Lenagan, and 
Rigal \cite[Theorem 3.7]{LLR} proved that every 
torsion free CGL extension is a noetherian unique factorization domain. 
The algebras $\UU^w_\pm$ are all torsion free CGL extensions
see \cite{MC}; thus they are all 
noetherian unique factorization domains. 

For $y \in W$, $y \leq w$ define the ideals
\begin{align}
\label{II}
I^+_w(y) &= \Span \{(c^+_{w,\la})^{-1} \xi \mid 
\la \in P, \xi \in (V^+_w(\la))^*, \,
\xi \perp (V^+_w(\la) \cap \UU_- V(\la)_{y \la} )
\},
\\
\label{III}
I^-_w(y) &= \Span \{(c^-_{w,\la})^{-1} \xi \mid 
\la \in P, \xi \in (V^-_w(\la))^*, \,
\xi \perp (V^-_w(\la) \cap \UU_+ V(-w_0 \la)_{-y \la} )
\}
\end{align}
of $S^+_w$ and $S^-_w$, respectively, 
using the identifications \eqref{ident}.

In \cite{Y3}, using results of Gorelik \cite{G}, we proved that 
the algebras $\UU^w_-$ (and thus $S^+_w$) are catenary and 
that the $H$-invariant height one prime ideals of $\UU^w_-$
(with respect to the conjugation action of $H$)
are precisely the ideals $\varphi^+_w(I^+_w(s_i))$ 
for $i \in \SS(w)$. The analogous fact for $\UU^w_+$
is proved 
by interchanging the role of plus and minus generators $X^\pm_i$.
Since $\UU^w_\pm$ are noetherian unique factorization domains 
and a normal element of  $\UU^w_\pm$ is homogeneous, if and only 
if it generates an $H$-invariant ideal, we have:
\ble{numberp} The number of pairwise nonproportional 
homogeneous prime elements
of $\UU^w_\pm$ is equal to $|\SS(w)|$. 
\ele
\subsection{A lemma on diagonal automorphisms of $\UU^w_\pm$}
\label{5.3}
Let 
\[
w= s_{i_1} \ldots s_{i_l}
\]
be a reduced expression
of $w$, $l = l(w)$. Let $\be_1, \ldots, \be_l$ and
$X^\pm_{\be_1}, \ldots, X^\pm_{\be_l}$ be 
the roots and root vectors, given by 
\eqref{beta} and \eqref{rootv}, respectively.

\ble{aut} If $\psi \in \Aut(\UU^w_\pm)$ is such that 
\[
\psi(X^\pm_{\be_j})  = q_{i_j}^{k_j} X^\pm_{\be_j}
\]
for some $k_1, \ldots, k_{l(w)} \in \Zset$, then 
there exists $\de_\pm \in P$ such that 
$\lcor \de_\pm, \pm \be_j\spcheck \rcor = k_j$, 
(i.e. $\lcor \de_\pm, \pm \be_j \rcor = d_j k_j$),
for all $j =1, \ldots, l(w)$.
\ele

Recall from \S \ref{2.1} that $q_i = q^{d_i}$, 
where $(d_1, \ldots, d_r)$ is the vector 
of relatively prime positive integers 
symmetrizing the Cartan matrix of $\g$.

\leref{aut} (and the statement in \reref{aut2} below) 
are well known and easy to prove 
for various special cases, e.g. the algebras of 
quantum matrices or $w= w_0$. The emphasis here is 
on the validity of the statement
for all $\g$ and $w \in W$.  
\\ \hfill \\
{\em{Proof of \leref{aut}.}} 
We argue by induction on $l(w)$, the case $l(w) =0$ 
being trivial. Assume that the statement of the lemma 
is true for $w \in W$ of length $l$. Let 
$w' \in W$, $l(w') =l+1$, and 
$s_{i_1} \ldots s_{i_l} s_{i_{l+1}}$ be a reduced
expression of $w'$. Denote $w=s_{i_1} \ldots s_{i_l}$.
For this reduced expression of $w$, denote 
by $\be_j$ and $X^\pm_{\be_j}$, $j=1, \ldots, l$ 
the roots and root vectors of $\UU^w_\pm$ 
given by \eqref{beta} and \eqref{rootv}.
Denote 
\[
\beta_{l+1} = s_{i_1} \ldots s_{i_l} (\al_{i_{l+1}})
\; \; 
\mbox{and}
\; \; 
X^\pm_{\beta_{l+1}} = T_{i_1} \ldots T_{i_l} (X^\pm_{i_{l+1}})
\in \UU^{w'}_\pm.
\]
Clearly $\UU^w_\pm \subset \UU^{w'}_\pm$.

Let $\psi \in \Aut(\UU^{w'}_\pm)$ and
\begin{equation}
\label{l1}
\psi(X^\pm_{\be_j})  = q_{i_j}^{k_j} X^\pm_{\be_j}, 
\quad \forall j = 1, \ldots, l+1,
\end{equation}
for some $k_1, \ldots, k_{l+1} \in \Zset$. Then 
$\psi$ restricts 
to an automorphism of $\UU^w_\pm$ satisfying the 
assumptions of the lemma. Applying the inductive 
assumption, we obtain that there exists
$\de_\pm \in P$ such that 
\begin{equation}
\label{we}
\lcor \de_\pm, \pm \be_j\spcheck \rcor = k_j, \; \; 
\mbox{i.e.} \; \; 
\lcor \de_\pm, \pm \be_j \rcor = d_j k_j, \; \; 
\forall j =1, \ldots, l.
\end{equation}
By \reref{qnormal1} we can assume 
that $\de \in P_{\SS(w)}.$

First, consider the case when there exists 
$j \in \{ 1, \ldots, l \}$ such that 
\[
X^\pm_{\be_j} X^\pm_{\be_{l+1}} - q^{\lcor \be_j, \be_{l+1} \rcor }
X^\pm_{\be_{l+1}} X^\pm_{\be_j} \neq 0. 
\]
The Levendorskii--Soibelman straightening rule
\eqref{LS0}, the fact that $\UU^{w'}_\pm$ is $Q$-graded 
by $\deg X^\pm_{\be_1} = \pm \be_1, \ldots, 
\deg X^\pm_{\be_{l+1}} = \pm \be_{l+1}$, and \eqref{we} 
imply:
\begin{multline*}
\psi(X^\pm_{\be_j} X^\pm_{\be_{l+1}} - q^{\lcor \be_j, \be_{l+1} \rcor }
X^\pm_{\be_{l+1}} X^\pm_{\be_j}) 
\\
=q^{ \lcor \de_-, \pm (\be_j + \be_{l+1}) \rcor } 
(X^\pm_{\be_j} X^\pm_{\be_{l+1}} - q^{\lcor \be_j, \be_{l+1} \rcor }
X^\pm_{\be_{l+1}} X^\pm_{\be_j}).
\end{multline*}
Since $\UU^w_\pm$ is a domain, from the above, \eqref{l1} and \eqref{we}, 
we obtain
\[
\lcor \de_\pm, \pm \be_{l+1} \rcor = d_{l+1} k_{l+1}.
\]  
Thus the weight $\de_\pm$ for $w$ also works for $w'$, 
which proves the statement of the lemma.

Now consider the case when 
\begin{equation}
\label{co}
X^\pm_{\be_j} X^\pm_{\be_{l+1}} - q^{\lcor \be_j, \be_{l+1} \rcor }
X^\pm_{\be_{l+1}} X^\pm_{\be_j} = 0, \quad
\forall j =1, \ldots, l.
\end{equation}
\thref{DKP} implies that $X^\pm_{\be_{l+1}}$ is a 
homogeneous prime element of $\UU^{w'}_\pm$. 
\thref{DKP} and \eqref{co} also imply that each 
homogeneous prime element of $\UU^w_\pm$ is a 
prime element of $\UU^{w'}_\pm$. Therefore
the number of homogeneous prime elements 
of $\UU^{w'}_\pm$ is strictly greater than that 
of $\UU^w_\pm$. \leref{numberp} and 
$|\SS(w')| \leq |\SS(w)|+ 1$ 
imply $|\SS(w')| = |\SS(w)|+ 1$. 
By part (ii) of \leref{Iw} 
\begin{equation}
\label{ni}
i_{l+1} \notin \SS(w), \; \;
\mbox{i.e.} \; \;  
\al_{i_{l+1}}\spcheck \notin Q_{\SS(w)}\spcheck.
\end{equation}
Since 
$\beta_{l+1} = s_{i_1} \ldots s_{i_l} (\al_{i_{l+1}})$,
\begin{equation}
\label{bel1}
\be_{l+1}\spcheck = \al_{i_{l+1}}\spcheck 
+ \sum_{i \in \SS(w)} m_i \al_i\spcheck,  
\end{equation}
for some $\{m_i \in \Zset \mid i \in \SS(w) \}$.
Set
\[
\de'_\pm = \de_\pm \pm
\Big( k_{l+1} \mp \sum_{i \in \SS(w)} m_i 
\lcor \de_\pm, \al_i\spcheck \rcor \Big)
\om_{l+1}. 
\]
Because $\de_\pm \in P_{\SS(w)}$, \eqref{ni} implies that 
$\lcor \de_\pm, \al_{i_{l+1}}\spcheck \rcor =0$. Therefore
\begin{equation}
\label{e0}
\lcor \de'_\pm, \pm \al_{i_{l+1}}\spcheck \rcor 
= k_{l+1} - \sum_{i \in \SS(w)} m_i
\lcor \de_\pm, \pm \al_i\spcheck \rcor.
\end{equation}
From \eqref{ni} we also obtain that
\begin{equation}
\label{e1}
\lcor \de'_\pm, \al_i\spcheck \rcor = 
\lcor \de_\pm, \al_i\spcheck \rcor, \quad \forall i \in \SS(w).
\end{equation}
The induction hypothesis and \leref{Iw} (ii) imply:
\[
\lcor \de'_\pm, \pm \be_j\spcheck \rcor = 
\lcor \de_\pm, \pm \be_j\spcheck \rcor = k_j,
\quad \forall j=1, \ldots, l.
\]
Combining \eqref{bel1}, \eqref{e0} and \eqref{e1}, 
we obtain 
\begin{align*}
&\lcor \de'_\pm, \pm \be_{l+1}\spcheck \rcor = 
\lcor \de'_\pm, \pm \al_{i_{l+1}}\spcheck \rcor 
+ \sum_{i \in \SS(w)} m_i 
\lcor \de', \pm \al_i\spcheck \rcor
\\
= &k_{l+1} - \sum_{i \in \SS(w)} m_i
\lcor \de_\pm, \pm \al_i\spcheck \rcor 
+ 
\sum_{i \in \SS(w)} m_i
\lcor \de'_\pm, \pm \al_i\spcheck \rcor =
k_{l+1},
\end{align*}
which completes the proof of the lemma.
\qed
\bre{aut2} Define an action of the torus
$\Tset^{|\SS(w)|}= (\KK^*)^{\times |\SS(w)|}$ 
on $\UU^w_\pm$ by 
\begin{equation}
\label{HT}
t \cdot X^\pm_{\be_j} = 
\Big( \prod_{ i \in \SS(w)}
t_i^{ \lcor \om_i\spcheck, \be_j \rcor } 
\Big) 
X^\pm_{\be_j},
\end{equation}
for $t = (t_i)_{i \in \SS(w)} \in \Tset^{|\SS(w)|}$,
in terms of the generators \eqref{rootv} of $\UU^w_\pm$.
Here $\om_1\spcheck, \ldots, \om_r\spcheck$
denote the fundamental coweights of $\g$.
This is an action by algebra automorphisms 
since the algebras $\UU^w_\pm$ are $Q_{\SS(w)}$-graded
by \eqref{non0U}.
Analogously to the proof of \leref{aut} one shows:

{\em{
If $\psi \in \Aut(\UU^w_\pm)$ is such that
\[
\psi(X^\pm_{\be_j})  = p_j X^\pm_{\be_j}, \quad
\forall  
j=1, \ldots, l,
\]
for some $p_j \in \KK^*$, then
there exists $t \in \Tset^r$ 
such that $\psi(x) = t \cdot x$, 
$\forall x \in \UU^w_\pm$.
}}
\ere
\subsection{Proof of \prref{normalp}}
\label{5.4}
\noindent
Assume that $u \in \UU^w_\pm$ is a nonzero homogeneous normal element. 
We will prove that there exists $\de_\pm \in P_{\SS(w)}$ such that 
\[
u X^\pm_{\be_j} = 
q^{ \lcor \de_\pm, \pm \be_j \rcor }
X^\pm_{\be_j} u, \quad \forall j=1, \ldots, l.
\]
Then applying the graded (anti)isomorphism from \thref{isom}
implies the statement of the proposition.

Fix $j \in \{ 1, \ldots, l \}$. Then 
\begin{equation}
\label{co1}
u X^\pm_{\beta_j} = Y_j u \; \; 
\mbox{for some} \; \; Y \in (\UU^w_\pm)_{\pm \beta_j}.
\end{equation}
Recall the notation \eqref{monomial}, and the
notions of highest term of a nonzero element of
$\UU^w_\pm$ and degree of a monomial from \S 2.7. 
Assume that the highest term of $u$ 
is has degree ${\bf{n}}$ for some
${\bf{n}} \in \Nset^{\times l}$.
Denote ${\bf{e}}_j = (0, \ldots, 0, 1, 0, \ldots, 0)$, where 
1 is in position $j$. Then, by \leref{LSmult} 
the highest term of the left hand side of \eqref{co1} 
has degree ${\bf{n}} + {\bf{e}}_j$.
Again applying \leref{LSmult}, we obtain that 
the highest term of $Y_j$ has degree ${\bf{e}}_j$, i.e. 
the highest term is a nonzero scalar multiple of
$X^\pm_{\be_j}$. 
At the same time $Y \in (\UU^w_\pm)_{\pm \beta_j}$; 
that is 
\begin{equation}
\label{Yj}
Y_j \in \Span \{(X^\pm)^{ {\bf{n}}' } \mid 
{\bf{n}}'= (n'_1, \ldots, n'_l) \in \Nset^{\times l}, 
n'_1 \be_1 + \ldots + n'_l \be_l = \be_j 
\}.
\end{equation}
It is well known that the ordering of the roots
\begin{equation}
\label{list}
\be_1, \ldots, \be_l
\end{equation}
of $\De^+ \cap w(\De^-)$ is convex, i.e. if a root in \eqref{list} 
is equal to the sum of two other roots in \eqref{list}, 
then it is listed in between.
Moreover if a root of $\g$ is the sum of two roots of 
$\De^+ \cap w(\De^-)$, then it belongs to 
$\De^+ \cap w(\De^-)$. This implies that if a root 
in the list \eqref{list} is a positive integral combination 
of several roots in \eqref{list}, then it is listed between
the leftmost and rightmost ones. This property 
and \eqref{Yj} imply that 
the highest term of $Y_j$ will not be 
a nonzero scalar multiple of $X^\pm_{\be_j}$ unless 
$Y_j$ is itself a scalar multiple of $X^\pm_{\be_j}$. Therefore
\[
u X^\pm_{\be_j} = p_j X^\pm_{\be_j} u
\]
for some $p_j \in \KK^*$. Comparing the 
highest terms of both sides and using \leref{LSmult}, 
we obtain
\[
u X^\pm_{\be_j} = q_{i_j}^{k_j} X^\pm_{\be_j} u, \quad
j=1, \ldots, l,
\]
for some $k_j \in \Zset$. Repeated applications 
of \eqref{LS0} give
\begin{align*}
k_j & = \sum_{k=1}^{j-1} \frac{n_k \lcor \be_j, \be_k \rcor}{d_{i_j}}
- \sum_{k=j+1}^{l} \frac{n_k \lcor \be_j, \be_k \rcor}{d_{i_j}}
\\
& = \sum_{k=1}^{j-1} n_k \lcor \be_j\spcheck, \be_k \rcor
- \sum_{k=j+1}^{l} n_k \lcor \be_j\spcheck, \be_k \rcor 
\in \Zset.
\end{align*}
Here we used the fact that 
$\lcor \be_j, \be_j \rcor = 
\lcor \al_{i_j}, \al_{i_j} \rcor = 
d_{i_j}$.

Applying \leref{aut}, we obtain that there exists 
$\de_\pm \in P_{\SS(w)}$ such that 
$\lcor \de_\pm, \pm \be_j \rcor$ $= d_{i_j} k_j$
for all $j=1, \ldots, l$; that is
\[
u X^\pm_{\be_j} = 
q^{ \lcor \de_\pm, \pm \be_j \rcor }
X^\pm_{\be_j} u, \quad \forall j=1, \ldots, l.
\]
This completes the proof of \prref{normalp}.
\qed
\subsection{Proof of \thref{normal1} }
\label{5.5}
Denote by 
\begin{equation}
\label{MM}
M^\pm_w = \{ p d^\pm_{w, \la} \mid
p \in \KK^*, \la \in P^+_{\SS(w)} \}
\end{equation}
the multiplicative subset of all
nonzero homogeneous normal elements of 
$S^\pm_w$, cf. \thref{normal1} (i).
We start with a lemma which narrows down the set of 
homogeneous normal elements of $S^+_w$. 
\ble{more-n} The set of homogeneous normal elements of 
$S^\pm_w$ consists of those elements of 
$S^\pm_w[(M^\pm_w)^{-1}]$ which have the form
\[
p \prod_{i \in \SS(w)} (d^\pm_{w, \om_i})^{n_i}
\]
for some $p \in \KK$, $n_i \in \Zset$ and belong to $S^\pm_w$.
The product over $i$ is taken in any fixed order 
as in \eqref{dprod}.
\ele

Each reduced expression $w =s_{i_1} \ldots s_{i_l}$ 
gives rise to a presentation of the algebra $\UU^w_\pm$ as 
an iterated skew polynomial algebra
\begin{equation}
\label{skew}
\KK[X^\pm_{\be_1}][X^\pm_{\be_2}; \tau_2, \theta_2] 
\ldots [X^\pm_{\be_l}; \tau_l, \theta_l] 
\end{equation}
where for $j=1, \ldots, l=l(w)$, 
$\tau_j$ is an automorphism of $(j-1)$-st
algebra in the extension and $\theta_j$ is a
$\tau_j$-derivation of the same algebra. 
(One constructs $\tau_j$ and $\theta_j$ from the 
Levendorskii--Soibelman straightening 
rule \eqref{LS0}, see \cite{MC}.)
Moreover the following conditions are trivially 
satisfied (and also follow from the property 
that $\UU^w_\pm$ are CGL extensions):

(i) All $X^\pm_{\be_1}, \ldots, X^\pm_{\be_l}$ 
are eigenvalues of $H$ under the conjugation action.

(ii) For $j = 1, \ldots, l$ there exist elements of 
$H_j \in H$ such that $\tau_j (X^\pm_{\be_k}) = H_j^{\pm 1} X^\pm_{\be_k} H_j^{\mp 1}$
for $j>k$ and the $H_j$ eigenvalue of 
$X^\pm_{\be_j}$ is not a root of unity for all $j$. 

Goodearl proved \cite[Theorem 6.4.II]{BG} 
that, if $A$ is an iterated skew polynomial algebra as 
in \eqref{skew} which satisfies the properties 
(i)--(ii) above, then every $H$-prime $I$ of $A$ 
is strongly rational, i.e. $Z(\Fract A/I)^H = \KK$.
Strictly speaking we need to use the 
extension of the conjugation action of $H$ 
on $\UU^w_\pm$ to the torus action \eqref{HT}
(the $H$-invariant ideals being the same as 
the $\Tset^{|\SS(w)|}$-invariant ideals).
Using the (anti)isomorphisms
$\varphi_w^\pm \colon S^\pm_w \to \UU^w_\mp$ 
(see \thref{isom}) and applying this result 
to the $\{0\}$ ideals of the algebras $\UU^w_\mp$, we obtain: 
\begin{equation}
\label{zero-cent}
Z(S^\pm_w[(M^\pm_w)^{-1}])_{0,0} = \KK.
\end{equation}
\\
\noindent
{\em{Proof of \leref{more-n}.}} Assume that $u$ is a 
nonzero homogeneous normal element of $S^\pm_w$. 
By \prref{normalp} it is $P$-normal. We then apply 
\thref{nS1} to obtain that there exists $\eta \in P_{\SS(w)}$ 
such that  $u \in (S^\pm_w)_{\pm(w - 1)\eta,0}$ 
and
\begin{equation}
\label{dede3}
u s = q^{\lcor - (w + 1 )\eta, \ga \rcor } s u, \quad
\forall s \in (S^\pm_w)_{-\ga,0}, \ga \in Q_{\SS(w)}.
\end{equation}
Let $\eta= \sum_{i \in \SS(w)} n_i \om_i$ 
for some $n_i \in \Zset$. Denote 
\[
u' = \prod_{i \in \SS(w)} (d^\pm_{w, \om_i})^{n_i},
\]
where the product over $i$ is taken in any order.
Then $u' \in (S^\pm_w)_{\pm (w - 1)\eta,0}$
and
\begin{equation}
\label{dede4}
u' s = q^{\lcor - (w + 1 )\eta, \ga \rcor } s u', \quad
\forall s \in (S^\pm_w)_{-\ga,0}, \ga \in Q_{\SS(w)},
\end{equation}
recall \eqref{ddcomm}.
Eq. \eqref{dede3} and \eqref{dede4} imply 
\[
u (u')^{-1} \in Z(S^\pm_w[(M^\pm_w)^{-1}])_{0,0}.
\]
From \eqref{zero-cent} we obtain that 
$u (u')^{-1} \in \KK^*$, i.e. 
\[
u = p u' = p \prod_{i \in \SS(w)} (d^\pm_{w, \om_i})^{n_i},
\]
for some $p \in \KK^*$. 
\qed
\\ \hfill \\
\noindent
{\em{Proof \thref{normal1}. Part (i)}}: Assume that 
$u \in S^\pm_w$ is a nonzero homogeneous normal element. 
\leref{more-n} implies that it is given by \eqref{dprod} for some 
$p \in \KK^*$, $n_i \in \Zset$. We claim that 
$n_i \in \Nset$ for all $i \in \SS(w)$. 
Assume that this is not the case. 
Then the element $u$ would be linearly independent from the set 
$d^\pm_{w, \la}$, $\la \in P^+_{\SS(w)}$. Indeed, if this is not 
the case, then after multiplying it with $d^\pm_{w, \mu}$ 
for some $\mu \in P^+_{\SS(w)}$ we will get a linear dependence
in the set $\{ d^\pm_{w, \la} \}_{\la \in P^+_{\SS(w)}}$, which 
contradicts with the first part of 
\thref{freeS}.
Therefore $u \notin N^\pm_w$ and for 
some $\mu \in P^+_{\SS(w)}$, $d^\pm_{w, \mu} u \in N^\pm_w$. 
This contradicts with the fact that $S^\pm_w$ is a free 
(left and right) $N^\pm_w$-module (by \thref{freeS} (ii)) and completes 
the proof of part (i) of \thref{normal1}.

{\em{Part (ii):}} By the first part of the theorem
each homogeneous normal element of $S^\pm_w$ has 
the form \eqref{dprod},  for some $p \in \KK^*$, 
$n_i \in \Nset$. Therefore the set of
homogeneous prime elements of $S^\pm_w$ 
is a subset of 
$\{ p d^\pm_{w, \om_i} \mid 
p \in \KK^*, i \in \SS(w) \}$.  
By \leref{numberp} $S^\pm_w$ has at least 
$|\II(w)|$ pairwise nonproportional homogeneous
prime elements. This is only possible if
$d^\pm_{w, \om_i}$ are 
prime elements of $S^\pm_w$, for all 
$i \in \SS(w)$. They are linearly 
independent because of \thref{freeS} (i).
\qed
\\

As a corollary of \thref{normal1} (ii), we obtain 
explicit formulas and generators for the height one $H$-primes 
$I^\pm_w(s_i)$, $i \in \SS(w)$ of $S^\pm_w$, 
(recall \eqref{II}--\eqref{III}) under the general conditions
on $\KK$ and $q$.

\bpr{generators} For any base field $\KK$, 
$q \in \KK^*$ not a root of unity, 
$w \in W$, and $i \in \SS(w)$ we have
\begin{equation}
\label{h1pr}
I_w^\pm(s_i) = S^\pm_w d^\pm_{w, \om_i}.
\end{equation}
\epr

\begin{proof}
Combining the Launois--Lenagan--Rigal 
result \cite[Theorem 3.7]{LLR} implying 
that $\UU^w_\mp$ is a unique factorization 
domain, the fact that the ideals $I^\pm_w(s_i)$
are height one prime ideals of $S^\pm_w$ and 
part (ii) of \thref{normal1}, we obtain that for 
each $i \in \SS(w)$ there exists $k \in \SS(w)$ 
such that 
\[
I^\pm_w(s_i) = S^\pm_w d^\pm_{w, \om_k}.
\]
Since $d^\pm_{w, \om_i} \in I^\pm_w(s_i)$, 
this is only possible if $k =i$, which 
establishes \eqref{h1pr}.
\end{proof}
\subsection{Prime and primitive ideals in the 
$\{0\}$-stratum 
of $\Spec S^\pm_w$.}
\label{5.6}
As an application of Theorems \ref{tfreeS} and \ref{tnormal1},
we obtain a formula for the prime (and a more explicit formula 
for the primitive) ideals of the algebras 
$S^\pm_w$ lying in the Goodearl--Letzter 
stratum over the $\{0\}$ ideal. As usual, some of the
results for primitive ideals require the base field $\KK$ to be 
algebraically closed. There will be no such
restriction for the results on prime ideals, 
which are valid for arbitrary base fields $\KK$.
Similar arguments are applied in the next subsection 
to obtain a classification 
of all prime elements of the algebras $S^\pm_w$ 
(and $\UU^w_\pm$). In particular, this gives 
explicit formulas for all height one prime 
ideals of $S^\pm_w$ .

Via the (anti)isomorphism of \thref{isom} these results give 
similar explicit formulas for the prime/primitive ideals 
in the 
$\{0\}$-stratum of $\Spec \UU^w_\mp$. The 
restatement is straightforward and will not be 
formulated separately.

Recall \cite{GL} that the 
Goodearl--Letzter
$\{0\}$-stratum of 
$\Spec S^\pm_w$ is defined by
\[
\Spec_{ \{ 0 \} } S^\pm_w = 
\{ I \in \Spec S^\pm_w \mid \cap_{t \in \Tset^r} 
t \cdot I = \{ 0 \},    
\]
where we use the rational action \eqref{Tract0} 
of $\Tset^r$ on $S^\pm_w$. Set
$\Prim_{ \{ 0 \} } S^\pm_w =  \Prim S^\pm_w \cap
\Spec_{ \{ 0 \} } S^\pm_w$. 

Recall also that $M^\pm_w$ denotes the 
multiplicative subset \eqref{MM}
of all nonzero homogeneous normal 
elements of $S^\pm_w$, see
\thref{normal1} (i).
First, we obtain a description of the center
$Z(S^\pm_w[ (M^\pm_w)^{-1} ])$.
Each $\mu \in P_{\SS(w)}$ can be represented in a unique way as $\mu = \mu_+ - \mu_-$
for some $\mu_+, \mu_- \in P^+_{\SS(w)}$ with disjoint support, see \eqref{supp}.
For $\mu \in P_{\SS(w)}$ define
\begin{equation}
\label{dla}
d^\pm_{w, \mu} = (d^\pm_{w, \mu_-})^{-1} d^\pm_{w, \mu_+} \in
N^\pm_w[ (M^\pm_w)^{-1}].
\end{equation}
It follows from \eqref{ddcomm} and \eqref{multd} 
that for all $\mu_1, \mu_2 \in P_{\SS(w)}$,
\begin{equation}
\label{dprod2}
d^\pm_{w,\mu_1} d^\pm_{w,\mu_2} = q^{j(\mu_1, \mu_2)}
d^\pm_{w,\mu_1 + \mu_2}, \; \;
\mbox{for some} \; \; 
j(\mu_1, \mu_2) \in \Zset.
\end{equation}
(Recall that $d^\pm_{w, \mu} \in \KK^*$ for all 
$\mu \in P_{\II(w)}^+$. Because of this and \eqref{multd}, 
one does not need to extend the definition \eqref{dla} 
to $\mu \in P$.)
\thref{freeS} (i) implies that
the localization $N^\pm_w[ (M^\pm_w)^{-1}]$
is isomorphic a quantum torus over $\KK$ 
of dimension $|\SS(w)|$. In particular,
\begin{equation}
\label{basisNM}
\{ d^\pm_{w, \mu} \mid \mu \in P_{\SS(w)} \}
\; \; 
\mbox{is a $\KK$-basis of} \; \;  
N^\pm_w[ (M^\pm_w)^{-1}].
\end{equation}
Applying \eqref{ddcomm} we also obtain 
\begin{equation}
\label{ddcomm2}
d_{w, \mu}^\pm s = q^{- \lcor (w + 1 ) \mu , \ga \rcor } 
s d_{w, \mu}^\pm, \quad \forall \mu \in P_{\SS(w)},
s \in (S_w^\pm[ (M^\pm_w)^{-1}])_{-\ga, 0}, \ga \in Q_{\SS(w)}.
\end{equation}

Define the lattice
\begin{align}
\label{kl}
\Kl (w) &=  \{ \mu \in P_{\SS(w)} \mid (w + 1) \mu  \in (Q_{\SS(w)})^\perp \}
\\
\nn
&=  \{ \mu \in P_{\SS(w)} \mid (w + 1) \mu  \in P_{\II(w)} \}
\\
\nn
&=  \{ \mu \in P_{\SS(w)} \mid \exists \nu \in P_{\II(w)} \; \; 
\mbox{such that} \; \; \mu  + \nu/2 \in \ker (w +1) \}.
\end{align}
The first equality follows from 
$(Q_{\SS(w)})^\perp \cap P = P_{\II(w)}$. The second equality 
follows from the fact that $w(\nu) = \nu$ for all 
$\nu \in P_{\II(w)}$, thus for any $\mu \in P$ and $\nu \in P_{\II(w)}$:
\[
(w+1) \mu = \nu, \; \; \mbox{if and only if} \; \; (w+1)(\mu - \nu/2) =0.
\]

The lattice $\Kl(w)$ has rank 
\begin{equation}
\label{mw}
m(w) := \dim \ker (w + 1).
\end{equation}
To see 
this, denote the projection 
$\sig \colon P_{\SS(w)} \oplus P_{\II(w)}/2 \to P_{\SS(w)}$
along $P_{\II(w)}/2$. 
The third equality in \eqref{kl} implies that 
\[
\Kl(w) = \sig \big( \ker (w+1) \cap (P_{\SS(w)} \oplus P_{\II(w)}/2) \big).
\]
The statement follows from the facts that 
$\ker (w+1) \cap (P_{\SS(w)} \oplus P_{\II(w)}/2)$ is a lattice of rank
$\dim \ker(w+1)$ and the restriction
\[
\sig \colon \ker (w+1) \cap (P_{\SS(w)} \oplus P_{\II(w)}/2) \to \Kl(w)
\]
is bijective.

Fix a basis $\mu^{(1)}, \ldots, \mu^{(m(w))}$ of $\Kl(w)$.
For $ \mu = k_1 \mu^{(1)} + \ldots + k_{m(w)} \mu^{(m(w))} \in \Kl(w)$ 
define
\[
e^\pm_{w, \mu} =
(d^\pm_{w,\mu^{(1)}})^{k_1} \ldots (d^\pm_{w, \mu^{(m(w))} })^{ k_{m(w)} }.
\]
We have $e^\pm_{w, \mu} e^\pm_{w, \mu'} = e^\pm_{w, \mu + \mu'}$, 
for all $\mu, \mu' \in \Kl(w)$.
By \eqref{dprod2} for all $\mu \in \Kl(w)$,
\begin{equation}
\label{de}
e^\pm_{w, \mu} = q^{j_\mu} d^\pm_{w, \mu},
\quad  \mbox{for some} \; \; i_\mu \in \Zset.
\end{equation}

Denote by $A_w^\pm$ the subalgebra of the quantum torus 
$N^\pm_w[ (M^\pm_w)^{-1} ]$
generated by 
\begin{equation}
\label{Zc}
d^\pm_{w, \mu^{(i)}}, (d^\pm_{w, \mu^{(i)}})^{-1}, \quad i=1, \ldots, m(w).
\end{equation}
\thref{freeS} (i) and \eqref{ddcomm2} imply that $A_w^\pm$ is a Laurent polynomial 
algebra over $\KK$ of dimension $m(w)$ with generators \eqref{Zc}.
The set
\begin{equation}
\label{MbM}
\{ e^\pm_{w, \mu } \mid \mu \in \Kl(w) \}
\end{equation}
is a $\KK$-basis of $A_w^\pm$. Of course, the same is true
for the set $\{ d^\pm_{w, \mu } \mid \mu \in \Kl(w) \}$.

\bre{wex} The one half in \eqref{kl} is needed as is shown by the 
following simple example. Let $\g = {\mathfrak{sl}}_3$, $w = s_1$. Then 
$\SS(s_1) = \{1\}$ and $\II(s_1) = \{ 2 \}$. Moreover,
$\Kl(s_1) = \Zset \om_1$ and 
\[
\om_1 - \om_2/2 \in \ker (s_1 +1),
\]
cf. the third equality in \eqref{kl}. In this 
case $S_{s_1}^\pm = \KK [d^\pm_{s_1, \om_1}]$ is a 
polynomial ring and 
$A_w^\pm = \KK [d^\pm_{s_1, \om_1}, (d^\pm_{s_1, \om_1})^{-1} ]$.
\ere

\ble{Z} For an arbitrary base field $\KK$, $q \in \KK^*$ not a root
of unity, and $w \in W$, the center 
$Z( S^\pm_w[ (M^\pm_w)^{-1} ])$ coincides with the algebra $A_w^\pm$.
\ele

The special case of \leref{Z} when $w$ is equal to the longest element
$w_0$ of $W$ is due to Caldero \cite{Ca0}. The special case of \leref{Z} 
for the algebras of quantum matrices is due to Launois and Lenagan 
\cite{LL}. The fact that 
$\Spec_{ \{ 0 \} } \UU^w_\pm$ has dimension equal to 
$\dim \ker (w+1)$ was obtained by Bell and Launois \cite{BL0}.

\bre{cent} \thref{freeS} (i) and \leref{Z} imply that the center of the 
algebra $S^\pm_w$ is the algebra
\[
Z(S^\pm_w) = \{ d^\pm_{w, \mu} \mid \mu \in \Kl(w) \cap P^+_{\SS(w)} \}.
\]
Often this algebra is much smaller than $A^\pm_w$. For instance 
in many cases it is trivial while $A^\pm_w$ is not. This is why 
it is more important to study the structure 
of $S^\pm_w$ as an $N^\pm_w$-module rather than $Z(S^\pm_w)$-module.
\ere
\noindent
{\em{Proof of \leref{Z}.}}
It follows from of \eqref{ddcomm2} and the first equality in \eqref{kl}
that $A_w^\pm \subseteq Z( S^\pm_w[ (M^\pm_w)^{-1} ])$.
To show the opposite inclusion, let 
$(d^\pm_{w, \la})^{-1} u \in Z( S^\pm_w[ (M^\pm_w)^{-1} ])$
for some homogeneous element 
$u \in S^\pm_w$ and $\la \in P^+_{\SS(w)}$. Then $u \in S^\pm_w$ should 
be a homogeneous normal element. Applying \thref{normal1} (i) 
and \eqref{dprod2} we obtain that 
\begin{equation}
\label{cen}
(d^\pm_{w, \la})^{-1} u = p d^\pm_{w, \mu}
\end{equation}
for some $p \in \KK^*$, $\mu \in P_{\SS(w)}$. Using 
\eqref{ddcomm2}, we obtain that 
$d^\pm_{w, \mu} \in Z( S^\pm_w[ (M^\pm_w)^{-1} ])$
if and only if $\mu \in \Kl(w)$, cf. the first 
equality in \eqref{kl}. Therefore 
$(d^\pm_{w, \la})^{-1} u$ is a scalar multiple
of one of the elements in \eqref{MbM}, and thus belongs 
to $A_w^\pm$. 
\qed
\\

In \cite{Y3} we prove that $\Spec \UU^w_\pm$ 
is normally separated, under the same 
general assumption on $\KK$ and $q$ as the ones in
this paper. Using the (anti)isomorphism from 
\thref{isom}, we obtain that the same 
is true for the algebras $S^\pm_w$.
By \cite[Theorems 5.3 and 5.5]{GK} every prime ideal in 
$\Spec_{ \{ 0 \} } S^\pm_w$ is of the form
\begin{equation}
\label{primUw}
(S^\pm_w[ (M^\pm_w)^{-1} ] . J^0) \cap S^\pm_w, 
\end{equation}
for some prime ideal $J^0$ of 
$Z( S^\pm_w[ (M^\pm_w)^{-1} ])$. Moreover 
each primitive ideal in $\Prim_{ \{ 0 \} } S^\pm_w$ 
has the form \eqref{primUw} for a maximal ideal $J^0$ of 
$Z( S^\pm_w[ (M^\pm_w)^{-1} ])$. Applying
the freeness result \thref{freeS} (ii) and \leref{Z} 
leads to the following result. 
Its proof is straightforward 
and will be omitted.

\bpr{idealsU} Assume that $\KK$ is an arbitrary 
base field, $q \in \KK^*$ is not a root of unity, 
and $w \in W$. Then:

(i) All prime ideals in $\Spec_{ \{ 0 \} } S^\pm_w$
have the form 
\begin{equation}
\label{ideal}
\big( (J^0. N^\pm_w[(M^\pm_w)^{-1}]) \cap N^\pm_w \big) . S^\pm_w 
\end{equation}
for some prime ideal $J^0$ of the Laurent polynomial ring $A_w^\pm$, 
see \leref{Z}.

(ii) The primitive ideals in $\Prim_{ \{ 0 \} } S^\pm_w$
are the ideals given by \eqref{ideal} for maximal
ideals $J^0$ of $A_w^\pm$.
\epr

The point of \prref{idealsU} is that it reduces the 
possibly complicated intersections from \eqref{primUw}
in the algebras $S^\pm_w$ to the intersections \eqref{ideal} 
inside the quantum affine space algebras $N^\pm_w$. The latter
intersections are obviously 
much simpler. Moreover the centers $Z( S^\pm_w[ (M^\pm_w)^{-1} ])$
are substituted by the explicit Laurent polynomial 
algebras $A_w^\pm$.

Next, we proceed with making the description from part (ii) 
of \prref{idealsU} even more explicit.
For $\la, \la' \in P^+_{\SS(w)}$, such that 
\[
\la' - \la = k_1 \mu^{(1)} + \ldots + k_{m(w)} \mu^{(m(w))} \in \Kl(w)
\]
set 
\begin{multline}
\label{lo}
n^\pm_{\la, \la'} = \mp 2 \sum_i k_i \lcor \mu^{(i)} , \la \rcor
\mp 2 \sum_{j<i} k_j \lcor \mu^{(j)} , \mu^{(i)} \rcor
\\ 
\mp \sum_i |k_i|(|k_i|-1) \lcor \mu^{(i)} , \mu^{(i)} \rcor
\pm 2 \sum_i k_i \lcor \mu^{(i)} , \mu^{(i)}_{-\sign(k_i)} \rcor.
\end{multline}
Applying repeatedly \eqref{ddcomm} and \eqref{multd}, and using the fact 
that $(w-1) \mu^{(i)}= - 2 \mu^{(i)}$, because $\mu^{(i)} \in \ker (w+1)$,
gives 
\begin{equation}
\label{dprod3}
d^\pm_{w, \la} e^\pm_{w, \la' - \la } = 
q^{ n^\pm_{\la, \la'} } d^\pm_{w, \la' }, \quad
\forall \la, \la' \in P^+_{\SS(w)} \; \; 
\mbox{such that} \; \; \la' - \la \in \Kl(w).
\end{equation}
We leave the details of this long but straightforward computation to the reader.
 
Denote by $J_{w, {\bf{1}} }^0$ the maximal ideal
of $A_w^\pm$ generated by
\[
d^\pm_{w,\mu^{(i)}} - 1, i=1, \ldots, m(w).
\]
Fix a set $\Lambda_{w} \subset P_{\SS(w)}$ 
of representatives of $P_{\SS(w)}/\Kl(w)$,
recall \eqref{kl}.
Let $\Lambda_{w}^+ \subset P^+_{\SS(w)}$ be a set 
of representatives of those cosets in 
$P_{\SS(w)}/\Kl(w)$ that intersect $P^+_{\SS(w)}$
nontrivially.

Denote by $J_{w, {\bf{1}} }$ the subspace of $N^\pm_w$, which is the span
over $\la \in \Lambda_w^+$ of all elements of the form
\[
\sum_{ \mu \in \Kl(w) \cap (- \la + P^+_{\SS(w)})}
p_{\la, \la+ \mu } d^\pm_{w, \la + \mu},
\]
for $p_{\la, \la + \mu} \in \KK$ such that 
\[
\sum_{ \mu \in \Kl(w) \cap (- \la + P^+_{\SS(w)}) } q^{- n^\pm_{\la, \la + \mu } }
p_{\la, \la + \mu } = 0.
\]
The next lemma proves that $J_{w, {\bf{1}} }$ is an ideal of 
$N^\pm_w$ and relates it to the setting of \prref{idealsU}.

\ble{Ind} Let $w \in W$. 

(i) One has
\begin{multline*}
J_{w, {\bf{1}} }^0 . N^\pm_w[(M^\pm_w)^{-1}] = \big\{ 
\sum_{\la \in \Lambda_w} \sum_{ \mu \in \Kl(w) } 
p_{\la, \la + \mu } d^\pm_{w, \la} e^\pm_{w, \mu}  \, \big|
p_{\la, \la + \mu } \in \KK,
\\
\sum_{ \mu \in \Kl(w) }
p_{\la, \la + \mu } = 0, 
\forall \la \in \Lambda_w \big\}.
\end{multline*}

(ii) The ideal
$\big( J_{w, {\bf{1}} }^0 . N^\pm_w[(M^\pm_w)^{-1}] \big) \cap N^\pm_w$
of $N^\pm_w$ equals $J_{w, {\bf{1}} }$.
\ele
\begin{proof} (i) Since $\{ d^\pm_{w, \mu} \mid \mu \in P_{\SS(w)} \}$
is a $\KK$-basis of $N^\pm_w[(M^\pm_w)^{-1}]$, from \eqref{dprod2},
\eqref{de}   
and the definition of $\Lambda_w$, we obtain that 
\begin{equation}
\label{llk}
\{ d^\pm_{w, \la} e^\pm_{w, \mu } \mid
\la \in \Lambda_w, \mu \in \Kl(w)  \} 
\; \; \mbox{is a $\KK$-basis of} \; \; N^\pm_w[(M^\pm_w)^{-1}].
\end{equation}
The statement of part (i) now follows from the definition 
of the ideal $J_{w, {\bf{1}} }^0$, and the facts that 
$e^\pm_{w, \mu} e^\pm_{w, \mu'} = e^\pm_{w, \mu + \mu'}$, 
$\forall \mu, \mu' \in \Kl(w)$ and 
$e^\pm_{w, \mu^{(i)}} = d^\pm_{w, \mu^{(i)}}$, $\forall i=1, \ldots m(w)$.

(ii) Part (i), \eqref{de}, \eqref{llk} and \thref{freeS} (i) imply that 
$\big( J_{w, {\bf{1}} }^0 . N^\pm_w[(M^\pm_w)^{-1}] \big) \cap N^\pm_w$
is equal to the ideal of $N^\pm_w$, which is the 
span over $\la \in \La_w^+$ of all elements of the form
\begin{equation}
\label{sum2}
\sum_{\mu \in \Kl(w) \cap (- \la + P^+_{\SS(w)})}
p_{\la, \la + \mu } d^\pm_{w, \la} e^\pm_{w, \mu},
\end{equation}
for $p_{\la, \la + \mu } \in \KK$ such that
\[
\sum_{ \mu \in \Kl(w) \cap (- \la + P^+_{\SS(w)} )}
p_{\la, \la+ \mu } = 0.
\]
It follows from \eqref{dprod3} that this is exactly the ideal 
$J_{w, {\bf{1}} }$.
\end{proof}

\bth{primitive} Assume that $\KK$ is an arbitrary base field, 
$q \in \KK^*$ is not a root of unity, and $w \in W$. Then 
$J_{w, {\bf{1}} } S^\pm_w$ is a primitive ideal in 
$\Prim_{ \{ 0 \} } S^\pm_w$. If the field $\KK$ is 
algebraically closed, then the primitive ideals of $S^\pm_w$ 
in the $\{0\}$-stratum of $\Prim S^\pm_w$ are the ideals
\[
t \cdot (J_{w, {\bf{1}} } S^\pm_w)
\]
for $t \in \Tset^r$ with respect to the action \eqref{Tract0}.
\eth

In the special case when $\g= {\mathfrak{sl}}_{r+1}$ and 
$w$ is equal to the longest element $w_0$ of $W$, 
\thref{primitive} and \coref{gen} below are due to Lopes, \cite{Lo2}.

We note that the freeness result of \thref{freeS2} provides an explicit 
formula for the primitive ideal $J_{w, {\bf{1}} } S^\pm_w$
of $S^\pm_w$. Indeed, we have that for each reduced expression 
$\vec{w}$ of $w$:
\begin{equation}
\label{idealJ}
J_{w, {\bf{1}} } S^\pm_w =
\bigoplus_{ {\bf{n}} \in \De(\vec{w}) } 
J_{w, {\bf{1}} }
\cdot (\varphi^\pm_w)^{-1} \big( (X^\mp)^{\bf{n}} \big), 
\end{equation}
cf. \eqref{De}.
\\ \hfill \\
\noindent
{\em{Proof of \thref{primitive}}}. The theorem follows from 
\prref{idealsU} (ii), \leref{Ind} (i), and the 
fact that $\Prim_{ \{ 0 \} } S^\pm_w$ is a single 
$\Tset^r$-orbit under \eqref{Tract0}, which
which is a consequence of \cite[Theorem 5.5]{GK}.
\qed
\\

The definition of the ideals $J_{w, {\bf{1}} }$ gives
immediately  efficient generating sets 
for the ideals $J_{w, {\bf{1}} } S^\pm_w$. Represent 
each $\mu \in \Kl(w)$ as $\mu = \mu_+ - \mu_-$ for 
$\mu_+, \mu_- \in P^+_{\SS(w)}$ with disjoint support, 
cf. \eqref{supp}. Then \eqref{dprod3} implies 
\begin{equation}
\label{Gel}
d^\pm_{w_\pm, \mu_-} ( 1 - e^\pm_{w_\pm, \mu} )
= d^\pm_{w_\pm, \mu_-} - q^{n^\pm_{\mu_-, \mu_+}} 
d^\pm_{w_\pm, \mu_-} \in J_{w, {\bf{1}} }.
\end{equation}
For all $\mu \in \Kl(w)$, the above are normal
elements of the algebras $N^\pm_w$ and $S^\pm_w$ 
since $e^\pm_{w_\pm, \mu} \in Z( S^\pm_w[(M^\pm_w)^{-1}])$.
From this and the definition of $J_{w, {\bf{1}} }$, 
we obtain
\begin{equation}
\label{Jid0}
J_{w, {\bf{1}} } = \sum_{\mu \in \Kl(w)} 
(d^\pm_{w_\pm, \mu_-} - q^{n^\pm_{\mu_-, \mu_+}} 
d^\pm_{w_\pm, \mu_-}) N^\pm_w  
\end{equation}
and
\begin{equation}
\label{Jid1}
J_{w, {\bf{1}} } S^\pm_w = \sum_{\mu \in \Kl(w)} 
(d^\pm_{w_\pm, \mu_-} - q^{n^\pm_{\mu_-, \mu_+}} 
d^\pm_{w_\pm, \mu_-}) S^\pm_w.  
\end{equation}
In each particular case one easily isolates a finite 
generating subset in \eqref{Jid0} (consisting 
of elements of the form \eqref{Gel} for $\mu$ 
in a finite subset of $\Kl(w)$).Then the same set 
(of normal elements of $S^\pm_w$)
generates the ideal $J_{w, {\bf{1}} } S^\pm_w$. 
Here is a simple general example of this.

\bco{gen} Assume that $w \in W$ is such that the lattice 
$\Kl(w)$ has a basis $\mu^{(1)}, \ldots \mu^{(m(w))}$, 
consisting of elements of $P^+_{\SS(w)}$ with pairwise
disjoint support. Then 
\begin{equation}
\label{genn}
J_{w, {\bf{1}} } S^\pm_w = \sum_{i=1}^{m(w)}
(1 - d^\pm_{w, \mu^{(i)}}) S^\pm_w.
\end{equation}
\eco
\begin{proof}
The condition on the element $w$ implies that for every 
$\la \in P^+_{\SS(w)}$ there exists $\la_{\min} \in P^+_{\SS(w)}$
such that
\[
(\la + \Kl(w)) \cap P^+_{\SS(w)} = \la_{\min} + 
\big( \Nset \mu^{(1)} \oplus \ldots \oplus \Nset \mu^{(m(w))} \big).
\]
We can then choose $\Lambda^+_w$ to be the set of all
such elements $\la_{\min}$. We have  
\[
P^+_{\SS(w)} = \Lambda^+_w \oplus 
\big( \Nset \mu^{(1)} \oplus \ldots \oplus \Nset \mu^{(m(w))} \big). 
\]
It follows from \eqref{sum2} that 
\[
J_{w, {\bf{1}} } = \sum_{i=1}^{m(w)}
(1 - d^\pm_{w, \mu^{(i)}}) N^\pm_w,
\]
which implies \eqref{genn}.
\end{proof}
\subsection{A classification of the 
prime elements of $S^\pm_w$}
\label{5.7}
As another application of \thref{normal1} 
we obtain a classification
of all inhomogeneous prime elements of the 
algebras $S^\pm_w$. When this is combined 
with \thref{normal1} (ii), it gives a classification 
of all prime elements of the algebras $S^\pm_w$.
The (anti)isomorphisms from \thref{isom} give an analogous 
classification of all prime elements of the 
algebras $\UU^w_\pm$. The formulation of the latter 
is straightforward and will not be stated separately.
From these results we obtain a classification 
of all normal elements of the algebras $S^\pm_w$, 
which via the (anti)isomorphism of \thref{isom} gives 
a classification of all normal elements of the 
algebras $\UU^w_\pm$.

For $n \in \Nset$ denote
\begin{multline}
\label{irr}
\Irr_n^0(\KK) = \{ f(x_1, \ldots, x_n)  \in \KK[x_1, \ldots, x_n] \mid  f(x) \; \; 
\mbox{is irreducible} 
\\
\mbox{and} \; \;  f(0, \ldots, 0) = 1 \}.
\end{multline}
Recall the definition \eqref{kl} of the lattice 
$\Kl (w) \subset P_{\SS(w)}$ and the notation 
$m(w) = \dim \ker(w + 1)$.
Recall from the previous subsection that 
$\mu^{(1)}, \ldots, \mu^{(m(w))}$ 
denotes a fixed basis of $\Kl(w)$.

For each $f(x_1, \ldots, x_{m(w)}) \in \Irr_{m(w)}^0(\KK)$ there exists a unique
$\la_f \in P^+_{\SS(w)}$ such that 
\[
d^\pm_{w, \la_f} f( d^\pm_{w, \mu^{(1)}}, \ldots, d^\pm_{w, \mu^{(m(w))}}) 
= \sum_{ \la' \in P^+_{\SS(w)}} p_{\la'} d^\pm_{w, \la'} \in N^\pm_{w}
\]
and
\begin{equation}
\label{cup}
\cap \{ \Supp \la' \mid \la' \in P^+_{\SS(w)}, p_{\la'} \neq 0 \} = \emptyset,
\end{equation}
recall \eqref{supp}. We denote 
\begin{equation}
\label{wh}
\wh{f} = d^\pm_{w, \la_f} f( d^\pm_{w, \mu^{(1)}}, \ldots, d^\pm_{w, \mu^{(m(w))}}) 
\in N^\pm_w.
\end{equation}
Since the second factor above belongs to the center of $S^\pm_w[(M^\pm_w)^{-1}]$,
we have from \eqref{ddcomm} that $\wh{f} \in S^\pm_w$ is normal and
more precisely:
\begin{equation}
\label{whcomm}
\wh{f} s = q^{- \lcor (w + 1 ) \la_{f} , \ga \rcor } 
s \wh{f}, \quad \forall 
s \in (S_w^\pm)_{-\ga, 0}, \ga \in Q_{\SS(w)}.
\end{equation}
It follows from \eqref{cup} that 
\begin{equation}
\label{divide}
(d^\pm_{w, \om_i})^{-1} \wh{f} \notin S^\pm_w, \quad
\forall i \in \SS(w),
\end{equation}
because of \eqref{dprod2}, \thref{freeS} (ii) 
and \eqref{basisNM}.

The next theorem contains our classification 
result for the inhomogeneous prime elements of $S^\pm_w$. 
Equivalently, it provides an explicit description 
of the height one prime ideals of 
$S^\pm_w$ which are not $\Tset^r$-invariant 
with respect to \eqref{Tract0}. The latter is an 
example of a case in which the formula 
from \prref{idealsU} (i) 
for the prime ideals in 
$\Spec_{ \{ 0 \} } S^\pm_w$ simplifies. 

\bth{inhompr} Assume that $\KK$ is an arbitrary base field, $q \in \KK^*$ 
is not a root of unity and $w \in W$. Then every inhomogeneous
prime element of $S^\pm_w$ is of the form $p \wh{f}$,
for some $f \in \Irr^0_{m(w)}(\KK)$ and $p \in \KK^*$, cf. 
\eqref{mw} and \eqref{irr}.
All such elements are distinct.

The height one prime ideals of $S^\pm_w$ which are not $\Tset^r$-invariant
with respect to \eqref{Tract0} have the form $\wh{f} S^\pm_w$, 
for some $f \in \Irr^0_{m(w)}(\KK)$. All such ideals are distinct.
\eth

The special case of \thref{inhompr} for the algebras of 
quantum matrices is due to Launois and Lenagan \cite{LL}.
\\ \hfill \\
{\em{Proof of \thref{inhompr}}}. 
All height one prime ideals of $S^\pm_w$ are 
the $\Tset^r$-invariant height one prime ideals 
with respect to \eqref{Tract0} (which as mentioned in \S \ref{5.2}
are the ideals $I^\pm_w(s_i)$ for $i \in \SS(w)$)
and the height one prime ideals in $\Spec_{ \{0\} } S^\pm_w$. 
The latter family consists of ideals which are not 
$\Tset^r$-invariant with respect to \eqref{Tract0}.  
By \cite[Theorem 5.3]{GK} every 
height one prime ideal in $\Spec_{ \{ 0 \} } S^\pm_w$ 
is of the form
\begin{equation}
\label{JJ}
(S^\pm_w[ (M^\pm_w)^{-1} ] . J^0) \cap S^\pm_w 
\end{equation}
for some height one prime ideal $J^0$ of 
$Z( S^\pm_w[ (M^\pm_w)^{-1} ])=A_w^\pm$, 
cf. \eqref{primUw}, and all such 
ideals are distinct. We have 
\[
Z( S^\pm_w[ (M^\pm_w)^{-1} ]) = A_w^\pm \cong
\KK[x^{\pm 1}_1, \ldots, x^{\pm 1}_{m(w)}], \quad 
d^\pm_{w, \mu^{(i)}} \mt x_i, i = 1, \ldots, m(w).
\]
Each height one prime ideal of the Laurent polynomial
ring $\KK[x^{\pm 1}_1, \ldots, x^{\pm 1}_{m(w)}]$
is generated by a prime element. Each prime 
element of $\KK[x^{\pm 1}_1, \ldots, x^{\pm 1}_{m(w)}]$
is uniquely represented as a product 
$u f(x_1, \ldots, x_{m(w)})$ where 
$u$ is a unit of the Laurent polynomial ring
(i.e. a nonzero Laurent monomial) and 
$f(x_1, \ldots, x_{m(w)}) \in \Irr^0_{m(w)}(\KK)$.
Therefore each height one prime ideal of $S^\pm_w$ 
has the form
\[
J(f) = \{ s \in S^\pm_w \mid \exists \la \in P^+_{\SS(w)} 
\; \; \mbox{such that} \; \; 
d^\pm_{w, \la} s \in S^\pm_w \wh{f} \}
\]
for some $f \in \Irr^0_{m(w)}(\KK)$,
because \eqref{wh} implies
$S^\pm_w[ (M^\pm_w)^{-1}] f = 
S^\pm_w[ (M^\pm_w)^{-1}] \wh{f}$. 
We claim that 
\[
J(f) = S^\pm_w \wh{f}, \quad 
\forall f \in \Irr^0_{w(w)}(\KK).
\]
To prove this, all we need to show is that 
for $s, s' \in S^\pm_w$,
\begin{equation}
\label{toshow}
d^\pm_{w, \om_i} s = s' \wh{f} \quad
\Rightarrow \quad
s \in S^\pm_w \wh{f}.
\end{equation}
Since $d^\pm_{w, \om_i} \in S^\pm_w$ is prime (see \thref{normal1} (ii)), 
if $d^\pm_{w, \om_i} s = \wh{f} s'$ then 
either $s'$ or $\wh{f}$ is a multiple 
of $d^\pm_{w, \om_i}$. The second is 
not possible because of \eqref{divide}.   
Hence $s' = d^\pm_{w, \om_i} s''$
for some $s'' \in S^\pm_w$ and thus 
$s = s'' \wh{f} \in S^\pm_w \wh{f}$.

Therefore all height one prime ideals have the 
form 
\[
J(f) = S^\pm_w \wh{f} 
\]
for some $f \in \Irr^0_{m(w)}(\KK)$ and all such 
ideals are distinct. This implies that every inhomogeneous
prime element of $S^\pm_w$ is of the form $p \wh{f}$,
for some $f \in \Irr^0_{m(w)}(\KK)$ and $p \in \KK^*$, and 
these elements are distinct. To show that the ideals $J(f)$
are completely prime for all $f \in \Irr^0_{m(w)}(\KK)$,
one either applies 
\cite[Theorem 3.7]{LLR} to conclude that all height 
one prime ideals of $S^\pm_w \cong \UU^w_\pm$ are 
generated by prime elements, or \eqref{whcomm} and the fact 
\cite[Theorem 2.1]{GLet0} that 
all prime ideals of $S^\pm_w \cong \UU^w_\mp$ are 
completely prime.
\qed
\\ \hfill \\
Next we use \thref{inhompr} to obtain a classification 
of all normal elements of the algebras $S^\pm_w$. 
Let $R$ be a unique factorization domain. Let $C$ 
be the set of all elements of $R$ which are not 
divisible by a prime element. Chatters proved 
\cite[Proposition 2.1]{Ch} that each nonzero
element of $R$ can be represented as a product 
$c p_1 \ldots p_n$ where $c \in C$ and $p_1, \ldots, p_n$
is a sequence of not necessarily distinct prime elements
of $R$. This factorization takes a particularly simple form 
in the case of normal elements of $R$. 
The proof of the second part of the following proposition 
was communicated to us by Tom Lenagan \cite{Le}.

\bpr{factor} Let $R$ be a noetherian unique factorization domain. Then:

(i) \cite{ChJ} Every nonzero normal element of $R$ can be represented as a 
product $u p_1 \ldots p_n$ where $u \in R$ is a unit, 
$p_1, \ldots, p_n$ is a sequence of not necessarily 
distinct prime elements of $R$ and $n \in \Nset$.

(ii) Assume in addition that a $\KK$-torus $T$ acts 
rationally on $R$ by algebra automorphisms
for an infinite field $\KK$. Then every nonzero normal 
element of $R$ which is a $T$-eigenvector
can be represented as a product $u p_1 \ldots p_n$ 
where $u \in R$ is a unit which is a $T$-eigenvector, 
$p_1, \ldots, p_n$ is a sequence of not necessarily distinct 
prime elements of $R$ which are $T$-eigenvectors and $n \in \Nset$.
\epr
\begin{proof} (i) The proof of this fact was sketched in 
\cite[p. 24]{ChJ}. We give a version of this proof for completeness.
Assume that $x \in R$ is a nonzero normal element. 
By Chatters' result \cite[Proposition 2.1]{Ch}
\[
x = c p_1 \ldots p_n
\]
for some $c \in C$ and prime elements $p_1, \ldots, p_n$. 
Since $R$ is a domain, $c$ is normal. Indeed, there exist
$\psi, \psi_1, \ldots, \psi_n \in \Aut(R)$ such that 
$x r = \psi(r) x$ and $p_k r = \psi_k(r) p_k$ for all $r \in R$, 
$k =1, \ldots, n$. Then $c r = (\psi \psi_n^{-1} \ldots \psi_1^{-1}(r)) c$, 
$\forall r \in R$. We are left with proving that $c$ has to be a unit.
Assume that this is not the case. Let $P$ be a minimal prime over $c R$. 
By the principle ideal theorem \cite[Theorem 4.1.11]{McR}, the height of 
$P$ is equal to 0 or 1. The height of $P$ cannot be equal to 0 since $c \neq 0$
and $R$ is a domain. Using again the fact that $R$ is a unique
factorization domain we obtain that $P = p R$ for some prime 
element $p$ of $R$. Since $c R \subseteq p R$, $c$ is a multiple of 
$p$ which contradicts the fact that $c \in C$.

(ii) Let $x \in R$ be a nonzero normal element which is a 
$T$-eigenvector. If $x$ is a unit, the statement is clear.
If this is not the case, let $P$ be a minimal prime 
of $R$ over $x R$. By the result of Brown and Goodearl 
\cite[Proposition II.2.9]{BG} all $T$-primes of $R$ are 
prime and thus $P$ is a $T$-prime ideal. 
Using the principle ideal theorem \cite[Theorem 4.1.11]{McR} 
as in part (i),
we find that $P = p R$ for some prime element $p \in R$ which is
a $T$-eigenvector. Then 
\[
x = x' p
\]
for some normal element $x' \in R$ which is a $T$-eigenvector.
The statement follows by induction using the fact that $R$ is noetherian.
\end{proof}

\prref{factor} (ii) and \eqref{multd} imply that the first part of 
\thref{normal1} follows from its second part. (The condition that 
$\KK$ is infinite follows from the assumption that $q \in \KK^*$ 
is not a root of unity.) At the same time one 
should point out that our proof of the second part of \thref{normal1} 
relies heavily on its first part. Thus one needs an independent
proof of \thref{normal1} (ii) in order to obtain a second proof 
of \thref{normal1} (i) via \prref{factor} (ii).

The next theorem uses \prref{factor} (ii), and Theorems 
\ref{tnormal1} (ii) and \ref{tinhompr} to 
extend the classification of \thref{normal1} (i) 
to a classification of all normal elements of the algebras $S_w^\pm$. 
Recall the definition \eqref{kl} of the lattice $\Kl(w)$.

\bth{allnorma} For all base fields $\KK$, $q \in \KK^*$ 
not a root of unity and $w \in W$, the normal elements of 
$S^\pm_w$ are precisely the elements of the form
\begin{equation}
\label{normal-form}
\sum_{\mu \in \Kl(w) \cap(- \la + P_{\SS(w)}^+)} p_\mu 
d^\pm_{w, \la}
\end{equation}
for some $\la \in P_{\SS(w)}^+$ and a family of scalars
$p_\mu \in \KK$ (only finitely many of which are nonzero).
\eth

Another way to formulate \thref{allnorma} is to say that all normal 
elements of the algebras $S^\pm_w$ are linear combinations of 
homogeneous normal elements (classified in \thref{normal1} (i)).
The special case of \thref{allnorma} for $w = w_0$ 
(the longest element of $W$), $\KK$ of characteristic $0$ and 
$q$ transcendental over $\Qset$ is due to Caldero \cite{Ca}.
\\ \hfill \\
\noindent
{\em{Proof of \thref{allnorma}.}} It follows from \eqref{ddcomm} 
and the definition \eqref{kl} of the lattice $\Kl(w)$ that 
all elements of the form \eqref{normal-form} are normal:
\[
\Big(
\sum_{\mu \in \Kl(w) \cap(- \la + P_{\SS(w)}^+)} p_\mu 
d^\pm_{w, \la}
\Big)
s = q^{\lcor (w+1) \la, \ga \rcor } s 
\Big(
\sum_{\mu \in \Kl(w) \cap(- \la + P_{\SS(w)}^+)} p_\mu 
d^\pm_{w, \la}
\Big)
\]
for all $s \in (S^\pm_w)_{- \ga, 0}$, $\ga \in Q_{\SS(w)}$.
The key point is to prove is that all normal elements of $S^\pm_w$
have this form.

\leref{LSmult} implies that the units of the algebras $\UU^w_\pm$ are precisely
the nonzero scalars in them. Because of the (anti)isomorphisms from \thref{isom}, the 
same is true for the algebras $S^\pm_w$. \prref{factor} (i), \thref{inhompr} 
and \eqref{multd} imply that each normal element of $S^\pm_w$ 
is a linear combination of $d_{w, \la}^\pm$
for some $\la \in P_{\SS(w)}$ (i.e. belongs to the subalgebras 
$N_w^\pm$). Thus each normal element of $S^\pm_w$ has the form
\[
x = \sum_{\mu \in (- \la + P_{\SS(w)}^+)} p_\mu
d^\pm_{w, \la}
\] 
for a finite family of scalars $p_\mu \in \KK$. It follows from 
\thref{inhompr}, \eqref{whcomm} and \thref{normal1} (ii) that all 
prime elements of $S^\pm_w$ are $P$-normal. Since all 
units of $S^\pm_w$ are scalars, \prref{factor} (i) implies 
that all normal elements of $S^\pm_w$ are $P$-normal. 
Thus there exits $\de \in P$ such that 
\begin{equation}
\label{ee1}
x s = s  \Big( 
\sum_{\mu \in (- \la + P_{\SS(w)}^+)} p_\mu q^{\lcor (w+1) \de, \ga \rcor } 
d^\pm_{w, \la}
\Big) 
\end{equation}
for all $s \in (S^\pm_w)_{- \ga, 0}$, $\ga \in Q_{\SS(w)}$.
At the same time, from \eqref{ddcomm} we have
\begin{equation}
\label{ee2}
x s = s  \Big( 
\sum_{\mu \in (- \la + P_{\SS(w)}^+)} p_\mu q^{\lcor (w+1)(\la+ \mu), \ga \rcor } 
d^\pm_{w, \la}
\Big) 
\end{equation}
for all $s \in (S^\pm_w)_{- \ga, 0}$, $\ga \in Q_{\SS(w)}$. 
By \eqref{gradS+-} the span of all $\ga$ such that 
$(S^\pm_w)_{- \ga, 0} \neq 0$ is precisely $Q_{\SS(w)}$.   
Since $\{ d^\pm_{w, \la} \}_{\la \in P^+_{\SS(w)}}$
are linearly independent \eqref{ee1} and \eqref{ee2}
imply that, if $p_\mu \neq 0$ then 
$\lcor (w+1) \mu, \ga \rcor = 0$ for all $\ga \in Q_{\SS(w)}$,
i.e. $\mu \in \Kl(w)$.
Thus all normal elements of $S^\pm_w$ have the form 
\eqref{normal-form}.
\qed

The part of the proof of \thref{allnorma} that 
each normal element of $S^\pm_w$
is a linear combination of homogeneous normal 
elements also follows from the following general 
fact for normal elements in graded domains.

\bpr{grdom} Assume that $R$ is a $\Zset^m$-graded domain over a field 
$\KK$ and that $y = \sum_{{\bf{n}} \in \Zset^m} y _{\bf{n}} \in R$ 
is a normal element such that $y_{\bf{n}} \in R_{\bf{n}}$. 
Then all $y_{\bf{n}}$ are normal and there exists a graded 
automorphism $\phi$ of $R$ such that 
\[
y_{\bf{n}} x = \phi(x) y_{\bf{n}}, \quad \forall x \in R, 
\, \, {\bf{n}} \in \Zset^m.
\]
In particular $y x = \phi(x) y$, $\forall x \in R$.
\epr

For the proof of \prref{grdom} we will need the following lemma.
\ble{gr-n} Assume that $R$ is a $\Zset^m$-graded domain over a field 
$\KK$. Then for every normal element $y \in R$ 
\[
y R_{\bf{n}} = R_{\bf{n}} y, \quad \forall {\bf{n}} \in \Zset^m.
\]
\ele
\begin{proof}
For simplicity we will restrict to the case $m=1$, leaving the 
general case to the reader. The statement is clear for $y=0$. 
Let $y \neq 0$.
Write $y = y_i + \ldots + y_j$ where $i \leq j$, $y_k \in R_k$, 
$y_i \neq 0$, $y_j \neq 0$. 
Let $x \in R_n$, $n \in \Zset$. Then 
\begin{equation}
\label{nn1}
yx = (\sum_{k= n_1}^{n_2} z_k )y
\end{equation}  
for some $z_k \in R_k$, $z_{n_1} \neq 0$, $z_{n_2} \neq 0$.
Since $R$ is a domain, the lowest degree term in 
the left hand side of \eqref{nn1} is $y_i x$ and sits in degree $i+n$. 
The lowest degree term in the right hand side of \eqref{nn1} 
is $z_{n_1} y_i$ and sits in degree $i+n_1$. Therefore 
$n_1=n$. Analogously, comparing the highest degree terms 
of the two sides of \eqref{nn1} leads to $n_2= n$. 
Therefore $n_1=n_2=n$, which proves the lemma.
\end{proof}
\noindent
{\em{Proof of \prref{grdom}.}} The statement is obvious for 
$y=0$. Let $y \neq 0$.
There exists an automorphism $\phi \in \Aut(R)$ such that 
\begin{equation}
\label{eq2}
y x = \phi(x) y, \quad \forall x \in R.
\end{equation}
\leref{gr-n} implies that $\phi$ is a graded automorphism, 
i.e. $\phi(R_{\bf{k}}) = R_{\bf{k}}$, 
$\forall {\bf{k}} \in \Zset^m$. We substitute 
$y = \sum_{ {\bf{n}} \in \Zset^m} y_{\bf{n}}$ 
in \eqref{eq2} and take $x \in R_{\bf{k}}$. 
Equating the components in degree ${\bf{n}} + {\bf{k}}$ 
and using the graded property of $\phi$, leads to 
\[
y_{\bf{n}} x = \phi(x) y_{\bf{n}}, \quad 
\forall {\bf{n}}, {\bf{k}} \in \Zset, \, x \in R_{\bf{k}},
\] 
which completes the proof of the proposition.
\qed
\sectionnew{Module structure of $R_{\bfw}$ 
over their subalgebras generated by Joseph's normal elements}
\lb{Module}
\subsection{Statement of the freeness result}
\label{6.1}
In this section we analyze the structure of $R_{\bfw}$ as a module 
over its subalgebra generated by the Joseph set 
of normal elements $E_{\bfw}^{\pm 1}$, recall \eqref{Ew}.
We prove that $S_{\bfw}$ is a free module over its subalgebra 
generated by the normal elements $y_{\om_i}$, $i=1, \ldots, r$.
We use this to prove that $R_{\bfw}$ is a free module over 
its subalgebra generated by the set $E_{\bfw}^{\pm 1}$. This result 
will be the main tool in classifying $\Max R_q[G]$ in the next 
section, which in turn will be used  to answer affirmatively a question 
of Goodearl and Zhang \cite{GZ}, that all 
maximal ideals of $R_q[G]$ have finite codimension.
The latter will be applied in the last section to prove that
$R_q[G]$ satisfies a stronger property than catenarity, namely
that all maximal chains of prime ideals in $R_q[G]$ have the same 
length, equal to $\GKdim R_q[G]$.

Denote by $L_{\bfw}$ the subalgebra of $R_{\bfw}$ generated by 
$E_{\bfw}^{\pm1}$, i.e. the subalgebra of $R_{\bfw}$ spanned
by $\{ c^+_{w_+, \mu_1} c^-_{w_-, \mu_2} \mid \mu_1, \mu_2 \in P \}$.
Then:

\bth{free} For an arbitrary base field $\KK$ and $q \in \KK^*$
which is not a root of unity, the algebra $R_{\bfw}$ is a free 
left and right $L_{\bfw}$-module in which $L_{\bfw}$ is a direct 
summand viewed as a module over itself.
\eth

Recall \eqref{y}. Denote by $N'_{\bfw}$ the subalgebra of 
$S_{\bfw}$, generated by $y_{\om_i}$, 
$i \in \{ 1, \ldots, r \}$.
We will prove the following result and deduce from 
it \thref{free}:

\bth{free2} For an arbitrary base field $\KK$ and $q \in \KK^*$
which is not a root of unity, the algebra $S_{\bfw}$ is a 
free left and right $N'_{\bfw}$-module in which $N'_{\bfw}$ 
is a direct summand viewed as a module over itself.
\eth

Explicit versions of the decompositions in Theorems 
\ref{tfree} and \ref{tfree2} will be obtained in 
Theorems \ref{tf1} and \ref{tf2}. 
Similarly to the algebras $S^\pm_w$, the centers 
of $R_{\bfw}$ and $S_{\bfw}$ are much smaller 
than the subalgebras generated by the homogeneous 
normal elements of the Joseph's set $E_{\bfw}$ and 
the multiplicative subset of $S_{\bfw}$ generated 
by $y_{\om_i}$, $i \in \{ 1, \ldots, r \}$, respectively. 
Because of this, one obtains stronger results 
when considering the module structure of $R_{\bfw}$ 
and $S_{\bfw}$ over their subalgebras $L_{\bfw}$ and 
$N'_{\bfw}$, rather than their centers.
It is this type of results that are 
eventually applicable to classify $\Max R_q[G]$.
\subsection{A $Q \times Q$-filtration of $S_{\bfw}$}
\label{6.2} 
The algebra $S_{\bfw}$ is only $Q$-graded by \eqref{gradS}.
In this subsection we prove that it has a nontrivial 
$Q \times Q$-filtration which reveals a 
richer structure than the grading. This and the freeness result 
of Section \ref{free1} will be the main tools in the proofs 
of Theorems \ref{tfree} and \ref{tfree2}.
For $w \in W$, denote
\begin{equation}
\label{Qpm}
Q^+_w = 
\sum_{\be \in w(\De_+) \cap \De_-} \Nset \be \subset Q^+ \cap Q_{\SS(w)}.
\end{equation}
Recall that 
\begin{equation}
\label{Splus}
S_{\bfw} = \bigoplus_{(\ga_+, \ga_-) \in Q^+_{w_+} \times Q^+_{w_-} }  
(S^+_{w_+})_{-\ga_+, 0} (S^-_{w_-})_{\ga_-, 0}
\end{equation}
and
\begin{equation}
\label{Stimes}
(S^+_{w_+})_{-\ga_+, 0} (S^-_{w_-})_{\ga_-, 0} 
\cong 
(S^+_{w_+})_{-\ga_+, 0} \otimes_\KK (S^-_{w_-})_{\ga_-, 0}
\end{equation}
as $\KK$-vector spaces (via the multiplication map),
see \eqref{Sisom} and \eqref{gradS+-}. 
For $(\ga_+, \ga_-) \in Q^+_{w_+} \times Q^+_{w_-}$
denote
\begin{equation}
\label{Sgaga}
(S_{\bfw})^{(\ga_+, \ga_-)}= 
(S^+_{w_+})_{-\ga_+, 0} (S^-_{w_-})_{\ga_-, 0}.
\end{equation}
Consider the induced partial order on $Q^+_{w_+} \times Q^+_{w_-}$
from the product partial order \eqref{po} of $Q\times Q$. Thus, 
for two pairs  
$(\ga'_+, \ga'_-)$, $(\ga_+, \ga_-) \in Q^+_{w_+} \times Q^+_{w_-}$
we set $(\ga'_+, \ga'_-) \prec (\ga_+, \ga_-)$ if $\ga'_+ < \ga_+$ and 
$\ga'_- < \ga_-$, i.e. if there exist $\be_\pm \in Q^+ \backslash \{0\}$ such that 
$\ga_+ = \ga'_+ + \be_+$ and $\ga_- = \ga'_- + \be_-$. 

For $(\ga_+, \ga_-) \in Q^+_{w_+} \times Q^+_{w_-}$ define
\begin{align*}
(S_{\bfw})^{\prec (\ga_+, \ga_-)} 
&= \bigoplus_{ (\ga'_+, \ga'_-) \in Q^+_{w_+} \times Q^+_{w_-}, 
(\ga'_+, \ga'_-) \prec (\ga_+, \ga_-) } 
(S_{\bfw})^{(\ga_+, \ga_-)}
\\
&= \bigoplus_{\ga'_\pm \in Q^+_{w_\pm}, \ga'_\pm < \ga_\pm}  
(S^+_{w_+})_{-\ga'_+, 0} (S^-_{w_-})_{\ga'_-, 0}
\end{align*}
and
\begin{equation}
\label{Sfilt}
(S_{\bfw})^{\preceq (\ga_+, \ga_-)} = (S_{\bfw})^{(\ga_+, \ga_-)} \oplus
(S_{\bfw})^{\prec (\ga_+, \ga_-)}.
\end{equation}
Since the algebras $S^\pm_{w_\pm}$ are $Q$-graded we have
\begin{equation}
\label{mult-gr}
s_+ (S_{\bfw})^{\prec (\ga'_+, \ga'_-)} s_- \subseteq
(S_{\bfw})^{\prec (\ga'_+ \footnote{}+ \ga_+, \ga'_- + \ga_-)}, 
\end{equation}
for all $s_\pm \in (S^\pm_{w_\pm})_{\ga_\pm}$, 
$(\ga'_+, \ga'_-)$, $(\ga_+, \ga_-) \in Q^+_{w_+} \times Q^+_{w_-}$.

Consider the $Q^+_{w_+} \times Q^+_{w_-}$ (exhaustive ascending) 
filtration of the space 
$S_{\bfw}$ by the subspaces $(S_{\bfw})^{\preceq (\ga_+, \ga_-)}$, 
$(\ga_+, \ga_-) \in Q^+_{w_+} \times Q^+_{w_-}$. The next result 
proves that this is an algebra filtration and computes the 
structure of the associate graded algebra. 

\bpr{Sbfw} For all 
$(\ga_+, \ga_-), (\ga'_+, \ga_-) \in Q^+_{w_+} \times Q^+_{w_-}$, 
$s_\pm \in (S^\pm_{w_\pm})_{\mp \ga_\pm}$ and
$s'_\pm \in (S^\pm_{w_\pm})_{\mp \ga'_\pm}$, we have
\[
\big(s_+ s_- \big). \big(s'_+ s'_- \big) = 
q^{-\lcor \ga_-, \ga'_+ \rcor} 
\big( (s_+ s'_+) (s_- s'_-) \big)
\mod (S_{\bfw})^{ \prec (\ga_+ + \ga'_+, \ga_- + \ga'_-) }.
\]
\epr
Note that in the setting of the proposition
\begin{align*}
&s_+ s_- \in (S_{\bfw})^{(\ga_+ ,\ga_-)},  \; \; 
s'_+ s'_- \in (S_{\bfw})^{(\ga'_+ ,\ga'_-)} \; \; 
\mbox{and}
\\
&(s_+ s'_+) (s_- s'_-) \in 
(S_{\bfw})^{(\ga_+ + \ga'_+,\ga_- + \ga'_-)}. 
\end{align*}
We will identify
\[
(S_{\bfw})^{\preceq (\ga_+, \ga_-)} / (S_{\bfw})^{\prec (\ga_+, \ga_-)}
\cong (S_{\bfw})^{(\ga_+, \ga_-)} \; \; \mbox{for} 
\; \; (\ga_+, \ga_-) \in Q^+_{w_+} \times Q^+_{w_-},
\]
(cf. \eqref{Sfilt}) and 
\begin{equation}
\label{grS}
\gr S_{\bfw} \cong \bigoplus_{ (\ga_+, \ga_-) \in Q^+_{w_+} \times Q^+_{w_-}}
(S_{\bfw})^{(\ga_+, \ga_-)},
\end{equation}
(cf. \eqref{Sgaga}).
Denote the multiplication in $\gr S_{\bfw}$ by $\odot$.

\bco{grmult} Under the identification of \eqref{grS} the 
multiplication in $\gr S_{\bfw}$ is given by 
\[
\big(s_+ s_- \big) \odot \big(s'_+ s'_- \big) = 
q^{-\lcor \ga_-, \ga'_+ \rcor} 
\big( (s_+ s'_+) (s_- s'_-) \big),
\]
for all
$(\ga_+, \ga_-), (\ga'_+, \ga_-) \in Q^+_{w_+} \times Q^+_{w_-}$, 
$s_\pm \in (S^\pm_{w_\pm})_{\mp \ga_\pm}$ and
$s'_\pm \in (S^\pm_{w_\pm})_{\mp \ga'_\pm}$.
\eco
\noindent
{\em{Proof of \prref{Sbfw}.}} It follows from \eqref{commRR} 
that 
\[
s_- s'_+ = q^{-\lcor \ga_-, \ga'_+ \rcor} s'_+ s_-  
+ \sum_{i=1}^k (s'_+)^{(i)} (s_-)^{(i)},  
\]
for some $(s'_+)^{(i)} \in (S^+_{w_+})_{ - \ga_+^{(i)}}$, 
$(s_-)^{(i)} \in (S^-_{w_-})_{ \ga_-^{(i)}}$,
$\ga_+^{(i)} \in Q^+_{w_+}$, $\ga_+^{(i)} < \ga'_+$, 
$\ga_-^{(i)} \in Q^-_{w_-}$, $\ga_-^{(i)} < \ga_-$, 
$i = 1, \ldots, k$. Therefore
\begin{equation}
\label{au}
s_- s'_+ = q^{-\lcor \ga_-, \ga'_+ \rcor} s'_+ s_- 
\mod (S_{\bfw})^{ \prec ( \ga'_+, \ga_-) }.
\end{equation}
Multiplying \eqref{au} on the left by $s_+$ and on 
the right by $s'_-$, and using \eqref{mult-gr} implies the 
statement of the proposition. 
\qed
\subsection{The action of $\gr N'_{\bfw}$ on $\gr S_{\bfw}$}
\label{6.3}
Next, we apply the results from the previous subsection to 
the $N'_{\bfw}$-module structure of $S_{\bfw}$.

First, observe from \eqref{x} that for all $i \in \II({\bfw})$ 
the image of $x_{\om_i}$ in $\wh{R}_{\bfw}$ is equal to 
$c^+_{1,\om_i} c^-_{1, \om_i} = c^+_{w_+,\om_i} c^-_{w_-, \om_i}$. 
Applying \eqref{n3}, we get 
\begin{equation}
\label{y-1}
y_{\om_i} = 
(c^+_{w_+,\om_i})^{-1} (c^-_{w_-, \om_i})^{-1} c^+_{w_+,\om_i} c^-_{w_-, \om_i} =
1, \quad \forall i \in \II({\bfw}).
\end{equation}
Recall the definition \eqref{dd} of the elements 
$d^\pm_{w_\pm, \la} \in (S^\pm_{w_\pm})_{\pm(w_\pm-1)\la, 0}$.
We have 
\begin{equation}
\label{dlala}
d^+_{w_+, \la} d^-_{w_-, \la} \in (S_{\bfw})^{( (1-w_+) \la, (1-w_-) \la)}. 
\end{equation}
Eqs. \eqref{x} and \eqref{n3} imply that  
\begin{align}
\label{d-filt}
y_{\om_i} &= 
(c^+_{w_+,\om_i})^{-1} (c^-_{w_-, \om_i})^{-1} c^+_{1,\om_i} c^-_{1, \om_i}
= 
q^{\lcor \om_i, (1-w_-) \om_i \rcor }
(c^+_{w,\om_i})^{-1} c^+_{1,\om_i} (c^-_{w, \om_i})^{-1} c^-_{1, \om_i}
\\
&=
q^{\lcor \om_i, (1-w_-) \om_i \rcor }
d^+_{w_+, \om_i} d^-_{w_-, \om_i} \mod
(S_{\bfw})^{\prec ( (1-w_+) \om_i, (1-w_-) \om_i )}, 
\nn
\end{align}
for all $i \in \SS({\bfw})$.

Recall the definition \eqref{yla} of the elements $y_\la \in N'_{\bfw}$,
$\la \in P^+$. Applying repeatedly \prref{Sbfw} and using the fact that 
$d^\pm_{w_\pm, \om_i} \in S^\pm_{w_\pm}$ are 
$P$-normal, we obtain:

\bco{exp} For every $\la \in P^+_{\SS(\bfw)}$, 
$s_\pm \in (S^\pm_{w_\pm})_{\mp \ga_\pm, 0}$ 
there exists $m \in \Zset$ such that
\[
y_\la (s_+ s_-) = q^{m} (s_+ d^+_{w_+,\la}) (d^-_{w_-, \la} s_-) 
\mod (S_{\bfw})^{\prec (\ga_+ + (1-w_+)\la, \ga_- + (1-w_-)\la)}.  
\] 
\eco
\noindent
Note that in the setting of \coref{exp}, 
\[
(s_+ d^+_{w_+,\la}) (d^-_{w_-, \la} s_-) 
\in (S_{\bfw})^{(\ga_+ + (1-w_+)\la, \ga_- + (1-w_-)\la)}.
\]
Setting $s_+=1$, $s_-=1$, we obtain that for all 
$\la \in P^+_{\SS({\bfw})}$
\begin{equation}
\label{y-mod}
y_\la = q^{m_\la} d^+_{w_+, \la} d^-_{w_-, \la} 
\mod (S_{\bfw})^{\prec ((1-w_+)\la, (1-w_-)\la)}, 
\end{equation}
for some $m_\la \in \Zset$ and 
$d^+_{w_+, \la} d^-_{w_-, \la} \in (S_{\bfw})^{((1-w_+)\la, (1-w_-)\la)}$.

Denote
\[
\Ga_{\bfw} = \{ ( (1-w_+)\la, (1-w_-) \la ) \mid \la \in 
P^+_{\SS({\bfw})} \}. 
\]
Eq. \eqref{y-mod} implies
\begin{equation}
\label{Ngr}
\gr N'_{\bfw} \cong \bigoplus_{ (\ga_+, \ga_-) \in \Ga_{\bfw} }
(N'_{\bfw})^{(\ga_+, \ga_-)},
\end{equation}
where for $\la \in P^+_{\SS({\bfw})}$
\begin{multline*}
(N'_{\bfw})^{( (1-w_+)\la , (1-w_-) \la )} =
N'_{\bfw} \cap (S_{\bfw})^{( (1-w_+)\la , (1-w_-) \la )}
\\
= \Span \{ y_\mu \mid \mu \in P_{\SS({\bfw})}, 
\mu - \la \in \ker (1-w_+) \cap \ker (1-w_-) \}
\end{multline*}
in the identification \eqref{grS}.
Denote by $\gr y_\la$ the image of $y_\la$ in $\gr N'_{\bfw}$. 
Eq. \eqref{y-mod} implies that for each $\la \in P_{\SS({\bfw})}^+$ 
there exists $m_\la \in \Zset$ such that
\begin{equation}
\label{gry}
\gr y_\la = q^{m_\la} d^+_{w_+, \la} d^-_{w_-, \la}
\end{equation}
in terms of the identification \eqref{grS}.

For $s \in (S_{\bfw})^{\preceq (\ga_+, \ga_-)}$ 
denote by $\gr s$ its image in $\gr S_{\bfw}$.
\coref{exp} implies that for all $\la \in P^+_{\SS(\bfw)}$, 
$s_\pm \in (S^\pm_{w_\pm})_{\mp \ga_\pm, 0}$ 
$\ga_\pm \in Q^+_{w_\pm}$
there exists $m \in \Zset$ such that
\begin{equation}
\label{sla}
(\gr y_\la) \odot \big( \gr (s_+ s_-) \big) = q^{m} (s_+ d^+_{w_+,\la})
(d^-_{w_-, \la} s_-),
\end{equation}
where in the right hand side we used the identification \eqref{grS}.
\subsection{Structure of the algebras $N'_{\bfw}$ and 
separation of variables for $S_{\bfw}$}
\label{6.4}
Recall that for ${\bfw} = (w_+, w_-) \in W \times W$ 
\[
\SS({\bfw}) =  \SS(w_+) \cup \SS(w_-). 
\]
We have 
\[
\SS({\bfw}) =  \big( \SS(w_+) \cap \SS(w_-) \big) 
\bigsqcup \big( \SS(w_+) \backslash \SS(w_-) \big)
\bigsqcup \big( \SS(w_-) \backslash \SS(w_+) \big)
\]
and the corresponding decomposition
\begin{equation}
\label{decom}
P^+_{\SS({\bfw})} =  
P^+_{ \SS(w_+) \cap \SS(w_-)} 
\bigoplus
P^+_{\SS(w_+) \backslash \SS(w_-)}
\bigoplus 
P^+_{\SS(w_-) \backslash \SS(w_+)}.
\end{equation}
For $\la \in P^+_{\SS({\bfw})}$, denote 
its components 
\begin{equation}
\label{la+-0}
(\la)_0 \in P^+_{ \SS(w_+) \cap \SS(w_-)},
\quad
(\la)_+ \in P^+_{\SS(w_+) \backslash \SS(w_-)},
\quad 
(\la)_- \in P^+_{\SS(w_-) \backslash \SS(w_+)}
\end{equation}
in the decomposition \eqref{decom}. 
For $\mu \in P^+_{\SS(w_\mp) \backslash \SS(w_\pm)}$,
$d^\pm_{w_\pm, \mu}$ is a nonzero scalar by 
\eqref{y-1} and \eqref{multd}. Using this 
and one more time \eqref{multd},
we obtain that for each 
$\la  \in P^+_{\SS({\bfw})}$ there exist
integers $m_\la$ and $m'_\la$ such that
\begin{equation}
\label{ti}
d^+_{w_+, \la} = q^{m_\la} d^+_{w_+,  (\la)_0 + (\la)_+ }
\end{equation} 
and
\begin{equation}
\label{ti2}
d^-_{w_-, \la} = q^{m'_\la} d^-_{ w_-, (\la)_0 + (\la)_-}.
\end{equation} 

It follows from \eqref{yc} that
\begin{equation}
\label{y-qcomm}
y_{\om_i} y_{\om_j} = 
q^{\lcor w_- \om_i, w_+ \om_j \rcor - \lcor w_+ \om_i, w_- \om_j \rcor} 
y_{\om_j} y_{\om_i}, \quad
i, j = 1, \ldots, r.
\end{equation}
The following result describes the structure of the algebra 
$N'_{\bfw}$.

\bpr{Nprime} For all ${\bfw} \in W \times W$ the algebra 
$N'_{\bfw}$ is isomorphic to the quantum affine space algebra
over $\KK$ of dimension $|\SS({\bfw})|$ 
with generators $y_{\om_i}$, $i \in \SS({\bfw})$ and
relations \eqref{y-qcomm}. The set $\{ y_\la \mid 
\la \in P^+_{\SS({\bfw})} \}$ is a $\KK$-basis of $N'_{\bfw}$. 
\epr
\begin{proof} By part (i) of \thref{freeS}
the elements $d^\pm_{w_\pm, \la} \in S^\pm_{w_\pm}$, 
$\la \in P^+_{\SS(w_\pm)}$ are linearly 
independent. Taking \eqref{Stimes} into account, 
we see that 
\begin{equation}
\label{multind}
\{ d^+_{w_+, \la_1} d^-_{w_-, \la_2} \mid
\la_1 \in P_{\SS(w_+)}^+, \la_2 \in P_{\SS(w_-)}^+  
\}
\subset S_{\bfw}
\end{equation}
is a linearly independent set.

Recall the discussion of quantum affine space algebras
from \S \ref{5a.1}. Because the relations \eqref{y-qcomm}
hold, all we need to show is that the ordered monomials in 
$y_{\om_i}$, $i \in \SS({\bfw})$ are linearly
independent. The latter are $y_\la$, $\la \in P_{\SS({\bfw})}^+$
up to a nonzero scalar multiple.  Applying 
\eqref{gry} and
\eqref{ti}--\eqref{ti2}, we obtain that there exist
integers $n_\la$, $n'_\la$ for $\la \in P_{\SS(w_\pm)}^+$ 
such that 
\begin{align}
\label{meq}
&\{ q^{n_\la} \gr y_\la \mid \la \in P_{\SS({\bfw})}^+ \}
\\
\nn
=&\{ q^{n'_\la} d^+_{w_+, \la} 
d^-_{w_-, \la} \mid
\la \in P_{\SS({\bfw})}^+ \}
\\
\nn
=& \{ d^+_{w_+, (\la)_0 + (\la)_+} 
d^-_{w_-, (\la)_0 + (\la)_-} \mid
\la \in P_{\SS({\bfw})}^+ \}
\end{align}
in the identification \eqref{grS}.
Since the third set is a subset of the set in \eqref{multind}, 
the elements $\{ \gr y_\la \mid \la \in P_{\SS({\bfw})}^+ \}$ 
are linearly independent.  
Therefore the elements 
$\{y_\la \mid \la \in P_{\SS({\bfw})}^+ \}$ are linearly independent.
\end{proof}

Denote 
\[
\Om_{\bfw} = \{ (\mu_1, \mu_2) \in P_{\SS(w_+) \cap \SS(w_-)}^+ 
\times P_{\SS(w_+) \cap \SS(w_-)}^+ 
\mid 
\Supp \mu_1 \cap \Supp \mu_2 = \emptyset \},
\]
recall \eqref{supp}.
For a set $Y$ denote by $\Diag(Y)$ the 
diagonal subset of $Y \times Y$. 

\ble{comp} Let ${\bfw} = (w_+, w_-) \in W \times W$. Then:

(i) Each element of  
$P_{\SS(w_+)}^+ \times P_{\SS(w_-)}^+$ can 
be uniquely represented as a sum of an element 
of $\Om_{\bfw}$ and an element of the set
\begin{align*}
&\{ (\la)_0 + (\la)_+, (\la)_0 + (\la)_-) \mid 
\la \in P_{\SS({\bfw})}^+ \}
\\
=& 
\Diag \big( P_{\SS(w_+) \cap \SS(w_-)}^+ \big)
\bigoplus
\big( P^+_{\SS(w_+) \backslash \SS(w_-)} \times
P^+_{\SS(w_-) \backslash \SS(w_+)} \big),
\end{align*}
cf. \eqref{decom}.

(ii) There exist integers 
$ \{ m_{\la_1, \la_2} \mid (\la_1, \la_2) \in P_{\SS(w_+)}^+ \times P^+_{\SS(w_-)}\}$ 
such that the set 
\[
\{ q^{m_{\la_1, \la_2}} d^+_{w_+, \la_1} d^-_{w_-, \la_2} \mid 
(\la_1, \la_2) \in P_{\SS(w_+)}^+ \times P_{\SS(w_-)}^+ \}
\]
coincides with the set 
\[
\big\{ d^+_{w_+, \mu_1} \big(d^+_{w_+, \la}   d^-_{w_-, \la} \big) 
d^-_{w_-, \mu_2}  \; \big| \: 
\la \in P_{\SS({\bfw})}^+,
(\mu_1, \mu_2) \in \Om_{\bfw} 
\big\}.
\]
\ele
\begin{proof} (i) We have 
\[
P^+_{\SS(w_\pm)} =  
P^+_{ \SS(w_+) \cap \SS(w_-)} 
\bigoplus
P^+_{\SS(w_\pm) \backslash \SS(w_\mp)}.
\]
Because of this, the statement of the first part 
is equivalent to 
\[
P_{\SS(w_+) \cap \SS(w_-)}^+ \times 
P_{\SS(w_+) \cap \SS(w_-)}^+ =
\Om_{\bfw} \bigoplus 
\Diag \big( P_{\SS(w_+) \cap \SS(w_-)}^+ \big).
\]
This fact is easy to verify and is left to the reader. 

(ii) Eq.  \eqref{multd} and the second equality in \eqref{meq}  
imply that there exist integers 
$ m_{\la, \mu_1, \mu_2}$, 
$\la \in P_{\SS({\bfw})}^+$, 
$(\mu_1, \mu_2) \in \Om_{\bfw}$ such that 
\begin{align*}
&\big\{ d^+_{w_+, \mu_1} \big(d^+_{w_+, \la}   d^-_{w_-, \la} \big) 
d^-_{w_-, \mu_2}  \; \big| \: 
\la \in P_{\SS({\bfw})}^+,
(\mu_1, \mu_2) \in \Om_{\bfw} 
\big\}
\\
=& 
\{ q^{m_{\la, \mu_1, \mu_2}} d^+_{w_+, (\la)_0 +(\la)_+ + \mu_1} 
d^-_{w_-, (\la)_0 + (\la)_- + \mu_2 } \mid
\la \in  P_{\SS({\bfw})}^+, 
(\mu_1, \mu_2) \in \Om_{\bfw} \}.
\end{align*}
Now the second part follows from the first one.
\end{proof}

Denote $l_\pm = l(w_\pm)$. Fix reduced expressions
$\vec{w}_\pm$ of $w_\pm$.  Recall the definition 
\eqref{De} of the sets $\De(\vec{w}_\pm) \subseteq \Nset^{l_\pm}$. 
Denote 
\begin{multline}
\label{B'w}
B'_{\bfw} = \big\{ 
(\varphi^+_{w_+})^{-1} \big( (X^-)^{{\bf{n}}_+} \big) d^+_{w_+, \mu_1} d^-_{w_-, \mu_2}
(\varphi^-_{w_-})^{-1} \big( (X^+)^{{\bf{n}}_-} \big) \; \big| \; 
\\
{\bf{n}}_+ \in \De(\vec{w}_+), 
{\bf{n}}_- \in \De(\vec{w}_-), 
(\mu_1, \mu_2) \in  \Om_{\bfw} 
\big\} \subset S_{\bfw}
\end{multline}
and
\begin{equation}
\label{D'w}
D'_{\bfw} = \Span B'_{\bfw}.
\end{equation}

The following theorem is our explicit freeness 
result for the module structure of $S_{\bfw}$
over their subalgebras $N'_{\bfw}$.
\thref{free} follows immediately from it.

\bth{f2} Assume that $\KK$ is an arbitrary base field, 
$q \in \KK^*$ is not a root of unity, 
${\bfw} = (w_+, w_-) \in  W \times W$, 
and that ${\vec{w}}_\pm$ are 
reduced expressions of $w_\pm$. 
Then we have the following freeness of 
$S_{\bfw}$ as a left and right $N'_{\bfw}$ module:
\[
S_{\bfw} \cong N'_{\bfw} \bigotimes_\KK D'_{\bfw}
\cong D'_{\bfw} \bigotimes_\KK N'_{\bfw}.
\]
\eth

To prove \thref{f2} it is sufficient to prove the corresponding
statement at the level of associated graded modules, 
which is established by the following result.

\bpr{fag} In the setting of \thref{f2}, the sets
\[
\{\gr y_\la \mid \la \in P^+_{\SS({\bfw})} \}  \odot \gr B'_{\bfw}
\quad
\mbox{and}  
\quad
\gr B'_{\bfw} \odot  \{\gr y_\la \mid \la \in P^+_{\SS({\bfw})} \} 
\]
are bases of $\gr S_{\bfw}$.  In other words
\[
\gr S_{\bfw} \cong \gr N'_{\bfw} \bigotimes_\KK \gr D'_{\bfw}
\cong \gr D'_{\bfw} \bigotimes_\KK \gr N'_{\bfw}
\]
by using the multiplication $\odot$ in $\gr S_{\bfw}$.
\epr
  
\begin{proof} Since the elements $\gr y_\la$ normalize 
each element of $\gr B'_{\bfw}$, it suffices to prove 
that $\{\gr y_\la \mid \la \in P^+_{\SS({\bfw})} \}  \odot \gr B'_{\bfw}$
is a basis of $S_{\bfw}$.

Applying \eqref{sla}, we obtain that for each 
$\la \in P_{\SS({\bfw})}^+$, 
$(\mu_1, \mu_2) \in  \Om_{\bfw}$, 
${\bf{n}}_\pm \in \De(\vec{w}_\pm)$
there exists an integer 
$m_{\la, \mu_1, \mu_2, {\bf{n}}_+, {\bf{n}}_-   }$
such that 
\begin{multline*}
\big( \gr y_\la \big) \odot 
\gr \big[ (\varphi^+_{w_+})^{-1} \big( (X^-)^{{\bf{n}}_+} \big) d^+_{w_+, \mu_1} d^-_{w_-, \mu_2}
(\varphi^-_{w_-})^{-1} \big( (X^+)^{{\bf{n}}_-} \big) \big]
\\
= 
q^{ m_{\mu_1, \mu_2, {\bf{n}}_+, {\bf{n}}_- } }  
(\varphi^+_{w_+})^{-1} \big( (X^-)^{{\bf{n}}_+} \big) d^+_{w_+, \mu_1} 
\big(d^+_{w_+, \la}   d^-_{w_-, \la} \big) 
d^-_{w_-, \mu_2}
(\varphi^-_{w_-})^{-1} \big( (X^+)^{{\bf{n}}_-} \big),
\end{multline*}
where in the right hand side we used the identification \eqref{grS}.
\leref{comp} (ii) now implies that for some integers 
$m_{\la_1, \la_2, {\bf{n}}_+, {\bf{n}}_- }$, 
$(\la_1, \la_2) \in P_{\SS(w_+)}^+ \times P_{\SS(w_-)}^+$, 
${\bf{n}}_\pm \in \De(\vec{w}_\pm)$, 
\begin{align*}
&\{\gr y_\la \mid \la \in P^+_{\SS({\bfw})} \}  \odot \gr B'_{\bfw}
\\
=&\big\{ 
q^{m_{\la, \la_1, \la_2, {\bf{n}}_+, {\bf{n}}_- }}
(\varphi^+_{w_+})^{-1} \big( (X^-)^{{\bf{n}}_+} \big) d^+_{w_+, \la_1} d^-_{w_-, \la_2}
(\varphi^-_{w_-})^{-1} \big( (X^+)^{{\bf{n}}_-} \big) \; \big| \; 
\\
&(\la_1, \la_2) \in P_{\SS(w_+)}^+ \times P_{\SS(w_-)}^+,
{\bf{n}}_\pm \in \De(\vec{w}_\pm)
\big\},
\end{align*}
where again the right hand side uses the identification \eqref{grS}.
By \thref{freeS} (ii)
\[
\big\{ (\varphi^+_{w_+})^{-1} \big( (X^-)^{{\bf{n}}_+} \big) d^+_{w_+, \la_1} \,
\big| \, {\bf{n}}_+ \in \De( \vec{w}_+), 
\la_1 \in P_{\SS(w_+)}^+ 
\big\}
\]
and 
\[
\big\{ 
d^-_{w_-, \la_2} (\varphi^-_{w_-})^{-1} \big( (X^+)^{{\bf{n}}_-} \big) \,
\big| \, {\bf{n}}_- \in \De( \vec{w}_-), 
\la_2 \in P_{\SS(w_-)}^+ 
\big\}
\]
are bases of $S^+_{w_+}$ and $S^-_{w_-}$ respectively. 
Since $S_{\bfw} \cong S^+_{w_+} \otimes_\KK S^-_{w_-}$ 
under the multiplication map, we obtain that 
$\{\gr y_\la \mid \la \in P^+_{\SS({\bfw})} \}  \odot \gr B'_{\bfw}$
is a basis of $\gr S_{\bfw}$.
\end{proof}
\subsection{Structure of the algebras $L_{\bfw}$ and 
freeness of $R_{\bfw}$ over $L_{\bfw}$}
\label{6.5}
In this subsection we obtain an explicit version of the 
freeness result in \thref{free} and describe the structure 
of the algebra $L_{\bfw}$. 

We begin with some implications of the results from the 
previous subsection to the structure of 
$S_{\bfw}[y^{-1}_{\om_i},  i=1, \ldots, r]$.
Denote by $L'_{\bfw}$ the subalgebra of 
$S_{\bfw}[y^{-1}_{\om_i},  i=1, \ldots, r]$, 
generated by $y^{\pm 1}_{\om_i}$, $i=1, \ldots, r$.
\thref{f2} immediately implies:

\bco{free3} Assume that $\KK$ is an arbitrary base 
field, $q \in \KK^*$ is not a root of unity,  
${\bfw} \in W \times W$, and $\vec{w}_\pm$ 
are reduced expressions of $w_\pm$. Then:

(i) $L'_{\bfw}$ is isomorphic to the quantum torus algebra
over $\KK$ of dimension $|\SS({\bfw})|$
with generators $(y_{\om_i})^{\pm 1}$, $i \in \SS({\bfw})$ and
relations \eqref{y-qcomm}.

(ii) The ring $S_{\bfw}[y^{-1}_{\om_i}, i=1, \ldots, r]$ is a 
free left and right $L'_{\bfw}$-module and more precisely:
\[
S_{\bfw}[y^{-1}_{\om_i}, i=1, \ldots, r]
\cong L'_{\bfw} \bigotimes_\KK D'_{\bfw}
\cong D'_{\bfw} \bigotimes_\KK L'_{\bfw}.
\]
\eco

Recall from \S \ref{6.1} that $L_{\bfw}$ denotes 
the subalgebra of $R_{\bfw}$ which is generated 
by $c^+_{w_+, \la}$, $c^-_{w_-, \la}$, 
$\la \in P$.
Its structure is described in
the following result.

\bpr{Lw} The algebra $L_{\bfw}$ is a  
quantum torus algebra over $\KK$
of dimension $r + |\SS({\bfw})|$  with generators 
$(c^+_{w_+, \om_i})^{\pm1}$, $i \in \SS({\bfw})$ and 
 $(c^-_{w_-, \om_j})^{\pm1}$, $j = 1, \ldots, r$, 
and relations
\begin{align*}
c^+_{w_+, \om_{i_1}} c^+_{w_+, \om_{i_2}} &=
c^+_{w_+, \om_{i_2}} c^+_{w_+, \om_{i_1}}, 
\; \; i_1, i_2 \in \SS({\bfw}), 
\\
c^-_{w_-, \om_{j_1}} c^-_{w_-, \om_{j_2}} &=
c^-_{w_-, \om_{j_2}} c^-_{w_-, \om_{j_1}}, 
\; \; j_1, j_2 =1, \ldots, r,
\\
c^+_{w_+, \om_i} c^-_{w_-, \om_j} &=
q^{ - \lcor w_+ \om_i, w_- \om_j \rcor }  
c^-_{w_+, \om_j} c^+_{w_+, \om_i}, 
\; \; i \in \SS({\bfw}), j = 1, \ldots, r.
\end{align*} 
\epr
\begin{proof} The inverse of the isomorphism 
\eqref{ySR} restricts to an algebra isomorphism 
\[
\psi^{-1}_{\bfw} \colon 
L_{\bfw} \to L'_{\bfw} \# \wh{L}^-_{w_-},
\]
where 
\begin{align}
\label{psi1}
\psi^{-1}_{\bfw} ( c^-_{w_-, \om_j}) &= 
c^-_{w_-, \om_j}, \; \; j=1, \ldots, r,
\\
\label{psi2}
\psi^{-1}_{\bfw} ( c^+_{w_+, \om_i}) &= 
q^{ - \lcor w_+ \om_i, w_- \om_i \rcor +1 }
(y_{\om_i})^{-1} \# (c^-_{w_-, \om_j})^{-1}, 
\; \; i=1, \ldots, r.
\end{align}
The second equality follows from \eqref{y} and 
\eqref{n3}.
By \coref{free3}, $L'_{\bfw}$ is a quantum torus algebra 
over $\KK$ of dimension $\SS({\bfw})$
with generators $(y_{\om_i})^{\pm1}$, $i \in \SS({\bfw})$
and by \eqref{y-1}, $y_{\om_i}=1$ for all 
$i \in \II({\bfw})$. Recall from 
\S \ref{3.4} that $\wh{L}^-_{w_-}$ 
is an $r$ dimensional Laurent 
polynomial algebra over $\KK$ with generators
$(c^-_{w_-, \om_j})^{\pm1}$, $j = 1, \ldots, r$.
The commutation relation \eqref{smash} implies that
$L'_{\bfw} \# \wh{L}^-_{w_-}$ is a 
quantum torus algebra over $\KK$ 
of dimension $r + |\SS({\bfw})|$ 
with generators  $(y_{\om_i} \# 1)^{\pm1}$, $i \in \SS({\bfw})$ and
$(c^-_{w_-, \om_j})^{\pm1}$, $j=1, \ldots, r$. 
Therefore $L'_{\bfw} \# \wh{L}^-_{w_-}$ is also a 
quantum torus algebra with generators 
$((y_{\om_i})^{-1} \# (c^-_{w_-, \om_j})^{-1})^{\pm1}$, 
$i \in \SS({\bfw})$ and $(c^-_{w_-, \om_j})^{\pm1}$, $j=1, \ldots, r$. 
It follows from \eqref{psi1}--\eqref{psi2} that the
algebra $L_{\bfw}$ is isomorphic to the 
quantum torus algebra over $\KK$ 
of dimension $r + |\SS({\bfw})|$ with generators 
$(c^+_{w_+, \om_i})^{\pm1}$, $i \in \SS({\bfw})$ and 
$(c^-_{w_-, \om_j})^{\pm1}$, $j = 1, \ldots, r$. The 
commutation relations between them are derived from 
\eqref{normalization} and \eqref{n3}.
\end{proof}

\bco{basis} For all ${\bfw}=(w_+, w_-) \in W \times W$,
the algebra $L_{\bfw}$ has a basis consisting of 
\begin{equation}
\label{ccb}
c^+_{w_+, \mu_1} c^-_{w_-, \mu_2}
\end{equation}
for $\mu_1 \in P_{\SS({\bfw})}$, $\mu_2 \in P$.
We have
\begin{equation}
\label{cc}
c^+_{w_+, \mu} c^-_{w_-, \mu} \in \KK^*, \; \; 
\forall \mu \in P_{\II({\bfw})}.
\end{equation}
In particular, the set \eqref{ccb} for 
$\mu_1 \in P$, $\mu_2 \in P_{\SS({\bfw})}$
is also a basis of $L_{\bfw}$.
\eco
\begin{proof} The corollary follows from \prref{Lw} 
and \eqref{mult}.
\end{proof}

Recall \eqref{ySR}, \eqref{B'w} and \eqref{D'w}, and denote
\[
B_{\bfw} = B'_{\bfw} \# 1 =
\{ b \# 1 \mid b \in B'_{\bfw} \} 
\subset
S_{\bfw}[y_{\om_1}^{-1}, \ldots, y_{\om_r}^{-1}]
\# \wh{L}^-_{w_-}, \quad
D_{\bfw} = \Span B_{\bfw}.
\]
The next theorem provides an explicit form of the freeness result 
from \thref{free}. 

\bth{f1} Assume that $\KK$ is an arbitrary base 
field, $q \in \KK^*$ is not a root of unity,  
${\bfw} = (w_+, w_-) \in W \times W$, and $\vec{w}_\pm$ 
are reduced expressions of $w_\pm$. Then
the algebra $R_{\bfw}$ is a free left and right 
$L_{\bfw}$-module via
\[
R_{\bfw} 
\cong L_{\bfw} \bigotimes_\KK (\psi_{\bfw})^{-1} (D_{\bfw})
\cong (\psi_{\bfw})^{-1} (D_{\bfw}) \bigotimes_\KK L_{\bfw}.
\]
\eth
\begin{proof}
The isomorphism \eqref{ySR}, \coref{free3} (ii)
and the fact that 
$\psi_{\bfw}$ restricts to an algebra isomorphism
$L'_{\bfw} \# \wh{L}^-_{w_-} \to L_{\bfw}$ 
imply 
\[
R_{\bfw} 
\cong (\psi_{\bfw})^{-1} (D_{\bfw}) \bigotimes_\KK L_{\bfw}.
\]
The equality 
\[
L_{\bfw} \bigotimes_\KK (\psi_{\bfw})^{-1} (D_{\bfw})
\cong (\psi_{\bfw})^{-1} (D_{\bfw}) \bigotimes_\KK L_{\bfw}
\]
follows from \eqref{smash}.
\end{proof}
\sectionnew{A classification of maximal ideals of $R_q[G]$
and a question of Goodearl and Zhang}
\lb{Max}
\subsection{A projection property of the ideal $I_{(1,1)}$}
\label{7.1}
In this section we classify all maximal ideals of $R_q[G]$ 
and derive an explicit formula for each of them. We apply this 
result to resolve a question of Goodearl and 
Zhang \cite{GZ} by showing that all maximal ideals of $R_q[G]$ have finite
codimension. In the next section we use this result to prove 
that $R_q[G]$ has the property that all maximal 
chains of prime ideals of it have the same length.
The main step in the proof of the classification 
theorem is to prove that only the highest stratum of the 
decomposition of $\Spec R_q[G]$ in \thref{J-thm}
contains maximal ideals. To obtain this, we combine 
the methods from Section \ref{Dixmier} giving formulas for 
the primitive ideals of $R_q[G]$ and the freeness
theorem for the algebras $R_{\bfw}$ from Section \ref{Module}.
We analyze the images of the primitive ideals and the ideal 
$I_{(1,1)}$ in the direct sum decomposition from \thref{f1}
and eventually deduce that none of the primitive ideals 
in $\Prim_{\bfw} R_q[G]$ is maximal for ${\bfw} \neq (1,1)$. 

We start with the statement of the key step of the classification 
result:

\bth{max} For an arbitrary base field $\KK$ and $q \in \KK^*$ 
which is not a root of unity, all maximal ideals of $R_q[G]$ 
belong to $\Spec_{(1,1)} R_q[G]$, i.e.
\[
\Max R_q[G] \subset \Spec_{(1,1)} R_q[G].
\] 
\eth

In the setting of \S \ref{6.4}--\ref{6.5} denote 
\[
B^\ci_{\bfw} = B_{\bfw} \backslash \{1\}
\]
and
\[
R^\ci_{\bfw} = (\psi_{\bfw})^{-1} (\Span B^\ci_{\bfw}) \bigotimes_\KK L_{\bfw}
=  L_{\bfw} \bigotimes_\KK (\psi_{\bfw})^{-1} (\Span B^\ci_{\bfw}).
\]
By \thref{f2} we have the direct sum decomposition of $L_{\bfw}$-bimodules
\[
R_{\bfw} = L_{\bfw} \oplus R^\ci_{\bfw}.
\]
We denote by 
\[
\pi_{\bfw} \colon R_{\bfw} \to L_{\bfw}
\]
the projection onto the first component (which is a homomorphism 
of $L_{\bfw}$-bimodul\-es). Denote by $N_{\bfw}$ the subalgebra of 
$L_{\bfw}$, generated by $c^+_{w_+, \la}$ and $c^-_{w_-, \la}$ for 
$\la \in P^+$. We will need two $N_{\bfw}$-(bi)submodules 
of $L_{\bfw}$
\[
M^{++}_{\bfw} \subset M^{+}_{\bfw} \subset L_{\bfw}, 
\] 
defined as follows. Denote the 
submonoids
\begin{align*}
Y^{++}_1 &= \{ \mu \in P \mid (1-w_+)\mu > 0 \},
\\
Y^{++}_2 &= \{ \mu \in P_{\SS({\bfw})} 
\mid (1-w_-)\mu > 0 \},
\end{align*}
where the inequalities are in terms of the partial order \eqref{po}. 
Denote also the following two submonoids of $P \times  P_{\SS({\bfw})}$:
\begin{equation}
\label{Y==}
Y_{\bfw}^{++} = 
\big( Y^{++}_1 \times Y^{++}_2  \big) 
\bigsqcup 
\big( Y^{++}_1 \times   P_{\SS({\bfw}) \cap \II(w_-)} \big)
\bigsqcup   
\big( P_{\II(w_+)} \times Y^{++}_2  \big)
\end{equation}
and
\begin{equation}
\label{Y=}
Y_{\bfw}^+ =  Y_{\bfw}^{++} \bigsqcup
\big( 
P_{\II(w_+)} \times P_{\SS({\bfw}) \cap \II(w_-)} 
\big).
\end{equation}
We have disjoint unions in \eqref{Y==}-\eqref{Y=}, 
because
\begin{align*}
\II(w_+)&= \{ i =1, \ldots, r \mid (1-w_+) \om_i = 0 \},
\\
\SS({\bfw}) \cap \II(w_-) &= \SS(w_+) \cap \II(w_-) 
= \{ i \in \SS({\bfw}) \mid (1-w_-) \om_i = 0 \}.
\end{align*}

\bre{sets} Note that in general $Y_{\bfw}^{++}$ and  $Y_{\bfw}^+$
are strictly contained in the sets
\begin{multline*}
\{ (\mu_1, \mu_2) \in P \times P_{\SS({\bfw})} 
\mid (1-w_+) \mu_1 \geq 0, (1-w_-)\mu_2 \geq 0 
\\
\mbox{and at lest one inequality is strict} \; \}
\end{multline*}
and
\[
\{ (\mu_1, \mu_2) \in P \times P_{\SS({\bfw})} 
\mid (1-w_+) \mu_1 \geq 0, (1-w_-)\mu_2  \geq 0 \}, 
\]
respectively. This is so, 
because $\ker (1- w_\pm) \cap P$ are generally
larger than $P_{\II(w_\pm)}$.
\ere

Let
\begin{align}
\label{set1}
M^{++}_{\bfw} &= \{ c^+_{w_+, \mu_1} c^-_{w_-, \mu_2} \mid 
(\mu_1, \mu_2) \in Y_{\bfw}^{++} \},
\\
\label{set2}
M^+_{\bfw} &= \{ c^+_{w_+, \mu_1} c^-_{w_-, \mu_2} \mid 
(\mu_1, \mu_2) \in Y_{\bfw}^+ \}.
\end{align}
Since $(1- w_\pm) \la >0$ for all 
$\la \in P^+_{\SS(w_\pm)}$, using \eqref{mult} 
we obtain that
$M^{++}_{\bfw} \subset M^+_{\bfw}$ are 
$N_{\bfw}$-(bi)submodules of $L_{\bfw}$.
Although we will not need this below, 
we note that $N_{\bfw} \subset M^+_{\bfw}$,
which follows from \eqref{cc} and the fact that 
$(1-w_\pm) \la > 0$ 
for all $\la \in P^+_{\SS(w_\pm)} \backslash \{ 0 \}$.

\coref{basis} implies:
\ble{basis2} For all ${\bfw} = (w_+, w_-) \in W \times W$:

(i) The algebra $N_{\bfw}$ has a $\KK$-basis consisting of
\[
c^+_{w_+, \mu} c^-_{w_-, \la},
\; \;  \mu \in P^+_{\SS({\bfw})} \oplus P_{\II({\bfw})},
\la \in P^+_{\SS({\bfw})}.
\]

(ii) The spanning sets in \eqref{set1} and \eqref{set2} are $\KK$-bases of the $N_{\bfw}$-modules 
$M^{++}_{\bfw}$ and $M^+_{\bfw}$, respectively.
\ele

The following proposition contains the main tool for the proof of \thref{max}.

\bpr{projj} For an arbitrary base field $\KK$, $q \in \KK^*$ which is not a root
of unity and ${\bfw} \in W \times W$, we have
\begin{equation}
\label{proj1}
\pi_{\bfw} (R_q[G]/I_{\bfw} ) \subset M^+_{\bfw}
\end{equation}
and
\begin{equation}
\label{proj2}
\pi_{\bfw} (I_{(1,1)} /I_{\bfw} ) \subset M^{++}_{\bfw}.
\end{equation}
\epr
\subsection{Proof of \prref{projj}}
\label{7.2}
\thref{f2} implies the direct sum decomposition of $N'_{\bfw}$-(bi)modules
\[
S_{\bfw} = N'_{\bfw} \oplus S^\ci_{\bfw}, \; \;  
\mbox{where} \; \; 
S^\ci_{\bfw} =
N'_{\bfw} \bigotimes_\KK \Span ( B'_{\bfw} \backslash \{1 \} ).
\] 
Denote by 
\[
\pi'_{\bfw} \colon S_{\bfw}  \to  N'_{\bfw} 
\]
the corresponding projection into the first summand. 
Recall the isomorphism \eqref{ySR} and eqs. \eqref{psi1}--\eqref{psi2}.
Clearly we have 
\begin{equation}
\label{equiv}
\psi_{\bfw} \big( (y_\la^{-1} \pi'_{\bfw}(s) ) \# c^-_{w_-, \mu}  \big)
= \pi_{\bfw} \big( 
\psi_{\bfw} (  (y_\la^{-1} s)  \# c^-_{w_-, \mu}) \big)
\end{equation}
for all $s \in S_{\bfw}$, $\la \in P^+$ and $\mu \in P$.

For simplicity of the exposition we will split the proof of \prref{projj} into 
two parts: the proofs of \eqref{proj2} and \eqref{proj1}. First note that 
\begin{multline}
\label{I}
I_{(1,1)} = \Span \{ c^{\la_1}_{\xi_1, \la_1}  c^{- w_0 \la_2}_{\xi_2, - \la_2} 
\mid \la_1, \la_2 \in P^+; \; \xi_1 \in (V(\la_1)^*)_{ - \mu_1}, 
\\
\xi_2 \in (V( - w_0 \la_2)^*)_{ \mu_2};\; \;  
\mu_1, \mu_2 \in P, \; \mu_1 < \la_1 \; \; \mbox{or} \; \; 
\mu_2 < \la_2 \}
\end{multline}
and
\begin{equation}
\label{RqG}
R_q[G] = I_{(1,1)} \oplus 
\Span \{ c^+_{1, \la_1}  c^-_{1, \la_2} 
\mid \la_1, \la_2 \in P^+ \}.
\end{equation}

For two elements $a$ and $b$ of a $\KK$-algebra $R$ we denote
\begin{equation}
\label{approx}
a \approx b,  \; \; \;  \mbox{if} \; \; \; 
a = q^m b \; \; \mbox{for some} \; \; 
m \in \Zset. 
\end{equation}
\\
\noindent
{\em{Proof of \eqref{proj2} in \prref{projj}}}. 
Recall 
that the images of the elements $c^\la_{\xi, v} \in R_q[G]$ in $R_q[G]/I_{\bfw}$ 
are denoted by the same symbols.

Fix $\la_1, \la_2 \in P^+$, 
$\mu_1, \mu_2 \in P$, and 
$\xi_1 \in (V(\la_1)^*)_{ - \mu_1}$, $\xi_2 \in (V( - w_0 \la_2)^*)_{ \mu_2}$.
In view of \eqref{I} we need to prove that
\begin{equation}
\label{incl}
\pi_{\bfw}( c^{\la_1}_{\xi_1, \la_1}  c^{- w_0 \la_2}_{\xi_2, - \la_2}  ) 
\in M^{++}_{\bfw}
\end{equation}
in the following three cases: case (1) $\mu_1 < \la_1$ and  $\mu_2 < \la_2$; 
case (2) $\mu_1 < \la_1$ and  $\mu_2 = \la_2$;
case (3) $\mu_1 = \la_1$ and  $\mu_2 > \la_2$. We will prove \eqref{incl} in 
cases (1) and (2). Case (3) is analogous to (2) and is left to the reader.

Recall the definition \eqref{g} of the projections 
$g^\pm_{w_\pm}$.
Denote for brevity
\begin{equation}
\label{cb}
c = c^{\la_1}_{\xi_1, \la_1}  c^{- w_0 \la_2}_{\xi_2, - \la_2}.
\end{equation}

Using the identification \eqref{ident} and eq. \eqref{n3} we obtain
\[
c \approx
c^+_{w_+, \la_1}  c^-_{w_-, \la_2} . 
\big(  (c^+_{w_+, \la_1})^{-1} g^+_{w_+} (\xi_1) \big) 
\big(  (c^-_{w_-, \la_2})^{-1} g^-_{w_-} (\xi_2) \big), 
\]
cf. \eqref{approx}. It follows form
\eqref{psi1}--\eqref{psi2} and \eqref{mult} that
\begin{equation}
\label{psic}
(\psi_{\bfw})^{-1} (c) \approx
\big[ y_{\la_1}^{-1} \big(  (c^+_{w_+, \la_1})^{-1} g^+_{w_+} (\xi_1) \big) 
\big(  (c^-_{w_-, \la_2})^{-1} g^-_{w_-} (\xi_2) \big) \big] \# c^-_{w_-, \la_2 - \la_1}. 
\end{equation}

{\em{Case (1):}} Since $\mu_1 < \la_1$ and  $\mu_2 < \la_2$ 
we have 
\begin{multline*}
\big(  (c^+_{w_+, \la_1})^{-1} g^+_{w_+} (\xi_1) \big) 
\big(  (c^-_{w_-, \la_2})^{-1} g^-_{w_-} (\xi_2) \big)  
\in (S_{\bfw})^{\mu_1 - w_+ \la_1, \mu_2 - w_- \la_2}
\\
\subset 
(S_{\bfw})^{ \prec ((1- w_+)\la_1, (1 - w_-) \la_2)}
\end{multline*}
in terms of the partial order $\prec$ on $Q \times Q$
from \S \ref{6.2}. \prref{fag} and \eqref{gry} imply 
\begin{multline*}
\pi'_{\bfw} \big[
\big(  (c^+_{w_+, \la_1})^{-1} g^+_{w_+} (\xi_1) \big) 
\big(  (c^-_{w_-, \la_2})^{-1} g^-_{w_-} (\xi_2) \big)  
\big]
\\
\in \Span \{ y_\la \mid \la \in P^+,
(1-w_+)\la < (1- w_+) \la_1, 
(1-w_-)\la < (1- w_-) \la_2
\}.
\end{multline*}
It follows from \eqref{equiv} that
\begin{multline*}
\pi_{\bfw}(c) \in 
\Span \{ 
\psi_{\bfw}( (y_{\la_1})^{-1} y_\la \#  c^-_{w_-, \la_2 - \la_1})
\mid \la \in P^+,
\\ (1-w_+)(\la_1 - \la) >0,
(1-w_-)(\la_2 - \la) > 0
\}.
\end{multline*}
Taking into account \eqref{psi1}--\eqref{psi2} and \eqref{cc}, we obtain
that $\pi_{\bfw}(c)$ belongs to 
\begin{align*}
&\Span \{ 
c^+_{w_+, \la_1 - \la}  c^-_{w_-, \la_2 - \la}
\mid \la \in P^+,
(1-w_+)(\la_1 - \la) > 0,   (1-w_-)(\la_2 -\la) > 0
\}
\\
& \subseteq
\Span \{ 
c^+_{w_+, \nu_1}  c^-_{w_-, \nu_2}
\mid \nu_1, \nu_2 \in P,
(1-w_+) \nu_1 > 0,   (1-w_-) \nu_2 > 0
\} \\
&= \Span \{ 
c^+_{w_+, \nu_1}  c^-_{w_-, \nu_2}
\mid  (\nu_1, \nu_2) \in 
Y^{++}_1 \times Y^{++}_2 \} 
\subset M^{++}_{\bfw}.
\end{align*}

{\em{Case (2):}} In this case $c^{- w_0 \la_2}_{\xi_2, - \la_2}$
is a scalar multiple of $c^-_{w_-, \la_2}$ and after
rescaling we can assume that 
\[
c = c^{\la_1}_{\xi_1, \la_1} c^-_{1, \la_2},
\] 
cf. \eqref{cb}. It follows from \eqref{psic} that
\[
(\psi_{\bfw})^{-1} (c) \approx
\big[ y_{\la_1}^{-1} \big(  (c^+_{w_+, \la_1})^{-1} g^+_{w_+} (\xi_1) \big) 
d^-_{w_-, \la_2} \big] \# c^-_{w_-, \la_2 - \la_1},
\]
recall \eqref{dd}. Then 
\[
\big(  (c^+_{w_+, \la_1})^{-1} g^+_{w_+} (\xi_1) \big) 
d^-_{w_-, \la_2} \in (S_{\bfw})^{(\mu_1 - w_+ \la_1, (1-w_-) \la_2)}.
\]
\prref{fag}, \eqref{B'w}, \eqref{gry} and the assumption $\mu_1 < \la_1$ imply that
\begin{multline*}
\big(  (c^+_{w_+, \la_1})^{-1} g^+_{w_+} (\xi_1) \big)  d^-_{w_-, \la_2} \in 
\big( S^\ci_{\bfw} \oplus 
\Span \{ y_\la \mid
\la \in P^+, 
\la -\la_2 \in P_{\II(w_-)}, 
\\
(1-w_+) \la < (1-w_+) \la_1 \}
\big) + 
(S_{\bfw})^{ \prec ( (1 - w_+) \la_1, (1-w_-) \la_2)}.
\end{multline*}
Therefore
\begin{align*}
&\pi'_{\bfw} \big[
\big(  (c^+_{w_+, \la_1})^{-1} g^+_{w_+} (\xi_1) \big) 
d^-_{w_-, \la_2}
\big]
\in
\\ 
&\Span \{ y_\la \mid
\la \in P^+, 
(1-w_+) \la < (1-w_+) \la_1, 
\la -\la_2 \in P_{\II(w_-)}
\}
\\
\oplus & \Span \{ y_\la \mid \la \in P^+,
(1-w_+)\la < (1- w_+) \la_1, 
(1-w_-)\la < (1- w_-) \la_2
\}.
\end{align*}
As in case (1), \eqref{equiv} implies that 
$\pi_{\bfw}(c)$ belongs to the span of 
the elements $\psi_{\bfw}( (y_{\la_1})^{-1} y_\la \#  c^-_{w_-, \la_2 - \la_1})$
where $\la \in P^+$ and either
\begin{equation}
\label{11}
(1-w_+) (\la_1 - \la) > 0, 
\la -\la_2 \in P_{\II(w_-)},
\end{equation}
or 
\begin{equation}
\label{22}
(1-w_+)(\la_1 - \la) >0,
(1-w_-)(\la_2 - \la) > 0.
\end{equation}
Eqs. \eqref{psi1}, \eqref{psi2}, and \eqref{cc} 
imply that the span of these elements is the 
space
\begin{align*}
&\Span \{ 
c^+_{w_+, \la_1 - \la}  c^-_{w_-, \la_2 - \la}
\mid \la \in P^+ \; \;
\mbox{satisfies either} \; \; 
\eqref{11} \; \; \mbox{or} \; \; \eqref{22} \}
\\
\subseteq & 
\Span \{ 
c^+_{w_+, \nu_1}  c^-_{w_-, \nu_2}
\mid  (\nu_1, \nu_2) \in 
(Y^{++}_1 \times P_{ \SS({\bfw}) \cap \II(w_-) }) \sqcup
(Y^{++}_1 \times Y^{++}_2), 
\}
\subset M^{++}_{\bfw}
\end{align*}
which completes the proof of \eqref{proj1}.
\qed
\\

Recall the notation \eqref{la+-0}. We have the 
decomposition
\[
P = P_{\SS({\bfw})} \oplus P_{\II({\bfw})}.
\]
For $\la \in P$ denote its components
\begin{equation}
\label{olo}
\ol{\la} \in P_{\SS({\bfw})}, 
\quad
\ol{\ol{\la}} \in P_{\II({\bfw})}
\end{equation}
with respect to this decomposition.
Denote the delta function on $P^+$: for 
$\la_1, \la_2 \in P^+$ 
\[
\de_{\la_1, \la_2} = 1, \; \; 
\mbox{if} \; \; \la_1 = \la_2
\quad
\mbox{and}
\quad
\de_{\la_1, \la_2} = 0, \; \; \mbox{otherwise}.
\]
\\
Eq. \eqref{proj1} in \prref{projj} follows 
from \eqref{proj2} and the following lemma.

\ble{zero} For all $\la_1, \la_2 \in P^+$,
\begin{multline*}
\pi_{\bfw} (c^+_{1, \la_1} c^-_{1, \la_2}) \in 
\Span 
\{ 
c^+_{w_+, \mu_1} c^-_{w_-, \mu_2} \mid 
(\mu_1, \mu_2) \in 
\\
(P_{\II(w_+)} \times P_{\SS({\bfw}) \cap \II(w_-)})
\sqcup ( Y^{++}_1 \times Y^{++}_2 )
\} \subset M^+_{\bfw},
\end{multline*}
cf. \eqref{Y=} and \eqref{set2}.
\ele
\begin{proof} Set 
\[
c = c^+_{1, \la_1} c^-_{1, \la_2}.
\]
Applying \eqref{psic}, \eqref{y-1}, \eqref{ti} and \eqref{ti2}, 
we obtain
\begin{align*}
(\psi_{\bfw})^{-1} (c) 
& \approx
\big( y_{\la_1}^{-1} d^+_{w_+, \la_1} d^-_{w_-, \la_2} \big) 
\# c^-_{w_-, \la_2 - \la_1} 
\\
& \approx
\big( y_{ \ol{\la_1} }^{-1} 
d^+_{w_+, (\ol{\la_1})_0 + (\ol{\la_1})_+ } 
d^-_{w_-, (\ol{\la_2})_0 + (\ol{\la_2})_- }
\big)
\# c^-_{w_-, \la_2 - \la_1}.
\end{align*}
\prref{fag}, and eqs. \eqref{gry} and \eqref{B'w} imply that
\[
d^+_{w_+, (\ol{\la_1})_0 + (\ol{\la_1})_+ } 
d^-_{w_-, (\ol{\la_2})_0 + (\ol{\la_2})_- } 
\in \KK y_{ (\ol{\la_1})_0 + (\ol{\la_1})_+
+ (\ol{\la_2})_-} +
(S_{\bfw})^{\prec ( (1-w_+) \la_1, (1-w_-) \la_2 ) },
\]
if $(\ol{\la_1})_0= (\ol{\la_2})_0$ and
\[ 
d^+_{w_+, (\ol{\la_1})_0 + (\ol{\la_1})_+ } 
d^-_{w_-, (\ol{\la_2})_0 + (\ol{\la_2})_- }
\in S_{\bfw}^\ci +
(S_{\bfw})^{\prec ( (1-w_+) \la_1, (1-w_-) \la_2 ) },
\] 
otherwise. Thus
\begin{align*}
& \pi'_{\bfw}(  
d^+_{w_+, (\ol{\la_1})_0 + (\ol{\la_1})_+ } 
d^-_{w_-, (\ol{\la_2})_0 + (\ol{\la_2})_- }
)
\in
\de_{(\ol{\la_1})_0, (\ol{\la_2})_0}
\KK y_{ (\ol{\la_1})_0 + (\ol{\la_1})_+
+ (\ol{\la_2})_-} \\
\oplus & \Span \{ y_\la \mid \la \in P^+,
(1-w_+)\la < (1- w_+) \la_1, 
(1-w_-)\la < (1- w_-) \la_2
\}.
\end{align*}
As in cases (1) and (2) of the proof of \eqref{proj1}, 
using \eqref{equiv}, \eqref{psi1}--\eqref{psi2} 
and \eqref{mult} we obtain:
\begin{align*}
\pi(c) \in & \de_{(\ol{\la_1})_0, (\ol{\la_2})_0}
\KK c^+_{w_+, (\ol{\la_1})_- -(\ol{\la_2})_-
+  \ol{\ol{\la_1}} } \, 
c^-_{w_-, (\ol{\la_2})_+ -(\ol{\la_2})_+
+  \ol{\ol{\la_2}} }
\\
\oplus 
& \Span 
\{ 
c^+_{w_+, \mu_1} c^-_{w_-, \mu_2} \mid 
(\mu_1, \mu_2) \in 
P_{\II(w_+)} \times P_{\SS({\bfw}) \cap \II(w_-)} \}.
\end{align*}
Since \eqref{cc} implies
\[
c^+_{w_+, (\ol{\la_1})_- -(\ol{\la_2})_-
+  \ol{\ol{\la_1}} } \, 
c^-_{w_-, (\ol{\la_2})_+ -(\ol{\la_2})_+
+  \ol{\ol{\la_2}} }
\approx
c^+_{w_+, (\ol{\la_1})_- -(\ol{\la_2})_-
+  \ol{\ol{\la_1}} - \ol{\ol{\la_2}}} \, 
c^-_{w_-, (\ol{\la_2})_+ -(\ol{\la_2})_+}
\]
and in addition
\[
\big(
(\ol{\la_1})_- -(\ol{\la_2})_- +  \ol{\ol{\la_1}} - \ol{\ol{\la_2}}, 
(\ol{\la_2})_+ -(\ol{\la_2})_+ \big) 
\in P_{\II(w_+)} \times P_{\SS({\bfw}) \cap \II(w_-)},
\]
we obtain the statement of the lemma.
\end{proof}
\subsection{Proof of \thref{max}}
\label{7.3}
We first analyze the projections of the ideals 
$\pi_{\bfw}(J_{ {\bfw}, \zeta, \theta}/I_{\bfw})$ 
for the primitive ideals defined in \S \ref{4.1}.
Our proof of \thref{max} relies on a combination
of this with \prref{projj}.

Fix ${\bfw} =(w_+, w_-) \in W \times W$.
Recall the definitions \eqref{L} and \eqref{Lred} of the lattices 
$\wt{\LL}({\bfw})$ and $\wt{\LL}({\bfw})_{\red}$. Recall from 
\S \ref{4.1} that $\{ \la^{(1)}, \ldots, \la^{(k)} \}$ is a basis
of $\wt{\LL}({\bfw})$. Denote by $J_{ {\bfw}, {\bf{1}}, {\bf{1}}}$
the ideal \eqref{Jw} corresponding to 
$\zeta_j=1$ for $j=1, \ldots, k$ and $\theta_i=1$ for $i \in \II({\bfw})$.   
Denote
\[
J^0_{ {\bfw}, {\bf{1}}, {\bf{1}} } =
\sum_{j=1}^k R_q[G] b_j(1) + 
\sum_{i \in \II({\bfw})} R_q[G](c^+_{w_+, \om_i} - 1)
+ I_{\bfw}. 
\]
Then 
\[
J_{ {\bfw}, {\bf{1}}, {\bf{1}}} = \{ r \in R_q[G] \mid 
cr \in J_{ {\bfw}, {\bf{1}}, {\bf{1}}}
\; \; 
\mbox{for some} \; \; c \in E_{\bfw} \},
\]
cf. \eqref{Ew}.
Recall from \leref{basis2} (ii) that $M^+_{\bfw}$ has a 
basis comprised of the elements in \eqref{set2}.
Denote by $(M^+_{\bfw})_{\bf{1}}$ the subspace of $M^+_{\bfw}$
which consists of those elements 
\[
\sum_{ (\mu_1, \mu_2) \in Y^+_{\bfw} } 
p_{\mu_1, \mu_2} c^+_{w_+, \mu_1} c^-_{w_-, \mu_2}, 
\quad p_{\mu_1, \mu_2} \in \KK,
\]
which have the property that for all 
$(\mu_1, \mu_2) \in Y^+_{\bfw}$
\[
\sum_{\la \in \wt{\LL}({\bfw})} p_{\mu_1 + \la, \mu_2 - \ol{\la} } =0,   
\]
recall \eqref{olo}. The subspace $(M^+_{\bfw})_{\bf{1}}$ is 
an $N_{\bfw}$ sub-bimodule of $M^+_{\bfw}$ by the following 
lemma.

\ble{MNprod} Let $\mu_1, \mu_2, \nu_1, \nu_2 \in P$ 
and $\{p_\la \in \KK \mid \la \in \wt{\LL}({\bfw}) \}$ 
be a collection of scalars of which only finitely many are 
nonzero. Then:
\begin{multline}
\label{M1p}
c^+_{w_+, \nu_1} c^-_{w_-, \nu_2} \Bigg(
\sum_{\la \in \wt{\LL}({\bfw}) } p_\la 
c^+_{w_+, \mu_1 + \la} c^-_{w_-, \mu_2- \ol{\la} } \Bigg) 
\\
=
q^{ \lcor w_- \nu_2 , ( w_+- w_-) \mu_1 \rcor } 
\Bigg(
\sum_{\la \in \wt{\LL}({\bfw}) } p_\la 
c^+_{w_+, \mu_1 + \nu_1 + \la} c^-_{w_-, \mu_2 + \nu_2 - \ol{\la} } \Bigg) 
\end{multline}
and
\begin{multline}
\label{M2p}
\Bigg(
\sum_{\la \in \wt{\LL}({\bfw}) } p_\la 
c^+_{w_+, \mu_1 + \la} c^-_{w_-, \mu_2- \ol{\la} } \Bigg) 
c^+_{w_+, \nu_1} c^-_{w_-, \nu_2} 
\\
=
q^{ \lcor w_+ \nu_1 , ( w_- - w_+) \mu_2 \rcor } 
\Bigg(
\sum_{\la \in \wt{\LL}({\bfw}) } p_\la 
c^+_{w_+, \mu_1 + \nu_1 + \la} c^-_{w_-, \mu_2 + \nu_2 - \ol{\la} } \Bigg). 
\end{multline}
\ele
\begin{proof} Since $\wt{\LL}({\bfw}) \subset \ker (w_+ - w_-)$, 
for $\la \in \wt{\LL}({\bfw})$ we have 
\begin{align*}
&\lcor w_- \nu_2, w_+ (\mu_1 + \la) \rcor - \lcor \nu_2, \mu_1 + \la \rcor
= 
\\
= &  
\lcor w_- \nu_2, w_+ (\mu_1 + \la) \rcor - \lcor w_- \nu_2, w_- ( \mu_1 + \la)  \rcor
= \lcor w_- \nu_2, (w_+ - w_-) \mu_1 \rcor.
\end{align*}
Eq. \eqref{n3} implies
\begin{align*}
&c^-_{w_-, \nu_2} c^+_{w_+, \mu_1 + \la} =
q^{ \lcor w_- \nu_2, w_+ (\mu_1 + \la) \rcor - \lcor \nu_2, \mu_1 + \la \rcor}
c^+_{w_+, \mu_1 + \la} c^-_{w_-, \nu_2} 
\\
=&
q^{\lcor w_- \nu_2, (w_+ - w_-) \mu_1 \rcor }
c^+_{w_+, \mu_1 + \la} c^-_{w_-, \nu_2}. 
\end{align*}
Now \eqref{M1p} follows from \eqref{mult}. 
Eq. \eqref{M2p} is proved in an analogous way using the fact 
that for $\la \in \wt{\LL}({\bfw})$,
$\ol{\la} \in \wt{\LL}({\bfw})_{\red}$, cf. \eqref{Lred}
and \eqref{LredL}.
\end{proof}

The following result relates the image of the 
$\pi_{\bfw}$-projection of 
$J_{ {\bfw}, {\bf{1}}, {\bf{1}} } /I_{\bfw}$ and the 
above defined subspace of $M^+_{\bfw}$.

\bpr{proj} For an arbitrary base field $\KK$, $q \in \KK^*$
which is not a root of unity, and ${\bfw}=(w_+, w_-) \in W \times W$,
we have
\[
\pi_{\bfw} ( J_{ {\bfw}, {\bf{1}}, {\bf{1}} } /I_{\bfw} )
\subseteq (M_{\bfw}^+)_{\bf{1}}. 
\]
\epr 
\begin{proof} First we will prove that
\begin{equation}
\label{J0}
\pi_{\bfw} (J^0_{ {\bfw}, {\bf{1}}, {\bf{1}} } /I_{\bfw})
\subseteq (M_{\bfw}^+)_{\bf{1}}. 
\end{equation}
Using \prref{projj}, we see that to prove \eqref{J0} 
it is sufficient to prove
\begin{equation}
\label{J0a}
M^+_{\bfw} b_j(1) \subseteq (M^+_{\bfw})_{\bf{1}}
\quad
\mbox{and}
\quad
M^+_{\bfw} (c^+_{w_+, \om_i} -1)
\subseteq (M^+_{\bfw})_{\bf{1}},
\end{equation}
for $j=1, \ldots, k$, $i \in \II({\bfw})$. 
From \eqref{M1p} it follows that 
\begin{multline*}
c^+_{w_+, \la_1} c^-_{w_-, \la_2} 
\Big(c^+_{w_+, \la^{(j)}_+} c^-_{w_-, \la^{(j)}_-} - 
c^+_{w_+, \la^{(j)}_-} c^-_{w_-, \la^{(j)}_+} \Big)
\in (M^+_{\bfw})_{\bf{1}}, 
\\
c^+_{w_+, \la_1} c^-_{w_-, \la_2}
(c^+_{w_+, \om_i} -1)
\in (M^+_{\bfw})_{\bf{1}}, \; \; 
\forall
j=1, \ldots, k, i \in \II({\bfw}), 
(\la_1, \la_2) \in Y^{++}_{\bfw},
\end{multline*}
recall \eqref{bzeta}. This proves \eqref{J0a} 
and thus \eqref{J0}.

\leref{MNprod} implies
\[
\big( (c^+_{w_+, \la_1} c^-_{w_-, \la_2})^{-1} (M^+_{\bfw})_{\bf{1}} \big)
\cap M^+_{\bfw} \subseteq (M^+_{\bfw})_{\bf{1}}, 
\quad \forall \la_1, \la_2 \in P^+,
\]
where the intersection in the left hand side 
is taken inside $L_{\bfw}$. Since 
\[
\pi_{\bfw} ( J_{ {\bfw}, {\bf{1}}, {\bf{1}} } /I_{\bfw} )
\subseteq \bigcup_{\la_1, \la_2 \in P^+}
\Big(
\big(
(c^+_{w_+, \la_1} c^-_{w_-, \la_2})^{-1}
\pi_{\bfw} ( J^0_{ {\bfw}, {\bf{1}}, {\bf{1}} } /I_{\bfw} ) \big)
\cap M^+_{\bfw} \Big),
\]
the proposition follows from \eqref{J0}.
\end{proof}

Next, we proceed with the proof of \thref{max}.
\\ \hfill \\
\noindent
{\em{Proof of \thref{max}.}} First we establish the validity of the 
theorem for algebraically closed fields $\KK$. Assume that 
the statement of the theorem is not correct, 
i.e. there exists ${\bfw} \in W \times W$, 
${\bfw} \neq (1,1)$ such that 
\begin{equation}
\label{contr}
\Max R_q[G] \cap \Spec_{\bfw} R_q[G] \neq \emptyset.
\end{equation}
Let $J \in \Spec_{\bfw} R_q[G]$ be a maximal ideal 
of $R_q[G]$.
\thref{J-thm} (iii) implies that there exits 
$t \in \Tset^r$ such that 
\begin{equation}
\label{Jt}
J \subseteq t . J_{ {\bfw}, {\bf{1}}, {\bf{1}} }, 
\end{equation}
where in the right hand side we use the action 
\eqref{Tract}. Since $J$ is a maximal ideal, 
we have an equality in \eqref{Jt}. Then
\[
J_{ {\bfw}, {\bf{1}}, {\bf{1}} } = t^{-1} . J 
\]
is also a maximal ideal of $R_q[G]$ since the $\Tset^r$-action
\eqref{Tract} is by algebra automorphisms. Because
$J_{ {\bfw}, {\bf{1}}, {\bf{1}} } \in \Spec_{\bfw} R_q[G]$
and ${\bfw} \neq (1,1)$, we have 
$I_{(1,1)} \nsubseteqq J_{ {\bfw}, {\bf{1}}, {\bf{1}} }$.
Therefore
\begin{equation}
\label{sumJI}
J_{ {\bfw}, {\bf{1}}, {\bf{1}} } + I_{(1,1)} = R_q[G].
\end{equation}
Thus there exits 
\[
c \in J_{ {\bfw}, {\bf{1}}, {\bf{1}} } \; \; 
\mbox{such that} \; \; 
c -1 \in I_{(1,1)}.
\]
Let
\[
\pi_{\bfw}(c) = \sum_{ (\mu_1, \mu_2) \in Y^+_{\bfw} } 
p_{\mu_1, \mu_2} c^+_{w_+, \mu_1} c^-_{w_-, \mu_2},
\]
for some $p_{\mu_1, \mu_2} \in \KK$. 
Observe that 
\[
(\la, - \ol{\la}) \notin Y^{++}_{\bfw}, \; \; \forall 
\la \in \wt{\LL}({\bfw}).
\]
Indeed, for all $\la \in \wt{\LL}({\bfw})$,
\[
(1- w_+) \la + (1-w_-) ( - \ol{\la} ) = 
(1- w_+) \la + (1-w_-) ( - \la ) = (w_- - w_+) \la =0,
\]
while every pair $(\mu_1, \mu_2) \in Y^{++}_{\bfw}$ 
has the property that 
\[
(1-w_+) \mu_1 \geq 0, \; 
(1-w_-) \mu_2 \geq 0, \; \; 
\mbox{and at least one of the inequalities is strict},
\]
see \reref{sets}. Since $c-1 \in I_{(1,1)}$, applying \prref{projj}
and \coref{basis}, we obtain
\[
p_{0,0} =1 \; \; 
\mbox{and} \; \;  p_{\la, - \ol{\la} } = 0, \; 
\forall \la \in \wt{\LL}({\bfw}), \la \neq 0.
\]
Therefore
\[
\sum_{\la \in \wt{\LL}({\bfw})} p_{\la, - \ol{\la} } = 1,   
\]
which contradicts with 
$c \in J_{ {\bfw}, {\bf{1}}, {\bf{1}} }$, see \prref{proj}.
This completes the proof of the proposition in the case when 
$\KK$ is algebraically closed.

The general case of the theorem is obtained 
by a base change argument.
Now assume that $\KK$ is an arbitrary base field. Denote
by $\ol{\KK}$ its algebraic closure. For a 
$\KK$-algebra $R$, denote 
$R_{\ol{\KK}} = R \otimes_\KK \ol{\KK}$.
The algebra $(R_q[G])_{\ol{\KK}}$ is isomorphic
to the analog of the algebra $R_q[G]$ defined over 
the base field $\ol{\KK}$.
It is well known and easy to verify that 
the counterparts of $I_{\bfw}$ and $Z(R_{\bfw})$ 
for $(R_q[G])_{\ol{\KK}}$ are 
$(I_{\bfw})_{\ol{\KK}}$ and $Z( (R_{\bfw})_{\ol{\KK}} )$.
Denote by 
$\ol{\iota}_{\bfw} \colon Z( (R_{\bfw})_{\ol{\KK}} )
\to \Spec_{\bfw} (R_q[G])_{\ol{\KK}}$ the counterpart of 
$\iota_{\bfw}$.

Let ${\bfw} \in W \times W$, ${\bfw} \neq (1,1)$. 
If $J \in \Spec_{\bfw} R_q[G]$, then by \thref{J-thm} (ii)
\[
J = \iota_{\bfw}(J^0), \; \; 
\mbox{for some} \; \; J^0 \in \Spec Z(R_{\bfw}).
\]
Moreover $(J^0)_{\ol{\KK}}$ is a proper two sided 
ideal of $Z( (R_{\bfw})_{\ol{\KK}} )$. Thus there 
exits a maximal ideal $\ol{J}^{0}$ of 
$Z( (R_{\bfw})_{\ol{\KK}} )$, containing
$(J^0)_{\ol{\KK}}$. By \thref{J-thm} (ii) 
\[
\ol{\iota}_{\bfw}( \ol{J}^{0} ) \in \Spec_{\bfw} (R_q[G])_{\ol{\KK}}
\] 
and by \thref{max} for algebraically closed base fields
\[
\ol{\iota}_{\bfw}( \ol{J}^{0} ) + (I_{(1,1)})_{\ol{\KK}} \subsetneq (R_q[G])_{\ol{\KK}},
\]
because $\ol{\iota}_{\bfw}( \ol{J}^{0} )$ is contained in a 
maximal ideal which is necessarily in the stratum 
$\Spec_{(1,1)} (R_q[G])_{\ol{\KK}}$. 
Therefore
\[
(J + I_{(1,1)}) \otimes_\KK \ol{\KK}
= \iota_{\bfw}(J^0) \otimes_\KK \ol{\KK} + (I_{(1,1)})_{\ol{\KK}}
\subseteq
\ol{\iota}_{\bfw}( \ol{J}^{0} ) + (I_{(1,1)})_{\ol{\KK}} \subsetneq (R_q[G])_{\ol{\KK}},
\]
so
\[
J + I_{(1,1)} \subsetneq R_q[G].
\]
Consequently $J \notin \Max R_q[G]$, since $J + I_{(1,1)}$
is a proper two sided ideal of $R_q[G]$, properly containing 
$J$. The latter holds because all ideals in $\Spec_{\bfw} R_q[G]$ 
do not contain $I_{(1,1)}$. 
We obtain that $\Max R_q[G] \cap \Spec_{\bfw} R_q[G] = \emptyset$
for all ${\bfw} \in W \times W$, ${\bfw} \neq (1,1)$, which 
proves the theorem for general base fields $\KK$.
\qed
\subsection{Classification of $\Max R_q[G]$ and a 
question of Goodearl and Zhang}
\label{7.4}
By \thref{max} each maximal ideal of $R_q[G]$ 
contains the ideal $I_{(1,1)}$. The structure 
of $R_q[G]/I_{(1,1)}$ is easy to describe;
it is isomorphic to an $r$ dimensional Laurent 
polynomial algebra over $\KK$. From this we obtain 
an explicit classification of $\Max R_q[G]$ 
and an explicit formula for all maximal 
ideals of $R_q[G]$. 

The quotient $R_q[G]/I_{(1,1)}$ is spanned 
by elements of the form $c^+_{1, \la_1} c^-_{1, \la_2}$, 
$\la_1, \la_2 \in P^+$. Since 
$\II(1,1) = \{1, \ldots, r\}$, \coref{basis} implies that 
\[
c^+_{1, \om_i} c^-_{1, \om_i} \in \KK^*, 
\quad \forall i=1, \ldots, r. 
\]
In fact
\[
c^+_{1, \om_i} c^-_{1, \om_i} =1,
\quad \forall i=1, \ldots, r,
\]
because $x_{\om_i} = c^+_{1, \om_i} c^-_{1, \om_i}$ 
in $\wh{R}_{(1,1)}$ (cf. \eqref{x}), and under the 
canonical homomorphism $\wh{R}_{(1,1)} \to R_{(1,1)}$,
$x_{\om_i} \mt 1$ (see \S \ref{3.4} for details). 
Define the algebra homomorphism
\[
\kappa \colon \KK[x_1^{\pm 1}, \ldots, x_r^{\pm 1}]
\to R_q[G]/I_{(1,1)}
\]  
by
\[
\kappa(x_i) = c^+_{1, \om_i}, \; \; 
\mbox{i.e.} \; \; 
\kappa(x_i^{-1}) = c^-_{1, \om_i}, \; \; 
i=1, \ldots, r.
\]
\ble{kappa} In the above setting, 
the map $\kappa \colon \KK[x_1^{\pm 1}, \ldots, x_r^{\pm 1}]
\to R_q[G]/I_{(1,1)}$ is an algebra isomorphism.
\ele
\begin{proof} From the above discussion we have that 
$\kappa$ is surjective. It is injective
by \prref{Lw}, which for ${\bfw}=(1,1)$ states that 
$L_{(1,1)}$ is an $r$ dimensional Laurent polynomial algebra 
over $\KK$ 
with generators $c^+_{w_+, \om_i}$, $i=1, \ldots, r$.
\end{proof}

Although we will not need this below, we note that 
the above arguments establish that
\[
L_{(1,1)} = N_{(1,1)} = R_q[G]/I_{(1,1)}.
\]

Denote the canonical projection
\[
\De_{(1,1)} \colon R_q[G] \to R_q[G]/I_{(1,1)}.
\]
The following result describes explicitly the maximal 
spectrum of $R_q[G]$ and provides an explicit formula
for each maximal ideal.
\bth{mid} Assume that $\KK$ is an arbitrary base field, 
$q \in \KK^*$ is not a root of unity. Then for each quantum group
$R_q[G]$ we have the homeomorphism 
\[
\Max R_q[G] \cong \Max \KK[x_1^{\pm 1}, \ldots, x_r^{\pm 1}],
\]
where both spaces are equipped with the 
corresponding Zariski topologies. Moreover the maximal ideal 
of $R_q[G]$ corresponding 
to $J' \in \Max \KK[x_1^{\pm 1}, \ldots, x_r^{\pm 1}]$ is 
\[
\De_{(1,1)}^{-1}( \kappa(J') ). 
\]
\eth
\begin{proof} The isomorphism $\kappa$ induces a 
homeomorphism
\[
\kappa \colon \Max \KK[x_1^{\pm 1}, \ldots, x_r^{\pm 1}]
\congto \Max (R_q[G]/I_{(1,1)}).
\]
The theorem now follows from the fact that each maximal
ideal of $R_q[G]$ contains $I_{(1,1)}$ by \thref{max}.
\end{proof}

The statement of \thref{mid} is even more explicit in the case 
of algebraically closed base fields $\KK$.

\bco{mid1} If the base field $\KK$ is algebraically closed and
$q \in \KK^*$ is not a root of unity, then each maximal ideal 
of $R_q[G]$ has the form
\begin{equation}
\label{Max-ID}
I_{(1,1)} + (c^+_{1, \om_1} - p_1) R_q[G]
+ \ldots + (c^+_{1, \om_r} - p_r) R_q[G],
\end{equation}
for some $p_1, \ldots, p_r \in \KK^*$.
\eco

Note that in \eqref{Max-ID} one only needs to multiply 
the terms $(c^+_{1, \om_i} - p_i)$, $i=1, \ldots, r$ 
by polynomials in $c^+_{1, \om_i}$, $c^-_{1, \om_i}$,
$i=1, \ldots, r$, because the rest is absorbed by $I_{(1,1)}$, 
see \leref{kappa}.

Another consequence of \thref{max} is the following result.

\bco{2} For an arbitrary base field $\KK$, 
$q \in \KK^*$ which is not a root of unity 
and for each quantum group $R_q[G]$ we have
\[
\Max R_q[G] = \Prim_{(1,1)} R_q[G].
\]
\eco

Finally using \thref{mid}, we settle a question 
of Goodearl and Zhang \cite{GZ}. 
The next section contains a detailed discussion
of the implications of this question.

\bco{GZq} Assume that $\KK$ is an arbitrary 
base field and $q \in \KK^*$ is not a root 
of unity. Then all maximal ideals of the 
quantum function algebras $R_q[G]$ 
have finite codimension. If the base field $\KK$ is algebraically closed, 
then all maximal ideals of $R_q[G]$ have codimension one.
\eco
\begin{proof}
By \thref{mid}, if $J \in \Max R_q[G]$, then there 
exists 
\[
J' \in \Max \KK[x_1^{\pm 1}, \ldots, x_r^{ \pm 1}]
\; \; \mbox{such that} \; \;  
J = \De_{(1,1)}^{-1} \kappa (J'). 
\]
Clearly
\[
R_q[G]/J \cong \KK[x_1^{\pm 1}, \ldots, x_r^{ \pm 1}]/J'.
\]
Since the latter is a quotient of a commutative algebra 
by a maximal ideal, it is finite dimensional. Thus 
$\dim R_q[G] /J  < \infty$.
\end{proof}
\sectionnew{Chain properties and homological applications}
\lb{Homology}
\subsection{Applications}
\label{8.1}
This section contains applications of the results from the 
previous sections to chain properties of ideals 
and homological properties of $R_q[G]$.

We start by recalling two theorems of Goodearl and Zhang, 
and Lu, Wu and Zhang.

\bth{GZ-thm} (Goodearl--Zhang \cite[Theorem 0.1]{GZ}) Assume that 
$A$ is a Hopf algebra over a field $\KK$ which satisfies the 
following three conditions:

(H1) $A$ is noetherian and $\Spec A$ is normally separated, i.e. 
for each two prime ideals $J_1 \subseteq J_2$, there exists a 
nonzero normal element in the ideal $J_2/J_1 \subseteq A/J_1$.

(H2) $A$ has an exhaustive ascending $\Nset$-filtration 
such that the associated graded algebra $\gr A$ is connected 
graded noetherian with enough normal elements, i.e. 
every simple graded prime factor algebra of $\gr A$ contains 
a homogeneous normal element of positive degree.

(H3) Every maximal ideal of $A$ is of finite codimension 
over $\KK$.   

Then $A$ is Auslander--Gorenstein and Cohen--Macauley, and 
has a quasi-Frobe\-nius classical quotient ring. Furthermore, 
$\Spec A$ is catenary and Tauvel's height formula holds.
\eth

\bth{LWZ-thm} (Lu--Wu--Zhang \cite[Theorem 0.4]{LWZ}) 
Assume that $A$ is a noetherian Hopf algebra which satisfies the 
condition (H2). Then $A$ is Auslander--Gorenstein and Cohen--Macauley.
\eth

Lu, Wu and Zhang also proved several other properties 
of noetherian Hopf algebras with the property (H2). We 
refer the reader to \cite[Theorem 0.4]{LWZ} for details.

We recall that a ring $R$ is {\em{catenary}} if each two chains 
of prime ideals between two prime ideals of $R$ have the 
same length. {\em{Tauvel's height formula}} holds for 
$R$, if for each prime ideal $I$ of $R$ its height is equal to 
\[
\GKdim R - \GKdim (R/I).
\]
Recall that a ring $R$ satisfies the first chain condition 
for prime ideals if all maximal chains in $\Spec R$
have the same length equal to $\GKdim R$. 
Such a ring is necessarily catenary. this notion was 
introduced by Nagata \cite{N} in the commutative 
case. We refer the reader to Ratliff's book \cite{R} 
for an exposition of chain conditions for prime ideals.

\bco{fcc} If $A$ is a Hopf algebra over the field $\KK$ which 
satisfies the conditions (H1),(H2) and (H3), then $A$ 
satisfies the first chain condition for prime ideals.
\eco
\begin{proof} If $I \in \Max A$, then its height is 
equal to 
\[
\GKdim R - \GKdim (R/I) = \GKdim R.
\]
Here $\GKdim (R/I)=0$ since $I$ has finite codimension.
\end{proof}

We now turn to applications to the quantum function algebras
$R_q[G]$ and their Hopf algebra quotients.
Among the three conditions (H1), (H2) and (H3) 
for $R_q[G]$, the third turned out to be the hardest to 
prove. It was shown by Hodges and Levasseur \cite{HL0} for 
$\g = {\mathfrak{sl}}_2$, and Goodearl and Lenagan \cite{GLen0}
for $\g = {\mathfrak{sl}}_3$, but was unknown 
for any other simple Lie algebra $\g$. \coref{GZq}
establishes the validity of (H3) for $R_q[G]$ 
in full generality.

Regarding condition (H1) for $R_q[G]$, Joseph proved \cite{J1,J}
that $R_q[G]$ is noetherian, 
and Brown and Goodearl \cite{BG0} proved that $\Spec R_q[G]$ is normally 
separated. These facts are true for an arbitrary 
base field $\KK$ and $q \in \KK^*$ which is not 
a root of unity. In this generality, the noetherianity 
of $R_q[G]$ was proved by Brown and 
Goodearl \cite{BG2}. The proof of \cite{BG0} 
of the normal separation of $\Spec R_q[G]$ 
works in this generality. We briefly sketch a
proof of this. Consider the action \eqref{Tract} of $\Tset^r$ 
on $R_q[G]$. Then $\Tset^r - \Spec R_q[G]$ 
is $\Tset^r$-normally separated. Indeed, assume that
$I'$ is a $\Tset^r$-stable prime ideal of $R_q[G]$
containing the ideal $I_{\bfw}$ for some 
${\bfw} \in W \times W$ (i.e. $I'$ is equal 
to one of the ideals $I_{{\bfw}'}$ for some 
${\bfw}' \in W \times W$). Then 
$E_{\bfw} \cap I' \neq \emptyset$ by 
\thref{J-thm} (i). Any element of 
$E_{\bfw} \cap I'$ provides the 
$\Tset^r$-normal separation of $I_{\bfw}$ and $I'$. 
The normal separation of $\Spec R_q[G]$
follows from the $\Tset^r$-normal separation 
of $\Tset^r - \Spec R_q[G]$ by \cite[Corollary 4.6]{GK}.

The condition (H2) for $R_q[G]$ was proved by Goodearl and Zhang
\cite{GZ} under the assumptions that $\KK= \Cset$ and $q$ is 
transcendental over $\Qset$. One can show that their proof 
works in the general situation by using the fact that 
all $H$-primes of the De Concini--Kac--Procesi algebras 
are polynormal, proved in \cite{Y3}. Instead, we offer a
new and elementary proof of this in the next subsection. 

Thus $R_q[G]$ satisfies all three conditions (H1), (H2) and (H3).
This implies that any quotient of $R_q[G]$ also satisfies them.
Applying \thref{GZ-thm} we obtain:

\bth{RqGt} Assume that $\KK$ is an arbitrary base field
and $q \in \KK^*$ is not a root of unity. Let $I$ be a 
Hopf algebra ideal of any of the quantized function algebras 
$R_q[G]$. Then $R_q[G]/I$ satisfies the first chain condition 
for prime ideals and Tauvel's height formula holds. In addition, 
it is Auslander--Gorenstein and Cohen--Macauley, and 
has a quasi-Frobenius classical quotient ring.
\eth

The Gelfand--Kirillov dimension of $R_q[G]$ is equal to 
$\dim \g$. The fact that $R_q[G]$ has the property that all 
maximal chains of prime ideals of $R_q[G]$ have the same length 
equal to $\dim \g$ is new for all 
$\g \neq {\mathfrak{sl}}_2, {\mathfrak{sl}}_3$
(the two special cases are in \cite{HL0, GL}
combined with \cite{GZ}).
Previously Goodearl and Lenagan \cite{GL0} proved 
that $R_q[SL_n]$ is catenary and Tauvel's height formula holds.
For $\KK= \Cset$ and $q$ transcendental over $\Qset$,
Goodearl and Zhang proved that $R_q[G]$ is catenary 
and Tauvel's height formula holds. The Auslander--Gorenstein 
and Cohen--Macauley properties of $R_q[SL_n]$ were established 
by Levasseur and Stafford \cite{LeSt}. In the case when 
$\KK= \Cset$ and $q$ is transcendental over $\Qset$, 
those properties of $R_q[G]$ were proved 
by Goodearl and Zhang, and for all 
Hopf algebra quotients $R_q[G]/I$ 
by Lu, Wu and Zhang, based on 
Theorems \ref{tGZ-thm} and \ref{tLWZ-thm},
respectively.
\subsection{$R^+ \circledast R^-$ is an algebra with enough normal 
elements.}
\label{8.2}
The algebra $R^+ \circledast R^-$ has a canonical $\Nset$-filtration
with respect to which it is connected. We prove that 
its augmentation ideal is polynormal from which we deduce that 
$R^+ \circledast R^-$ is an algebra with enough normal elements. 
We use this to give an elementary proof of the fact 
that $R_q[G]$ satisfies the condition (H2) from the previous
subsection, under the assumption that $q \in \KK^*$ is 
not a root of unity and without any restrictions on the 
characteristic of $\KK$.

Recall the definition of $R^+ \circledast R^-$ from 
\S \ref{3.4}. It is a noetherian algebra, see 
\cite[Proposition 9.1.11]{J}. One can also prove this
analogously to \cite[Theorem I.8.18]{BG}. 
For $\la \in P$ denote 
\[
\htt 
(\la) = \lcor \la, 
\al_1\spcheck + \ldots + \al_r\spcheck 
\rcor.
\]
Because of \eqref{Cartan}, 
the algebras $R^\pm$ are connected $\Nset$-graded by imposing
\[
\deg c^{\la_1}_{\xi_1, \la_1} = \htt (\la_1), \; \; \xi_1 \in V(\la_1)^*, \quad 
\deg c^{- w_0 \la_2}_{\xi_2, - \la_2} = \htt (\la_2), \; \; \xi_2 \in V(- w_0 \la_2)^*,
\]
for all $\la_1, \la_2 \in P^+$, recall \eqref{special-c}. It follows from \eqref{commR}
and \eqref{R+-} that
\begin{equation}
\label{degR+-}
\deg ( c^{\la_1}_{\xi_1, \la_1} c^{- w_0 \la_2}_{\xi_2, - \la_2} )= 
\htt (\la_1) + \htt (\la_2), \; \; \la_1, \la_2 \in P^+,
\xi_1 \in V(\la)^*, \xi_2 \in V(- w_0 \la)^*
\end{equation}
makes $R^+ \circledast R^-$ a connected $\Nset$-graded algebra.
Denote by $I_{++}$ its augmentation ideal, spanned by elements 
of positive degree.

Recall that an ideal $J$ of a ring $R$ is called polynormal if it 
has a sequence of generators $c_1, \ldots, c_n$ such that 
$c_i$ is normal modulo the ideal generated by $c_1, \ldots, c_{i-1}$, 
for $i=1, \ldots, n$. 

For $i=1, \ldots, r$ fix a basis $B_i$ of $V(\om_i)^*$ consisting
of weight vectors. Let $B = B_1 \sqcup \ldots \sqcup B_r$
and $\ol{B} = B \times \{+, -\}$, where the second term is the set with 
two elements $+$ and $-$. For $\eta = \{\xi, s \} \in \ol{B}$, 
$\xi \in (V(\om_i)^*)_\mu$, denote $[\eta]_1 = i$, 
$[\eta]_2 = \mu $, $[\eta]_3 = s$ and
\begin{align*}
c(\eta) & = c^{\om_i}_{\xi, \om_i}, \quad \mbox{if} \; \; s = +, 
\\
c(\eta) & = c^{\om_i}_{\xi, - w_0 \om_i}, \quad \mbox{if} \; \; s = -.
\end{align*}
(We recall that for each $i= 1, \ldots, r$, there exists 
$j =1, \ldots, r$ such that $-w_0(\om_i)= \om_j$.)
It is well known that $\{ c(\eta) \mid \eta \in \ol{B} \}$ is 
a generating set of the algebra $R^+ \circledast R^-$. In particular,
this set generates the augmentation ideal $I_{++}$. Fix any 
linear order on $\ol{B}$ with the properties that:
\begin{align}
&\mbox{if} \; \; \eta, \eta' \in B \; \; 
\mbox{and} \; \; [\eta]_1 = [\eta']_1, [\eta]_2 > [\eta']_2, 
[\eta]_3 = [\eta']_3=+, \quad \mbox{then} \; \; 
\eta < \eta',
\label{order1}
\\
&\mbox{if} \; \; \eta, \eta' \in B \; \; 
\mbox{and} \; \; [\eta]_1 = [\eta']_1, [\eta]_2 < [\eta']_2, 
[\eta]_3 = [\eta']_3=-, \quad \mbox{then} \; \; 
\eta < \eta',
\label{order2}
\end{align}
where we use the order \eqref{po} on $P$.
\bth{I++} Assume that $\KK$ is an arbitrary base field and 
$q \in \KK^*$ is not a root of unity. If $\{ \eta_1 < \ldots < \eta_N \}$ 
is a linear order on $\ol{B}$ satisfying \eqref{order1}--\eqref{order2}, 
then 
\[
c(\eta_1), \ldots, c(\eta_N)
\]
is a polynormal generating sequence of $I_{++}$.
\eth
\begin{proof} Fix $\eta_k = (\xi, s)$. We will prove that 
$c(\eta_k)$ is normal modulo the ideal generated 
by $c(\eta_1), \ldots, c(\eta_{k-1})$ in the case $s=+$. 
The case $s=-$ is treated analogously.  

We have $\xi \in V(\om_i)^*_\mu$, where 
$i = [\eta_k]_1$ and $\mu = [\eta_k]_2$.
It follows from
\eqref{order1} that there exists a subset 
$\{ j_1 ,\ldots, j_l \} \subseteq \{1, \ldots, k-1\}$ such that
\[
\eta_{j_m} = (\xi_{j_m}, +), m =1, \ldots, l 
\; \; \mbox{for some} \; \; \xi_{j_m} \in B_i
\] 
with the property that
\[
\{ \xi_{j_1}, \ldots, \xi_{j_l} \} \; \; 
\mbox{is a basis of} \; \; 
\bigoplus_{\mu' \in P, \mu'>\mu} V(\om_i)^*_{\mu'}.
\]
\leref{comm} (i) and \eqref{commR} imply that for all 
$a \in (R^\pm)_{\nu, \la} \subset R^+ \circledast R^-$, 
$\nu, \la \in P$:
\[
a c(\eta_k) = q^{\lcor \om_i, \la \rcor - \lcor \mu, \nu \rcor } c(\eta_k) a \mod
\sum_{m=1}^l c( \eta_{j_m} ) (R^+ \circledast R^-). 
\]   
This completes the proof of the theorem.
\end{proof}

The second part of the following corollary
proves that $R_q[G]$ satisfies the property (H2) 
from the previous subsection.
 
\bco{R+-norm} Assume that $\KK$ is an arbitrary base field
and $q \in \KK^*$ is not a root of unity.

(i) Then the algebra $R^+ \circledast R^-$ is a 
connected $\Nset$-graded noetherian algebra with 
respect to the grading \eqref{degR+-} with 
enough normal elements.

(ii) Consider the induced ascending $\Nset$-filtration on 
$R^+ \circledast R^-$ from the grading \eqref{degR+-}
and the induced $\Nset$-filtration on $R_q[G]$
from the canonical surjective homomorphism 
$R^+ \circledast R^- \to R_q[G]$, recall \S \ref{3.4}.
Then the associated graded
$\gr R_q[G]$ is a connected $\Nset$-graded noetherian
algebra with enough normal elements.
\eco
\begin{proof} (i) Let $J$ be a graded ideal of 
$R^+ \circledast R^-$ of codimension strictly 
greater than 1. Then $I_{++}$ is not contained 
in $J$. Let $c_1, \ldots, c_n$ be a polynormal 
generating sequence of $I_{++}$. 
If $c_k$ is the first element in the sequence 
which has a nonzero image in 
$(R^+ \circledast R^-)/J$, then this image is 
a nonzero normal element of the quotient.
This proves (i).

Part (ii) follows from part (i), because 
$\gr R_q[G]$ is a graded quotient 
of 
$\gr (R^+ \circledast R^-) \cong R^+ \circledast R^-$.
\end{proof}

\end{document}